\documentclass{article}
\usepackage[english]{babel}
\usepackage{amsfonts}
\usepackage{amscd}

\title
{\bf  On all numbers great and small (Topological fields of Conway's numbers and their completions)}

\author{Ju. T. Lisica\thanks {{\it Dedicated to the memory of John Horton Conway}. 2010 Mathematics Subject Classification. Primary 03E65, 03C62, 54G99. Keywords and expressions: Conway's numbers, $\zeta$-completions.}\\
St. Tikhon's Orthodox   University of Humanities\\
Likhov Lane 6, p. 1, Moscow, 127051,\\
e-mail: jutlisica@yandex.ru}

\date{}
\begin{document}
\maketitle

\begin{abstract}
The proper Class $\bf{No}$ of all Conway's numbers $\cite{l3}$ is considered as a region of investigation.  
It turns out to be a total ordered Field (i.e., a field whose domain is a proper Class) and  this totally, or linear ordered Class, containing the real numbers ${\mathbb R}$ and the ordinal numbers {\bf On}.

For any subfield $F$ of $\bf{No}$, i.e., $F$ is a set nor proper class, considered with topology induced by a linear ordering on $F$ a completion $\tilde F$ is constructed; in particular, for $\zeta=\omega^{\omega^\mu}$, $0\leq\mu<\Omega$, and for  a specially defined subfield $F={\mathbb P}_\zeta\subset{\bf No}$ a complete subfield ${\mathbb R}_\zeta\subset{\bf No}$ is defined as $\tilde {\mathbb P}_\zeta$. 

Fundamental (Cauchy) sequences $(x_\alpha)_{0\leq\alpha<\zeta}$ are considered in a subfield $F\subset {\mathbb P}_\zeta\subset{\bf No}$, where $\zeta$ is the smallest ordinal number which does not belong to $F$,   and they are the main instrument in the paper. 

A fragment of Mathematical Analysis in ${\mathbb R}_\zeta$ is given and two of its non-trivial results are presented: every positive number $x\in{\mathbb R}_\zeta$ has a unique $n$-th root in ${\mathbb R}_\zeta$, for each positive integer $n$ and every odd-degree polynomial with coefficients in ${\mathbb R}_\zeta$ has a root in ${\mathbb R}_\zeta$. Hence so-called fundamental theorem of algebra: the ring ${\mathbb R}_\zeta[i]\stackrel{def}{=}{\mathbb C}_\zeta$ of all numbers of the form $x+iy$ ($x,y\in{\mathbb R}_\zeta$), $i^2=-1$, is an algebraically closed field. 

Trigonometry functions $\sin\,x$, $\cos\,x$ and $\tan\,x$ are defined in ${\mathbb R}_\zeta$ and a proof of the fundamental theorem of algebra in ${\mathbb C}_\zeta$ is given.

Conway's problem of a natural definition of a function $y=a^x$ for all numbers $x\in{\bf No}$ and for positive finite numbers $a$ and even for some infinite numbers $a$ are solved. The corresponding inverse function $y=\log_ax$ are also constructed, which turns out to be continuous in topology of ${\mathbb R}_\zeta$. In Appendix a problem of Foundation of Set Theory is considered.
\end{abstract}

\begin{center}
{\bf 1. Necessary knowledge from Conway's numbers theory}
\end{center}

For definiteness we shall work through the paper inside a von Neumann-Bernays-G\"{o}del-type set theory ($NBG$ for short).

We begin with the best description of {\it numbers} {\bf great} and {\bf small} containing all real numbers ${\mathbb R}$ and all ordinal numbers ${\mathbb O}{\mathbb N}$ as well as other magnutudes like $\omega-1$, $\frac{\omega}{2}$, $\sqrt{\omega}$, $\sqrt[3]{\omega}$, ... $\omega^{\frac{1}{\omega}}$, $\omega^{\sqrt 2}$, etc., great numbers and $\frac{1}{\omega-1}$, $\frac{1}{\omega^2}$, $\omega^{-\sqrt{\omega}}$, $\omega^{-\frac{1}{\omega}}$, etc., small numbers, which was inductively given  by Conway in $\cite{l3}$ (1979). 

Conway's construction of numbers was a modification of two well-known ideas. 

One of them was the Mirimanoff's representation $\cite{l12}$ (1917), later in 1923 repeated by von Neumann $\cite{l090}$, of the ordinal numbers: 
$0=\{\}$, $1=\{\{\}\}=\{0\}$, $2=\{\{\},\{\{\}\}\}=\{0,\{0\}\}$, $3=\{\{\},\{\{\}\},\{\{\},\{\{\}\}\}\}=\{0,1,\{0,1\}\}=\{0,1,2\}$, 
..., $\omega=\{0,1,2,...,n,...\}$, $\omega+1=\{0,1,2,...,n,...,\omega\}$, ... and so on. 

Another one  was the construction of real numbers ${\mathbb R}$ $\cite{l4}$ (1888) via Dedekind sections $\{A\,|\,B\}$ of  rational numbers  ${\mathbb Q}$, i.e., $\xi=\{A\,|\,B\}\in{\mathbb R}$ if and only if  $A$ and $B$ are non-empty sets of rationals such that ${\mathbb Q}=A\cup B$, $A\cap B=\emptyset$, $A$ has no greatest number, and for every $a\in A$ and $b\in B$ one has $a<b$.

Conway's definition of a number is the following.

If $L$ and $R$ are any two {\it sets} of numbers, and {\it no member of} $L$ {\it is} $\geq$ {\it any member of} $R$, then there is a number $\{L\,|\,R\}$.

If $x=\{L\,|\,R\}$ is a number, then, for short, $x=\{x^L\,|\,x^R\}$, where   $x^L$ is a typical member of $L$ and $x^R$ is a typical member of $R$, i.e., $L=\{x^L\}$ and $R=\{x^R\}$ are {\bf sets}.

Further,

$x\geq y$ iff (no $x^R\leq y$ and $x\leq$ no $y^L$); 

$x\leq y$ iff $y\geq x$; as well as 

$x=y$ iff ($x\geq y$ and $y\geq x$); 

$x>y$ iff ($x\geq y$ and $y\not\geq x$); 

$x<y$ iff $y>x$.

 Different representations $\{L'\,|\,R'\}$ and $\{L\,|\,R\}$ can define the same number; that is why one must distinguish between the form $\{L\,|\,R\}$ of a number and the number itself. Conway called $L$ the lower class and $R$ the upper class of number $x=\{L\,|\,R\}$.

Operations $+$, $\cdot$ and $/$ on all Conway's numbers $x,y\in{\bf No}$ are the following:

\begin{equation}
\label{f01}
x+y=\{x^L+y,x+y^L\,|\,x^R+y,x+y^R\},
\end{equation}
\begin{equation}
\label{f02}
-x=\{-x^R\,|\,-x^L\},
\end{equation}
\begin{equation}
\label{f03}
xy=\{x^Ly+xy^L-x^Ly^L,x^Ry+xy^R-x^Ry^R\,|\,x^Ly+xy^R-x^Ly^R,x^Ry+xy^L-x^Ry^L\},
\end{equation}
and when $x>0$ the inverse to it, $xy=1$, is given by the following formula:
\begin{equation}
\label{f04}
y=\{0,\frac{1+(x^R-x)y^L}{x^R},\frac{1+(x^L-x)y^R}{x^L}\,|\,\frac{1+(x^L-x)y^L}{x^L},\frac{1+(x^R-x)y^R}{x^R}\}.
\end{equation}

The formula $(\ref{f04})$ is an inductive definition of an inverse number $y=\{y^L\,|\,y^R\}$ to the number $x=\{x^L\,|\,x^R\}$, where naturally $x\not=0$ and $x$ is positive number; when $x<0$, then $y=-y'$, where $y'$ is the inverse to the number $-x$.

${\bf No}$ is a Field.  Note also that for $x\in{\bf No}$  its {\it absolute value}, denoted by  $|x|$,  is equal to $x$ if $x>0$, and is equal to $-x$ if $x<0$, and is equal to $0$ if $x=0$.

{\bf Theorem 1.}
{\it If $x\geq y$ and $y\geq z$, then $x\geq z$.}

{\bf Theorem 2.}
{\it For any  number $x$ one has $x^L<x<x^R$ for all  $x^L$, $x^R$. Also, for two numbers $x$ and $y$ one must have $x\leq y$ or $x\geq y$.}

So, the totally ordered proper Class $\bf{No}$ of all numbers is the following:

 $0=\{\,|\,\}$ (born on day $0$),

$1=\{0\,|\,\}$ and $-1=\{\,|\,0\}$ (born on day $1$),

$2=\{0,1\,|\,\}$, $\frac{1}{2}=\{0\,|\,1\}$, $-\frac{1}{2}=\{-1\,|\,0\}$, $-2=\{\,|\,-1,0,\}$ (born on day $2$),

 ...

$\omega=\{0,1,2,3,...\,|\,\}$, $\pi$, $e$, $\sqrt{2}$, $1+\frac{1}{\omega}$, $1-\frac{1}{\omega}$,
$1/3=\{0,\frac{1}{4},\frac{5}{16},...\,|\,\frac{1}{2},\frac{3}{8},...\}$,... $1/\omega=\{0\,|\,1,\frac{1}{2},\frac{1}{4},\frac{1}{8},...\}$, $-\omega=\{\,|\,0,-1,-2,-3,...\}$ (born on day $\omega$),

$\omega+1=\{0,1,2,3,...\omega\,|\,\}$, $\omega-1=\{0,1,2,3,...\,|\,\omega\}$, $\sqrt{2}+\frac{1}{\omega}$, $\sqrt{2}-\frac{1}{\omega}$, 
... $-\omega-1$ (born on day $\omega+1$),

 ...

$\omega 2=\{0,1,2,3,...\omega,\omega+1,...\,|\,\}$, $\frac{\omega}{2}=\{0,1,2,3,...\,|\,\omega,\omega-1,\omega-2,...\}$, $\frac{2}{\omega}=\{\frac{1}{\omega}\,|\,1,\frac{1}{2},\frac{1}{4},...\}$, $\frac{1}{2\omega}=\{0\,|\,\frac{1}{\omega}\}$,... $\frac{1}{\omega^2}=\{0\,|\,\frac{1}{\omega},\frac{1}{2\omega},\frac{1}{4\omega},...\}$,
...$-\omega 2=\{\,|\,0,-1,-2,=3,...=\omega,=\omega-1,...\}$ (born on day $\omega 2$),
... and so on.

More precisely,
for each ordinal number $\alpha\in{\bf On}$ Conway defined a set $M_\alpha$ of numbers by setting $x=\{x^L\,|\,x^R\}$ in $M_\alpha$ if all the $x^L$ and $x^R$ are in the union of all the $M_\beta$ for $\beta<\alpha$. Then he putted $O_\alpha=\bigcup\limits_{\beta<\alpha}M_\beta$ and $N_\alpha=M_\alpha\setminus O_\alpha$. Then in the terminology of numbers' birth to which he  adhered:

$M_\alpha$ is the set of numbers born on or before $\alpha$ (Made numbers),

$N_\alpha$ is the set of numbers born first on day $\alpha$ (New numbers), and

$O_\alpha$ is the set of numbers born before day $\alpha$ (Old numbers).

Each $x\in N_\alpha$ defines a Dedekind section $\{L|R\}$ of $O_\alpha$, if one sets $L=\{y\in O_\alpha:y<x\}$ and $R=\{y\in O_\alpha:y>x\}$, then $x=\{L\,|\,R\}$. Moreover, $M_\alpha=O_\alpha\cup N_\alpha$ as the union of $O_\alpha$ together with all its sections, in the natural order. Notice also that for each number $x\in M_\alpha$ there are inequalities inequality $-\alpha\leq x\leq\alpha$ and only two numbers $-\alpha$ and $\alpha$ in $M_\alpha$ have forms $\{\,|O_\alpha\}$ and $\{O_\alpha\,|\,\}$, respectively, which dropped the Dedekind restriction on sets $L$ and $R$  to be nonempty.

Every number $x$ is in a unique set $N_\alpha$ (see $\cite{l3}$ p. 30). Taking into account the above decription, we call a {\it birthday form} $\{x^L\,|\,x^R\}$ of $x$ when the birthdays of all $x^L$, $x^R$ are less than $\alpha$.

{\bf Definition 1}. Two positive Conway's numbers $x$ and $y$ are {\it commensurate} if for some positive integer $n$ we have $x<ny$, $y<nx$; a number $x\in{\bf No}$ is {\it small compared to} $y\in{\bf No}$ if $-y<nx<y$ for all integer $n$; a number $x\in{\bf No}$ is {\it great compared to} $y\in{\bf No}$ if $y$ is small compared to $x$.

To be commensurate with is evidently an equivalence relation whose equivalent Classes are {\it convex}, i.e., if $x<z<y$ and $x$ and $y$ are commensurate, then $z$ is commensurate with both. Thus, there is a unique simplest number in each commensurate Class,  i.e., it is the simplest number among others which commensurate with $x$ or, what is the same, the birthday of it precedes the birthdays of others numbers in its equivalent Class, and such numbers Conway called {\it leaders}. They are powers $\omega^x$ for some numbers $x\in{\bf No}$, i.e., $\omega$ to the power of $x$, and are given  formally by the following formula:
\begin{equation}
\label{f314}
\omega^x=\{0,r\omega^{x^L}\,|\,r\omega^{x^R}\},
\end{equation}
  where $x=\{x^L\,|\,x^R\}$ and $r$ denotes a variable ranging over all positive real numbers.

{\bf Theorem 3.} {\it Each positive number $x\in{\bf No}$ is commensurate with some $\omega^y$, $y\in{\bf No}$}

 Conway proved (Theorem 21 in $\cite{l3}$, p. 33) that each number $x\in{\bf No}$ defines a unique expression
\begin{equation}
\label{f0201}
x=\sum\limits_{0\leq\beta<\alpha}\omega^{y_\beta}r_\beta,
\end{equation}
in which $\alpha$ denotes some ordinal, the numbers $r_\beta$ $(0\leq\beta<\alpha)$ are non-zero reals, and the numbers $y_\beta$ form a descending sequence of numbers. Moreover, normal forms for distinct $y$ are distinct, and every form satisfying these conditions occurs.

{\bf Theorem 4.} 
{\it Every positive number $x$ in ${\bf No}$ has a unique positive $n$th root in ${\bf No}$, for each positive integer $n$.}

{\bf Theorem 5.} 
{\it Every odd-degree polynomial $f(x)=x^n+Ax^{n-1}+Bx^{n-2}+...+K$ with coefficients $A,B,...,K$ in ${\bf No}$ has a root in ${\bf No}$.}

Proofs of these theorems  see in $\cite{l3}$, p. 31-33, 40-41.

\begin{center}
{\bf 2. Definition of a $\zeta$-Archimedean subfield ${\mathbb P}_\zeta$ of ${\bf No}$} 
\end{center}

Let ${\bf On}$ denote by $[0,\Omega)$, where $\Omega=\{{\bf No}\,|\,\}$ is a gap (see below) in ${\bf No}$. Later in this article we will deal with ordinals of the following type: $\zeta=\omega^{\omega^\mu}$, where $0\leq\mu<\Omega$, is $\zeta$ is the {\it main ordinal number} in the sense of Jacobstahl $\cite{l15}$. 
 
 Let  ${\mathbb P}_\zeta$ be a localization in zero of the ring $R_\zeta$ of all Conway's numbers generated by the set $O_\zeta\subset{\bf No}$ of all Conway's numbers born before   day $\zeta$. It is known that $R_\zeta$ is the intersection of all subrings of ${\bf No}$ containing $O_\zeta$ or what is the same each element of $R_\zeta$ is a linear combination with integer coefficients in ${\mathbb Z}={\mathbb Z}_\omega$ of finite products of powers of elements in $O_\zeta$ with natural exponents. 
 
 One can see that $R_\zeta$ is a ring without zero divisors and hence ${\mathbb P}_\zeta=S^{-1}R_\zeta$, where $S=R_\zeta\setminus\{0\}$, is actually a field. In fact it is a topological field whose base ${\cal B}$ of topology is the set of all intervals $(a,b)\subset{\mathbb P}_\zeta$, $-\zeta\leq a<b\leq\zeta$.

 {\bf Definition 2.} A subfield $F$ of ${\bf No}$ is called a $\zeta$-{\it field}, if it contains all numbers of $[0,\zeta)$ and does not contain ordinals $\alpha'\in[\zeta,\Omega)$.

 {\bf Definition 3.} A field $F\subset{\bf No}$ is called $\zeta$-{\it Archimedean} if it satisfies the following $\zeta$-Archimedean property: for each pair of its positive numbers $x,y$ there exists an ordinal number $\alpha\in[0,\zeta)$ such that $\alpha\cdot x>y$.

 {\bf Proposition 1}. {\it For any pair of numbers $x,y>0$  in a $\zeta$-field $F$ such that $x<y$ there exists an ordinal number $\alpha\in [0,\zeta)$ such that $x\cdot\alpha>y$ and thus, every $\zeta$-field $F$ is $\zeta$-Archimedian}.
  
{\bf Proof}.  Let $x,y$ be any positive numbers in a $\zeta$-field $F$ such that $x<y$. By Definition 2, $[0,\zeta)$ is co-final in $F$ and thus, for each number $a\in F$, there exists an ordinal $\gamma<\zeta$ such that $a<\gamma$. In particular, for numbers $\frac{1}{x},y\in F$ there are ordinal numbers $\gamma'$ and $\gamma''$ such that $\gamma'>\frac{1}{x}$ and $\gamma''>y$. Consequently, there is  an ordinal number $\beta$ such that $\frac{1}{x}<\beta$ and $y<\beta$. Indeed,  $\beta=\max\{\gamma',\gamma''\}$. Put $\alpha=\beta^2$. Since $\frac{1}{\beta}<x$ we obtain $x\cdot\alpha=\alpha \cdot x>\alpha\cdot\frac{1}{\beta}=\beta^2\cdot\frac{1}{\beta}=\beta>y$. 
 $\Box$

 {\bf Corollary 1.} {\it The field ${\mathbb P}_\zeta$  is a $\zeta$-field and hence satisfies the $\zeta$-Archimedean property.}
 
 {\bf Proof.} We see that $\zeta\notin O_\zeta$, $\zeta\notin P_\zeta$ and  $\zeta\notin {\mathbb P}_\zeta$. Otherwise, $\frac{1}{\zeta}$ is also in ${\mathbb P}_\zeta$ ans well as all nubers $\frac{1}{\alpha}-\frac{1}{\zeta}\in{\mathbb P}_\zeta$ with $\lim\limits_{0<\alpha<\zeta}(\frac{1}{\alpha}-\frac{1}{\zeta})=\frac{1}{\zeta}$. Contradicition with $0<\frac{1}{\zeta}<\frac{1}{\alpha}$, $0<\alpha<\zeta$. Clearly, that $\zeta$ is the smalles ordial such that it not in ${\mathbb P}_\zeta$. Otherwise, if $\zeta'\notin{\mathbb P}_\zeta$ and  $\zeta'<\zeta$, then $\zeta'\in O_\zeta\subset{\mathbb P}_\zeta$. Contradiction. 
 
  Now it is enough to prove  that $\bigcup\limits_{0\leq\alpha<\zeta}\{\alpha\}$ is co-final in ${\mathbb P}_\zeta$. 
 
 If $y$ is a positive number in $R_\zeta\subset {\mathbb P}_\zeta$, then as we have already noticed $y$  a linear combination with integer coefficients in ${\mathbb Z}={\mathbb Z}_\omega$ of finite products of powers of elements in $O_\zeta$ with natural exponents. Consequently, each $i$-summand of this combination is $<$ than some ordinal number $\alpha_i<\zeta$, $1\leq i\leq n$, because $\zeta$ is the main ordinal number in the sense of Jacobstahl and for $\alpha=\max\{\alpha_i\}_{1\leq i\leq n}$ we obtain $y<\alpha$ and $\alpha\in[0,\zeta)$.

 At last, let $y$ be a positive number in ${\mathbb P}_\zeta\setminus R_\zeta$. Then $y=[\frac{u}{z}]$, where $u,z$ are positive numbers in $R_\zeta$. One can easily see that $z\geq 1$. If it is not so, then
we can find an ordinal number $\alpha'<\zeta$ such that $z\cdot\alpha'>1$ because we can prove it like in the proof of Proposition 1, i.e., for numbers $z,1\in R_\zeta$, and take $\frac{u\cdot\alpha'}{z\cdot\alpha'}$ which represents the  number $y=[\frac{u}{z}]$ because evidently $[\frac{u}{z}]=[\frac{u\cdot\alpha'}{z\cdot\alpha'}]$.

 As above, $\frac{v}{z}=\frac{u+u}{z}>\frac{y}{x}$ and hence $v\cdot x=x\cdot v=x\cdot(\frac{v}{z}\cdot z)=(x\cdot\frac{v}{z})\cdot z\geq (x\cdot\frac{v}{z})\cdot 1=x\cdot\frac{v}{z}=x\cdot\frac{v}{z}> x\cdot\frac{y}{x}=y$. 
 
 Thus, ${\mathbb P}_\zeta$ is a $\zeta$-field and, by Proposition 1, satisfied the $\zeta$-Archimedean property. $\Box$
 
 {\bf Corollary 2.} {\it If a subfield $F$ of ${\bf No}$ is a subset of $M_\zeta$, where $\zeta$ is the smallest ordinal which is not in $F$, then $F$ is a $\zeta$-field and hence satisfies the $\zeta$-Archimedean property.}

 {\bf Proof}. It is enough to show that $[0,\zeta)$ is co-final in $F$.
 
 First of all, $\alpha<\zeta$ implies $\alpha\in F$ because $\zeta$ is the smallest ordinal which is not in $F$.
 
  Let $x$ be an arbitrary positive number of $F$. If $x$ is a number born on day $\alpha<\zeta$, then it is obvious that $x\leq\alpha$ and thus, $x<\alpha+1<\zeta$ because $\zeta$ is a limit ordinal.  
 
 If $x$ is a number born on day $\zeta$, then there is an ordinal number $\alpha\in[0,\zeta)$ such that $x<\alpha$. Otherwise, there should be the following inequality $\alpha\leq x$, or more precisely, $\alpha<x$, for all $\alpha\in N_\zeta$. Note that $x$ is not an ordinal number, say $\alpha_0$ because $[0,\zeta)$ has no maximal element. 
 
 Then $x=\{L\,|\,R\}$, where $L\subset O_\zeta$ and $R\subset O_\zeta$, and $(L,R)$ is a Dedekind section on $O_\zeta$ we conclude that $R\not=\emptyset$ ($x$ is not an ordinal)  and there is a number $x'\in R$ such that $x'$ was born on or before some ordinal $\alpha_1<\zeta$. We have already noticed that each number $x'$ born on or before day $\alpha_1$ there is an inequality $x'\leq\alpha_1$. Hence $x<x'\leq\alpha_1<\zeta$ what is in contradiction with inequalities $\alpha<x$, for all $0\leq\alpha<\zeta$. $\Box$

\begin{center}
{\bf 3. Convergence of $\zeta$-infinite sequences and functions defined on subsets of $\zeta$-fields  $F$ of numbers in ${\bf No}$}
\end{center}

{\bf Definition 4}. For a fixed subfield $F$ of numbers in ${\bf No}$, which is a $\zeta$-field, a mapping $x:(0,\zeta)\rightarrow F$ is called a {\it infinite sequence} of type $\zeta$ of Conway's numbers in $F$, or $\zeta$-{\it infinite sequence} in $F$, or shortly a $\zeta$-{\it sequence} in $F$. The range of $x$ is the set $\{x(1),x(2),...,x(\alpha),...\}$, conveniently written $\{x_1,x_2,...,x_\alpha,...\}$ or simply $x_1,x_2,...,x_\alpha,...$ , where $x_\alpha=x(\alpha)$, $0<\alpha<\zeta$, or even shorter $(x_\alpha)_{0<\alpha<\zeta}$. The elements of the range are called {\it terms}.

{\bf Definition 5}.
We say that $\zeta$-sequence $(x_\alpha)_{0<\alpha<\zeta}$ in a $\zeta$-field $F$ of numbers in ${\bf No}$ {\it converges} to $a\in F$, and we write $\lim\limits_{0<\alpha<\zeta}x_\alpha=a$,  if for each positive  number $\varepsilon\in F$ there is an ordinal number $\alpha_0\in(0,\zeta)$ such that $|x_\alpha-a|<\varepsilon$ for all $\alpha_0\leq\alpha<\zeta$.  
In this case we also say that $\zeta$-sequence $(x_\alpha)_{0<\alpha<\zeta}$ is {\it convergent} in  $F$.  

{\bf Definition 6}. A $\zeta$-sequence is {\it almost-stationary}, {\it almost-positive} and {\it almost-negative}, {\it almost-stationary}, {\it strictly increasing}  and {\it strictly decreasing} if for some $a\in  F$ there exists an ordinal number $0\leq\alpha_0<\zeta$ such that $x_\alpha=x_{\alpha_0}=a$, $x_\alpha>0$,   $x_\alpha<0$, $x_{\alpha}<x_{\alpha'}$ and $x_{\alpha}>x_{\alpha'}$  for all $\alpha_0\leq\alpha<\alpha'<\zeta$, respectively.

{\bf Definition 7}.  A $\zeta$-sequence $(x_\alpha)_{0<\alpha<\zeta}$ in a $\zeta$-field $F$ is caled a $\zeta$-{\it infinitely small sequence}  if $\lim\limits_{0<\alpha<\zeta}x_\alpha=0$.

{\bf Definition 8}.  A $\zeta$-sequence $(x_\alpha)_{0<\alpha<\zeta}$ in a $\zeta$-field $F$ is caled a $\zeta$-{\it infinitely great sequence}  if for each positive  number $E\in F$ there is an ordinal number $\alpha_0\in(0,\zeta)$ such that $|x_\alpha|>E$ for all $\alpha_0\leq\alpha<\zeta$. If it is {\it almost positive} or {\it almost negative}, then we will conditionally denote this by $\lim\limits_{0<\alpha<\zeta}x_\alpha=+\infty_\zeta$ and $\lim\limits_{0<\alpha<\zeta}x_\alpha=-\infty_\zeta$, respectively.

{\bf Remark 1.} Notice that $+\infty_\zeta$ and $-\infty_\zeta$, $\zeta=\omega^{\omega^\mu}$, $0\leq\mu<\Omega$,  are not numbers but formal {\it boundary symbols}, i.e., there is no algebraic operations $+$ and $\cdot$ on them, but formally we can write and understand inequalities $x<+\infty_\zeta$, $x'>+\infty_\zeta$ or $x>-\infty_\zeta$, $x<-\infty_\zeta$, for some (not all) Conway's numbers $x,x'\in {\bf No}$. Nevertheless, such symbols are very useful   and convenient for briefly writing verbose definitions like $(-\infty_\zeta,+\infty_\zeta)={\mathbb R}_\zeta$ and $(-\infty_\zeta,+\infty_\zeta)\subset{\mathbb R}_{\zeta'}$, where $\zeta<\zeta'\leq\Omega$ denotes the set of all numbers $x\in{\mathbb R}_{\zeta'}$ such that there is a number $r\in{\mathbb R}_\zeta$ such that $-r<x<r$, in particular, $(-\infty_\omega,+\infty_\omega)\subset{\bf No}$ is the set of all finite Conway's numbers as well as $(-\infty_{\Omega},+\infty_{\Omega})={\bf No}$. We will add two more symbols $\frac{1}{+\infty_\zeta}$ and $\frac{1}{-\infty_\zeta}$ to denote $(\frac{1}{+\infty_\zeta},+\infty_\zeta)$ as the set of all positive numbers in ${\mathbb R}_\zeta$ and $(-\infty_\zeta,\frac{1}{-\infty_\zeta})$ as the set of all negative numbers in ${\mathbb R}_\zeta$ as well as $(-\infty_\zeta,\frac{1}{-\infty_\zeta})$ and $(-\infty_\zeta,\frac{1}{-\infty_\zeta})$ in ${\mathbb R}_{\zeta'}$, $\zeta<\zeta'<\Omega$ like above. Also by $(\frac{1}{-\infty_\omega},\frac{1}{+\infty_\omega})$ we denote the class of all infinitesimal numbers in ${\bf No}$ in the sense of Conway; and by $(-\infty_\Omega,\frac{1}{-\infty_\Omega})$ and $(\frac{1}{+\infty_\Omega},+\infty_\Omega)$ we denote the classes of all negative and positive numbers in ${\bf No}$, respectively. See Appendix below.

{\bf Proposition 2}. {\it $\lim\limits_{0<\alpha<\zeta}x_\alpha=a\in F$ if and only if $(y_\alpha)_{0<\alpha<\zeta}=(x_\alpha-a)_{0<\alpha<\zeta}$ is $\zeta$-infinitely small in $F$. A $\zeta$-sequence is a $\zeta$-infinitely small sequence $(x_\alpha)_{0<\alpha<\zeta}$, $x_\alpha\not=0$  for all $\alpha\ge\alpha_0$ for some $\alpha_0$, if and only if $(\frac{1}{x_\alpha})_{\alpha_0\leq\alpha<\zeta}$ is a $\zeta$-infinitely great seacuence.}

{\bf Proof} is an immediate consequence of Definitions 6 and 7.

Thus we have the following formula for members of $\zeta$-sequence $(x_\alpha)_{0<\alpha<\zeta}$, which converges to $a\in F$:
\begin{equation}
\label{f31}
x_\alpha=a+y_\alpha, 
\end{equation}
where $(y_\alpha)_{0<\alpha<\zeta}$ is $\zeta$-infinitely small in $F$.

\begin{center}
{\bf 4. Fundamental  $\zeta$-sequences in $\zeta$-fields  $F$ of Conway's numbers}
\end{center}

{\bf Definition 9.}
 A $\zeta$-sequence $(x_\alpha)_{0\leq\alpha<\zeta}$ in a  $\zeta$-field $F$ is called  {\it fundamental}  or a {\it Cauchy} $\zeta$-{\it sequence}  if for each positive  number $\varepsilon\in F$  there is an ordinal number $\alpha_0$ such that $|x_\alpha-x_{\alpha'}|<\varepsilon$, for all $\alpha_0\leq\alpha<\alpha'<\zeta$. 

{\bf Definition 10.}
Two $\zeta$-fundamental sequences $(x_\alpha)_{0<\alpha<\zeta}$ and $(y_\alpha)_{0<\alpha<\zeta}$ in a  $\zeta$-field $F$ are $\zeta$-{\it equivalent}, denoted by $(x_\alpha)_{0<\alpha<\zeta}\sim(y_{\alpha})_{0<\alpha<\zeta}$,  if for each  positive number $\varepsilon\in F$  there are ordinal numbers $\alpha_0$ and $\alpha'_0$ such that $|x_\alpha-y_{\alpha'}|<\varepsilon$, for all $\alpha_0\leq\alpha<\zeta$ and all $\alpha'_0\leq\alpha'<\zeta$. 

It is clear that the $\zeta$-equivalence  is an equivalence relation $\sim$ on the set of all fundamental $\zeta$-sequences $(x_\alpha)_{0<\alpha<\zeta}$ in $F$. One has to verify only transitivity what is an easy exerciee. Note also that there can be a case that $\lambda$-sequence is $\zeta$-fundamental for a limit ordinal $\omega\leq\lambda<\zeta$,
but it is always $\zeta$-equivalent to some $\zeta$-fundamental $\zeta$-sequence $(x_\alpha)_{0<\alpha<\zeta}$. That is why we consider here only $\zeta$-sequences.

{\bf Lemma 1}.  {\it If $\zeta$-sequences $(x_\alpha)_{0<\alpha<\zeta}$ and $(y_\alpha)_{0<\alpha<\zeta}$ are $\zeta$-fundamental in $F$, then their sum $(x_\alpha+y_\alpha)_{0<\alpha<\zeta}$ and product $(x_\alpha\cdot y_\alpha)_{0<\alpha<\zeta}$ are also $\zeta$-fundamental in $F$. If in addition there is a positive number $r\in F$ such that $|y_\alpha|>r$, for  almost all $\alpha$, then the quotient  $(\frac{x_\alpha}{y_\alpha})_{\alpha_0\leq\alpha<\zeta}$ is also $\zeta$-fundamental in $F$  for some ordinal $\alpha_0$.}

{\bf Proposition 3}. {\it Each convergent $\zeta$-sequence in a $\zeta$-field $F\subset{\bf No}$ is $\zeta$-fundamental; converse is not true in general. }

Proofs are also an easy exercise. But an example of $\zeta$-fundamental sequence, $\zeta=\omega^{\omega^\mu}$, $0<\mu<\Omega$, which does not converge, is not trivial but we omit this construction here how not essential.

\begin{center}
{\bf 5. $\zeta$-completions $\tilde F$ of $\zeta$-fields  $F$ of Conway's numbers}
\end{center}

In spite of fact that there can be a fundamental $\zeta$-sequence $(x_\alpha)_{0<\alpha<\zeta}$ in a $\zeta$-field $F$ which is not convergent in $F$, it  defines a Conway's number $a=\{L\,|\,R\}$, born on day $\zeta$, with definite $L$ and $R$ such that  $\lim\limits_{0<\alpha<\zeta}x_\alpha=a$.

{\bf Definition 11}. Let $F$ be a subfield of ${\bf No}$ and $\{F\,|\,\}=\zeta=\omega^{\omega^\mu}$, $0\leq\mu<\Omega$.
By ${\tilde F}$ we define the set $\{[(x_\alpha)_{0<\alpha<\zeta}]\}$ of all classes $[(x_\alpha)_{0<\alpha<\zeta}]$ of equivalent fundamental $\zeta$-transfinite sequences $(x_\alpha)_{0<\alpha<\zeta}$ in $F$ and call it the $\zeta$-{\it completion} of $F$.

We define the following operations on ${\tilde F}$.

{\bf Definition 12}. For any two elements $[(x_\alpha)_{0<\alpha<\zeta}]$ and $[(y_\alpha)_{0<\alpha<\zeta}]$ in ${\tilde F}$ we put $[(x_\alpha)_{0<\alpha<\zeta}]+[(y_\alpha)_{0<\alpha<\zeta}]=[(x_\alpha+y_\alpha)_{0<\alpha<\zeta}]$ and $[(x_\alpha)_{0<\alpha<\zeta}]\cdot[(y_\alpha)_{0<\alpha<\zeta}]=[(x_\alpha\cdot y_\alpha)_{0<\alpha<\zeta}]$; if in addition $|y_\alpha|>r>0$, for some ordinal $\alpha_0$, then we put $\frac{[(x_\alpha)_{\alpha_0<\alpha<\zeta}]}{[(y_\alpha)_{\alpha_0<\alpha<\zeta}]}=[(\frac{x_\alpha}{ y_\alpha})_{\alpha_0<\alpha<\zeta}]$ for some ordinal $\alpha_0$.

{\bf Proof}.
One can easily verify the independence of the definition of addition, multiplication and division operations on $\zeta$-equivalent classes from a selection of representatives of the class of equivalent fundamental sequences. It is also clear that ${\tilde F}$ is a field. Moreover, it is  a totally ordered field, $[(x_\alpha)_{0<\alpha<\zeta}]\leq[(y_\alpha)_{0<\alpha<\zeta}]$ if and only if there exists an ordinal $\alpha_0$ such that $x_\alpha\leq y_\alpha$ for every $\alpha_0\leq\alpha<\zeta$  and hence it is a subfield of ${\bf No}$ because we define natural operations on numbers in $\tilde F$ which actually coincide with operations in ${\bf No}$ but not in a bit simple way. 

Now we have to identify each class $[(x_\alpha)_{0<\alpha<\zeta}]$ of ${\tilde F}$ with the unique Conway's number and show that the operations on them are the same as in ${\bf No}$.

Choose any fixed element $(x_\alpha)_{0<\alpha<\zeta}$ of the class $[(x_\alpha)_{0<\alpha<\zeta}]$. Notice also that  $x_\alpha\in F$, for all ${0<\alpha<\zeta}$. Since $(x_\alpha)_{0<\alpha<\zeta}$ is fundamental then for each positive number $\varepsilon\in F$ there is an ordinal $\alpha_0$ such that
\begin{equation}
\label{f071}
x_{\alpha'}-\varepsilon<x_\alpha<x_{\alpha'}+\varepsilon
\end{equation}
for $\alpha_0\leq\alpha<\alpha'<\zeta$.

Let $L'$ be the subset of $F$ of all $l'\in F$ such that there exists an $\alpha_0$ and inequalities $ l'<x_\alpha$ for all $\alpha_0\leq\alpha<\zeta$ and put $R'= F\setminus L'$.

Evidently $L'\not=\emptyset$ because $x_{\alpha'}-\varepsilon\in L'$ and on the other hand $x_{\alpha'}+\varepsilon\notin L'$ and thus, by construction, $x_{\alpha'}+\varepsilon\in R'$ and hence $R'\not=\emptyset$, too.

Notice that for each numbers $l'\in L'$ and $r'\in R'$ we have inequality $l'<r'$; otherwise, $l'\geq r'$, for all $\alpha_0<\alpha<\zeta$ we obtain $x_\alpha>l\geq r'$ and $r'\notin R'$. Contradiction.

Now we can define a Conway's number $a=\{L\,\,|\,R\}$ of ${\bf No}$ such that $a=[(x_\alpha)_{0<\alpha<\zeta}]$. If $L'$ has an extremal element (supremum) $a$ or $R'$ has an extremal
(infimum), which in this case is the same element $a$ of $F$, then we identify this Conway's number $a=\{L\,\,|\,\, R\}$ with the class $[(x_\alpha)_{0<\alpha<\zeta}]\in\tilde F$, where $L=L'\setminus\{a\}$ and $R=R'$ in the first case and $L=L'$ and $R=R'\setminus \{a\}$ in the second case. Notice that these cases depend on the choice of $(x_\alpha)_{0<\alpha<\zeta}$ in $[(x_\alpha)_{0<\alpha<\zeta}]$. Indeed, if the supremum $a$ in $L'$ will be when a chosen $(x_\alpha)_{0<\alpha<\zeta}$ is almost in $R'$ and the minimum $a$ in $R'$ when a chosen $(x_\alpha)_{0<\alpha<\zeta}$ is almost in $L'$. Moreover, in these cases $\lim\limits_{{0<\alpha<\zeta}}(x_\alpha)=a$.  In such a case $a$ can be call  a $\zeta$-rational number of $\tilde F$, by analogous with rational real numbers.

Another possible case when neither $L'$ nor $R'$ have extremal element. Then we obtain a Conway's number $a=\{L\,\,|\,\, R\}$, where $L=L'$ and $R=R'$ and  identify this Conway's number $a$ with the class $[(x_\alpha)_{0<\alpha<\zeta}]\in\tilde F$. Moreover, in this case also $\lim\limits_{{0<\alpha<\zeta}}(x_\alpha)=a$ and we also call $a$ a $\zeta$-irrational numbers of $\tilde F$, call by analogous with irrational real numbers.

By definitions and constructions we see that a Conway's number $a$ does not depend on the choice of element $(x_\alpha)_{0<\alpha<\zeta}$ of the class $[(x_\alpha)_{0<\alpha<\zeta}]$. Thus $a\in \tilde F$ is unique for the class $[(x_\alpha)_{0<\alpha<\zeta}]$.

Moreover, we have the following inequalities:
\begin{equation}
\label{f070}
x_{\alpha'}-\varepsilon\leq a\leq x_{\alpha'}+\varepsilon.
\end{equation}

If we consider another fundamental sequence $(y_\alpha)_{0<\alpha<\zeta}$ and its class $[(y_\alpha)_{0<\alpha<\zeta}]\in\tilde F$,
then we know that $(x_\alpha)_{0<\alpha<\zeta}+(y_\alpha)_{0<\alpha<\zeta}=(x_\alpha+y_\alpha)_{0<\alpha<\zeta}$ and thus $[(x_\alpha)_{0<\alpha<\zeta}]+[(y_\alpha)_{0<\alpha<\zeta}]=[(x_\alpha)+y_\alpha)_{0<\alpha<\zeta}]\in\tilde F$. We have already identified $[(x_\alpha)_{0<\alpha<\zeta}]$ and $[(y_\alpha)_{0<\alpha<\zeta}]$ in $\tilde F$ with Conway's numbers $a$ and $b$, respectively.  Thus, we know that $a+b\in\tilde F$ as sum of two classes of $\zeta$-fundamental sequences. We also know that $a\cdot b\in\tilde F$ and $\frac{a}{b}\in\tilde F$, $b\not=0$ But Conway's definition of sum and product of two numbers $a$ and $b$ is different. Thus, we need the following Proposition to show that these operations coinsied with ours.

{\bf Proposition 4}. {\it Suppose we identified $\tilde F$ as a subset if ${\bf No}$. Then for each $a,b\in \tilde F$ sum $a+b\in\tilde F$ and product $a\cdot b\in\tilde F$ in the sense of $\zeta$-sequences are the same number $a+b\in{\bf No}$ and $a\cdot b\in{\bf No}$ in the sense of Conway's sum and product of numbers $a$ and $b$ in ${\bf No}$.}

{\bf Proof.} Suppose first, that we defined as above $\lim\limits_{{0<\alpha<\zeta}}x_\alpha)=a$, $\lim\limits_{{0<\alpha<\zeta}}(y_\alpha)=b$ and $a+b=\lim\limits_{{0<\alpha<\zeta}}(x_\alpha)+\lim\limits_{{0<\alpha<\zeta}}(y_\alpha)=\lim\limits_{{0<\alpha<\zeta}}(x_\alpha+y_\alpha)$. Conway's  sum, by Definition, is $a+b=\{x^L+y,x+y^L\,|\,x^R+y,x+y^R\}$, where $a=\{x^L\,|\,x^R\}$ and $b=\{y^L\,|\,y^R\}$. Thus, we have to show that $a+b$ in the sense of $\zeta$-sequences is the same number $a+b$ in the sense of Conway's sum of numbers $a$ and $b$. For this purpose we define the Conway's number $a+b$, defined as above by classes $[(x_\alpha)_{0<\alpha<\zeta}]+[(y_\alpha)_{0<\alpha<\zeta}]=[(x_\alpha+y_\alpha)_{0<\alpha<\zeta}]$. Choose a representive element $(x_\alpha+y_\alpha)_{0<\alpha<\zeta}$ of it and define a Conway's number, generated by it as above. That is we consider the following set $L'=\{l'\in F\,|\,\exists\alpha_0\,\&\,l'<x_\alpha\forall\alpha\,|\,\alpha_0\leq\alpha<\zeta\}$ and $R'=F\setminus L'$.

There can be two cases: first $L'$ has maximum or $R'$ has minimum and neither $L'$ nor $R'$ have exitemal element. In the first case maximum of $L'$ is $a+b$ or minimum of $R'$ is $a+b$ and the Conway's number $a+b=\{L\,|\,R\}$, where $L=L'\setminus\{a+b\}$ and $R=R'$ or $L=L'$ and $R=R'\setminus\{a+b\}$, respectively, second when $L'$ and $R'$ have no extremal elements. Then $a+b=\{L'\,|\,R'\}$. Notice that in this case for $a=\{L\,|\,R\}$ and $b=\{L\,|\,R\}$ their above auxiliary sets $L'$ and $R'$ also have no extremal elements or if for $a$ its set $L'$ has maximum, then for $b$ its auxiliary set $R'$ has minimum or vice versa if for $a$ its auxiliary set $R'$ has minimum, then for $b$ its auxiliary set $L'$ has maximum.

In all these cases the sum $a+b$ and product $a\cdot b$ in the sense of $\zeta$-sequences is the Conway's sum $a+b=\{x^L+y,y^L+x\,|\,x^R+y,y^R+x\}$ and product $\{a\cdot b=\{x^Ly+xy^L-x^Ly^L,x^Ry+xy^R-x^Ry^R\,|\,x^Ly+xy^R-x^Ly^R,x^Ry+xy^L-x^Ry^L\}$ because, evidently, for all numbers $x^L+y,y^L+x$ there are inequalities $x^L+y<a+b$ and $y^L<a+b$ as well as for all $x^R+y,y^R+x$ there are inequalities $a+b<x^R+y$ and $a+b<y^R+x$; also for all numbers $x^Ly+xy^L-x^Ly^L,x^Ry+xy^R-x^Ry^R$ there are inequalities $x^Ly+xy^L-x^Ly^L<a\cdot b$, $x^Ry+xy^R-x^Ry^R<a\cdot b$ and all numbers $x^Ly+xy^R-x^Ly^R,x^Ry+xy^L-x^Ry^L$ there are inequalities $x^Ly+xy^R-x^Ly^R<a\cdot b$ and $x^Ry+xy^L-x^Ry^L<a\cdot b$. The latters follow from the evident inequalities: 

$a-x^L>0$ and $b-y^L>0$, hence $(a-x^L)(b-y^L)>0$ and thus $a\cdot b>x^Lb+ay^L-x^Ly^L$;

$x^R-a>0$ and $y^R-b>0$, hence $(x^R-a)(y^R-b)>0$ and thus $a\cdot b>x^Rb+ay^R-x^Ry^R$;
 
 $a-x^L>0$ and $b-y^R<0$,  hence $(a-x^L)(b-y^R)<0$ and thus $a\cdot b<x^Lb+ay^R-x^Ly^R$;
 
 $x^R-a>0$ and $b-y^L>0$,  hence $(x^R-a)(b-y^L)>0$ and thus $a\cdot b<x^Rb+ay^L-x^Ry^L$.
$\Box$

{\bf Theorem 6.}
{\it A $\zeta$-sequence $(x_\alpha)_{0<\alpha<\zeta}$ in ${\tilde F}$ converges in ${\tilde F}$ if and only if it is $\zeta$-fundamental.}

{\bf Proof of necessity.} Let $\zeta$-sequence $(x_\alpha)_{0<\alpha<\zeta}$ in ${\tilde F}$ be convergent in ${\tilde F}$ and thus $\lim\limits_{0<\alpha<\zeta}x_\alpha=a\in \tilde F$. Let $\varepsilon$ be  a positive number in $\tilde F$. For number $\frac{\varepsilon}{2}$ we find an ordinal $\alpha_0$ such that for each $\alpha_0<\alpha<\zeta$ there is inequality:
$$
|x_\alpha-a|<\frac{\varepsilon}{2}
$$

Then for all $\alpha_0<\alpha<\zeta$ and $\alpha_0<\alpha'<\zeta$ we have
$$
|x_\alpha-a|<\frac{\varepsilon}{2}\,\, and\, \, |x_{\alpha'}-a|<\frac{\varepsilon}{2}
$$
and obtain for all $\alpha_0<\alpha<\zeta$ and $\alpha_0<\alpha'<\zeta$ the following inequality:
$$
|x_\alpha-x_{\alpha'}|=|(x_\alpha-a) + (a-x_{\alpha'})|\leq|x_\alpha-a|+|a-x_{\alpha'}|<\frac{\varepsilon}{2}+\frac{\varepsilon}{2}=\varepsilon.
$$

{\bf Proof of sufficiency.} Let $\zeta$-sequence $(x_\alpha)_{0<\alpha<\zeta}$ in ${\tilde F}$ be a $\zeta$-fundamental sequence. We shall prove that there is a number $a\in\tilde F$ such that $\lim\limits_{0<\alpha<\zeta}x_\alpha=a$.

Let us suppose that $(x_\alpha)_{0<\alpha<\zeta}$ is not almost  stationary  $\zeta$-sequence  in ${\tilde F}$, otherwise, $(x_\alpha)_{0<\alpha<\zeta}$ would be converging to this constant and the poof would be done. Then for each positive number $\varepsilon\in{\tilde F}$ there is an ordinal $\alpha_0$ such that
\begin{equation}
\label{f7}
x_{\alpha'}-\varepsilon<x_\alpha<x_{\alpha'}+\varepsilon
\end{equation}
for $\alpha_0\leq\alpha<\alpha'<\zeta$.

Let $\tilde A$ be the subset of ${\tilde F}$ of all $a\in {\tilde F}$ such that there exists an $\alpha_0$ and inequalities $ a<x_\alpha$ for all $\alpha_0\leq\alpha<\zeta$ and put $\tilde A'=\tilde F\setminus \tilde A$.

Evidently $\tilde A\not=\emptyset$ because $x_{\alpha'}-\varepsilon\in \tilde A$ and on the other hand $x_{\alpha'}+\varepsilon\notin \tilde A$ and thus, by construction, $x_{\alpha'}+\varepsilon\in \tilde A'$ and hence $\tilde A'\not=\emptyset$ as well.

Notice that for each numbers $a\in \tilde A$ and $a'\in \tilde A'$ we have inequality $a<a'$; otherwise, $a\geq a'$, for all $\alpha_0<\alpha<\zeta$ we obtain $x_\alpha>a\geq a'$ but $a'\notin \tilde A$.

So we have a Dedekind section $\{\tilde A\,|\,\tilde A'\}$ of $\tilde F$.

 Put now $A=\tilde A\cap F$ and $A'=\tilde A'\cap F$. We consider a case when $A$ and $A'$ have no extremal elements in $F$, otherwise, this extremal element would be the limit of $(x_\alpha)_{0<\alpha<\zeta}$. Thus we conclude that $X$ is co-final in ${\tilde X}$ and $Y$ is co-initial in ${\tilde Y}$.  Now we define a fundamental $\zeta$-transfinite sequence $(y_\alpha)_{0<\alpha<\zeta}$ in $F$ \grqq co-final\grqq\, in $(x_\alpha)_{0<\alpha<\zeta}$ in the following sense: for each $0<\alpha<\zeta$ if $x_\alpha \in {\tilde Y}$, then $y_\alpha\in Y$ and $x_{\alpha'}<y_\alpha<x_\alpha$, where $\alpha'$ is the smallest index such that $\alpha'>\alpha$ and $x_{\alpha'}<x_\alpha$,  and  if $x_\alpha\in {\tilde X}$, $y_\alpha\in X$ and $x_{\alpha'}>y_\alpha>x_\alpha$, where $\alpha'$ is the smallest index such that $\alpha'>\alpha$ and $x_{\alpha'}>x_\alpha$. This is possible, firstly, because either $(x_\alpha)_{0<\alpha<\zeta}$ is co-final in ${\tilde X}$, or co-initial in ${\tilde Y}$, or has the corresponding $\zeta$-subsequences in both; secondly, each $x_\alpha$ in ${\tilde F}$ is a Dedekind section of $F$; thirdly, because $(x_\alpha)_{0<\alpha<\zeta}$ is fundamental and hence such constructed $\zeta$-transfinite sequence $(y_\alpha)_{0<\alpha<\zeta}$ is also fundamental. By Proposition 4, $\lim\limits_{0<\alpha<\zeta}y_\alpha= x\in {\tilde F}$. Since $(x_\alpha)_{0<\alpha<\zeta}$ is also \grqq co-final\grqq\, for $(y_\alpha)_{0<\alpha<\zeta}$ in the above sense we conclude that $\lim\limits_{0<\alpha<\zeta}x_\alpha= a$.$\Box$

 \begin{center}
 {\bf 6. Continuous functions in ${\mathbb R}_\zeta$, $\zeta=\omega^{\omega^\mu}$, $0<\mu<\Omega$}
 \end{center}

The above theory of $\zeta$-fundamental $\zeta$-sequences and their limits in ${\tilde P}_\zeta={\mathbb R}_\zeta$, $\zeta=\omega^{\omega^\mu}$, $0<\mu<\Omega$, allows us to build Mathematical Analysis of continuous functions $f(x)$ defined on $X\subseteq{\mathbb R}_\zeta$ which implies many (not all) classical results of Calculus when $\mu=0$. Let's briefly outline a sketch of such a theory.

{\bf Definition 13.} A point $x_0\in X\subseteq{\mathbb R}_\zeta$ is a {\it $\zeta$-limit point} of $X$ if there is a $\zeta$-sequence $(x_\alpha)_{0<\alpha<\zeta}$ in $X$ such that $\lim\limits_{0<\alpha<\zeta}x_\alpha=x_0$.

{\bf Definition 14.} A function $y=f(x)$ with domain  $X\subseteq{\mathbb R}_\zeta$ is called $\zeta$-continuous at $\zeta$-limit point  $x_0\in X$  if there for each $\zeta$-sequence $(x_\alpha)_{0<\alpha<\zeta}$ in $X$ such that $\lim\limits_{0<\alpha<\zeta}x_\alpha=x_0$ one has $\lim\limits_{0<\alpha<\zeta}f(x_\alpha)=f(x_0)$. We denote it by $\lim\limits_{x\rightarrow x_0}f(x)=0$. If every point in $X$ is a $\zeta$-limit point and $y=f(x)$ is $\zeta$-continuous at $x$, then we say that $y=f(x)$ is a $\zeta$-continuous function in $X$.

 First of all, it is clear that the sum $f(x)+g(x)$, product $f(x)g(x)$ and quotient $\frac{f(x)}{g(x)}$ of two $\zeta$-continuous functions $f(x)$ and $g(x)$ on $X\subseteq{\mathbb R}_\zeta$ are also $\zeta$-continuous; in the latter case we assume that $g(x)\not=0$, for each $x\in X$. Next the composition $g(y)=g(f(x))$ of two $\zeta$-continuous functions $y=f(x)$ and $z=g(y)$ on $X\subseteq{\mathbb R}_\zeta$ and on $Y=f(X)\subseteq {\mathbb R}_\zeta$, respectively, is also $\zeta$-continuous on $X$.

One can also define the derivative, the central notion of differential calculus and then investigate functions for monotony, extremes, convexity, concavity and inflaction, asymptotic behavior, etc. The most surprising thing is that $dim{\mathbb R}_\zeta=0$, where $\zeta=\omega^{\omega^\mu}$ for $\mu>0$ and 1 for $\mu=0$. We omit this easy part.

\begin{center}
{\bf 7. On an $n$th root of a positive number in ${\mathbb R}_\zeta$}
\end{center}

{\bf Theorem 7}. {\it Every positive number $x\in {\mathbb R}_\zeta$, $\zeta=\omega^{\omega^\mu}$, $0<\mu<\Omega$ has a unique positive $n$th root $y\in {\mathbb R}_\zeta$, i.e., $y^n=x$, for each integer $1<n<\omega$.}

We shall show that Theorem 4 implies Theorem 7. Indeed, Theorem 4, proved by Conway, says that every  positive number $x\in {\bf No}$ has a unique positive $n$th root $y\in {\bf No}$, i.e., $y^n=x$, for each integer $1<n<\omega$. It is enough to show that if $x\in{\mathbb R}_\zeta$, then $y$ is also in ${\mathbb R}_\zeta$.

For this purpose we shall describe the Class ${\bf No}\setminus{\mathbb R}_\zeta$. It is clear that if the birthday of $x\in{\bf No}$ is $\alpha<\zeta$, then $x\in{\mathbb R}_\zeta$, and   if the birthday of $x\in{\bf No}$ is $\alpha>\zeta$, then $x\notin{\mathbb R}_\zeta$. The hard case is when  if the birthday of $x\in{\bf No}$ is equal to $\zeta$. In this very case some of numbers are in ${\mathbb R}_\zeta$ and some are not in ${\mathbb R}_\zeta$. The following lemmas describe the state of affairs.

{\bf Lemma 2}. {\it Let ${\mathbb R}_\zeta$ be a $\zeta$-field with $\zeta=\omega^{\omega^\mu}$, $0\leq\mu<\Omega$. Then the set ${\bf No}\setminus{\mathbb R}_\zeta$ is a union of the following intervals of numbers in ${\bf No}$}:

$1)$ {\it $\{y\in{\bf No}\,|\,y>x\, \forall x\in{\mathbb R}_\zeta\}$ and $\{y\in{\bf No}\,|\,y<x\, \forall x\in{\mathbb R}_\zeta\}$};

$2)$ {\it for every strictly increasing $\lambda$-sequence $(x_\alpha)_{0<\alpha<\lambda}$ in ${\mathbb R}_\zeta$ $(\lambda$ is a limit ordinal $\omega\leq\lambda\leq\zeta)$ and the set $X^+=\{x\in{\mathbb R}_\zeta\,|\, x>x_\alpha\, \forall\alpha\in(0,\lambda)\}$ this interval is defined by the following inequalities: $x_\alpha<y$, for all $\alpha\in(0,\lambda)$, and  $y<x$, for all $x\in X^+$};

$3)$ {\it for every strictly decreasing $\lambda$-sequence $(x_\alpha)_{0<\alpha<\lambda}$ in ${\mathbb R}_\zeta$ $(\lambda$ is a limit ordinal $\omega\leq\lambda\leq\zeta)$ and the set $X^-=\{x\in{\mathbb R}_\zeta\,|\, x<x_\alpha,\, \forall\alpha\in(0,\lambda)\}$ this interval is defined by the following inequalities: $x<y$, for all $x\in X^-$ and  $y<x_\alpha$, for all $\alpha\in(0,\lambda)$}.

{\it Moreover, in the first case when $|y|>\alpha$ for all $0\leq\alpha<\zeta$, then $y\in{\bf No}\setminus{\mathbb R}_\zeta$. In the second and third cases when a strictly increasing or strictly decreasing $\lambda$-sequence $(x_\alpha)_{0<\alpha<\lambda}$ is convergent in ${\mathbb R}_\zeta$, say to $\hat x\in{\mathbb R}_\zeta$, then every number $y=\hat x\mp\delta$, where $0<\delta<\frac{1}{\alpha}$ for all $0<\alpha<\zeta$, is an element of ${\bf No}\setminus{\mathbb R}_\zeta$. And when a strictly increasing or strictly decreasing $\lambda$-sequence $(x_\alpha)_{0<\alpha<\lambda}$ is not convergent in ${\mathbb R}_\zeta$, then $y\in{\bf No}\setminus{\mathbb R}_\zeta$ if in the normal form $y=\sum\limits_{\beta<\gamma}\omega^{y_\beta}r_\beta$ there is a $\beta$-term $\omega^{y_\beta}r_\beta$ such that $\omega^{y_\beta}\notin{\mathbb R}_\zeta$}.

{\bf Proof}. The first case   $\{y\in{\bf No}\,|\,y>x\, \forall x\in{\mathbb R}_\zeta\}$ and $\{y\in{\bf No}\,|\,y<x\, \forall x\in{\mathbb R}_\zeta\}$ is evident. Indeed, these sets are evidently subsets of ${\bf No}$ and no element $y$ of it is an element of ${\mathbb R}_\zeta$. Moreover each of them is convex, i.e., if $y_1<y_2<y_3$ and $y_1$ and $y_3$ are elements of one of them, then $y_2$ is also an element of this set and thus they are intervals of elements in ${\bf No}$.

The second and third cases when a convergent strictly increasing or strictly decreasing $\lambda$-sequences $(x_\alpha)_{0<\alpha<\lambda}$ is convergent to $\hat x\in{\mathbb R}_\zeta$   is also evident because $|\hat x-y|=|\hat x-x_\alpha+x_\alpha-y|<|\hat x-x_\alpha|+|x_\alpha-y|<\frac{2}{\alpha}$. for all $0<\alpha<\lambda$. Put $\delta=\hat x-y$ or $\delta=y-\hat x$, respectively, and thus $\hat x-\delta$ or $y=\hat x+\delta$, respectively, with $0<\delta<\frac{1}{\alpha'}$ for all $0<\alpha'<\lambda$, where $\alpha=2\alpha'$.

The second and third cases when a  strictly increasing or strictly decreasing $\lambda$-sequences $(x_\alpha)_{0<\alpha<\lambda}$ is not convergent in ${\mathbb R}_\zeta$   is also evident because because otherwise, $y=\sum\limits_{\beta<\gamma}\omega^{y_\beta}r_\beta$ with all $\beta$-terms $\omega^{y_\beta}r_\beta$, $0<\beta<\gamma$, such that $\omega^{y_\beta}\in{\mathbb R}_\zeta$,  should be in ${\mathbb R}_\zeta$. $\Box$

{\bf Lemma 3}. For each  number $x\in{\bf No}$ such that $x\not=0$ and all positive real numbers $r\in{\mathbb R}\subset{\bf No}$ there are the following inequalities:
\begin{equation}
\label{f0701}
\omega^x<r,\qquad\mbox{if}\,\,x<0,\qquad\mbox{and}\qquad
\omega^x>r,\qquad\mbox{if}\,\,x>0.
\end{equation}

{\bf Proof}.
First of all, we shall prove that for any ordinal $0<\alpha<\Omega$ there is an inequality $r<\omega^{\frac{1}{\alpha}}$ for each real number $r>0$. Consider the number $a=\{0\,|\,\frac{1}{\alpha}\}$.  By formula $(\ref{f314})$, $\omega^a=\{0,r\cdot\omega^{ a^L}\,|\,r\cdot\omega^{a^R}\}=\{0,r\cdot\omega^0\,|\,r\cdot\omega^{\frac{1}{\alpha}}\}=\{0,r\,|\,r\cdot\omega^{\frac{1}{\alpha}}\}$, where $r$ denotes a variable ranging over all positive reals. Thus, $r$ is a typical member of a left option of $\omega^a$ and $r\cdot\omega^{\frac{1}{\alpha}}$ is a typical member of a right option of $\omega^a$. In particular, for each member $r$ of a left option  of $\omega^a$ and for the member $\omega^{\frac{1}{\alpha}}$ of a right  option of $\omega^a$  when $r=1$, by Theorem 2, we obtain $r<\omega^a<\omega^{\frac{1}{\alpha}}$ what implies  an inequality $r<\omega^{\frac{1}{\alpha}}$. 

Let $x$ be an arbitrary positive number in ${\bf No}$. 

If $1\leq x$, then trivially $\omega^x\geq\omega^1=\omega$ and $\omega>r$ for each  positive  real number $r$. 

If $x<1$, then there exists an ordinal $1<\alpha<\Omega$ such that $\frac{1}{\alpha}<x$. Indeed, put $\alpha'=\{\mbox{ordinals}\,\,\beta<\frac{1}{x}\,|\,\}$ and $\alpha=\alpha'+1$. Clearly, $\alpha>\frac{1}{x}$ and $\frac{1}{\alpha}<x$. Then $\omega^{\frac{1}{\alpha}}<\omega^x$ and with above result $r<\omega^{\frac{1}{\alpha}}$ for each positive real number $r$ we obtain an inequality $r<\omega^x$ for each  positive real number $r$.

Let $x$ be an arbitrary negative number in ${\bf No}$. Then $-x$ is a positive number in ${\bf No}$. We have already proved that there is an inequality $r<\omega^{-x}$ for each positive real number $r$. Thus $\omega^x<\frac{1}{r}$ for each positive real number $r$. $\Box$

{\bf Lemma 4}. {\it Let $g_m(x)=b_0x^m+b_1x^{m-1}+...+b_{m-1}x+b_m$, $0<m<\omega$, be an arbitrary polynomial with coefficients $b_0,b_1,...,b_m$ in ${\mathbb R}_\zeta$, $b_0\not=0$ and $\bar x$ be any fixed number in {\bf No}. If $|\bar x|>\alpha$ for all $0<\alpha<\zeta$, then  there is an inequality
\begin{equation}
\label{f0409}
|g(\bar x)|>\alpha,
\end{equation}
for all $0<\alpha<\zeta$}.

{\bf Proof}. We shall prove it by induction. For $m=1$ one has $g_1(x)=b_0x+b_1$, $b_0\not=0$, and $b_0,b_1\in{\mathbb R}_\zeta$, suppose the contrary, i.e., $|\bar x|>\alpha$ for all $0<\alpha<\zeta$ and $|g_1(\bar x)|\leq \alpha_0$ for some $0<\alpha_0<\zeta$. Then $|b_0\bar x+b_1|\leq \alpha_0$ and hence $-\alpha_0-b_1\leq b_0\bar x\leq \alpha_0-b_1$. If  $b_0<0$ it implies $\frac{-\alpha_0-b_1}{b_0}\geq\bar x\geq \frac{\alpha_0-b_1}{b_0}$ and if  $b_0>0$ it implies $\frac{-\alpha_0-b_1}{b_0}\leq\bar x\leq \frac{\alpha_0-b_1}{b_0}$. Since $\alpha_0,b_0,b_1\in{\mathbb R}_\zeta$ then $\frac{-\alpha_0-b_1}{b_0}$ and $\frac{\alpha_0-b_1}{b_0}$ are elements of ${\mathbb R}_\zeta$. 

Consider now any ordinal $\alpha_1\in{\mathbb R}_\zeta$ such that $\alpha_1>\max\{|\frac{-\alpha_0-b_1}{b_0}|,|\frac{\alpha_0-b_1}{b_0}|\}$. Then, by both cases of above inequalities we obtain $|\bar x|<\alpha_1$ what contradicts with supposition that $|\bar x|>\alpha$ for all $0<\alpha<\zeta.$

Suppose now that for each $1\leq n\leq m-1$ there are inequalities $|g_n(\bar x)|> \alpha$ for all $0<\alpha<\zeta$, and prove that for $n=m$ there is the same inequality. Indeed, suppose the contrary, i.e., $|\bar x|>\alpha$ for all $0<\alpha<\zeta$ and $|g_m(\bar x)|=|b_0\bar x^m+b_1\bar x^{m-1}+...+b_{m-1}\bar x+b_m|\leq \alpha_0$ for some $0<\alpha_0<\zeta$, what is the same that $|\bar x\cdot g_{m-1}(\bar x)+b_m|\leq \alpha_0$, where $g_{m-1}(x)=b_0 x^{m-1}+b_1 x^{m-2}+...+b_{m-1}$. Hence $-\alpha_0-b_m\leq \bar x\cdot g_{m-1}(\bar x)\leq \alpha_0-b_m$. Since $-\alpha_0-b_m,\alpha_0-b_m\in{\mathbb R}_\zeta$ we consider any ordinal $\alpha_1\in{\mathbb R}_\zeta$ such that $\alpha_1>\max\{|-\alpha_0-b_m|,|\alpha_0-b_m|\}$. Then $|\bar x\cdot g_{m-1}(\bar x)|<\alpha_1$ what contradicts with supposition that $|\bar x|>\alpha$  for all $0<\alpha<\zeta$ and inductive supposition $|g_{m-1}(\bar x)|>\alpha$ for all $0<\alpha<\zeta$, in particular, for $\alpha=1$, because multiplication of latter inequalities $|\bar x|>\alpha$ and $g_{m-1}(\bar x)>1$ imply $|\bar x\cdot g_{m-1}(\bar x)|>\alpha\cdot 1=\alpha$ for all $0<\alpha<\zeta$.

{\bf Lemma 5}. Let $x$ be a number in $\mathbb R_\zeta$, $\zeta=\omega^{\omega^\mu}$, $0\leq\mu<\Omega$, such that  $\alpha$, $0\leq\alpha\leq\omega^\mu$ is the birthday of $x$. Then $\omega^x\in\mathbb R_\zeta$. In particular, when $\zeta$ is an $\varepsilon$-number, then for each number $x\in{\mathbb R}_\zeta$ one has $\omega^x\in{\mathbb R}_\zeta$.

{\bf Proof}. It is easy to see that if $\alpha$ is the birthday of $x\in{\bf On}$, then  $\omega^\alpha$ is the birthday of $y=\omega^x\in{\bf On}$ (see $\cite{l3}$, p. 8-14, 31). 

By supposition of Lemma 2, $x\in{\mathbb R}_\zeta$ and its birthday is $0\leq\alpha<\omega^\mu$, $0\leq\mu<\Omega$. Then $\omega^x\in{\bf No}$ was born on day $\omega^\alpha$. Consequently, $\omega^\alpha<\omega^{\omega^\mu}=\zeta$ and thus $\omega^x\in{\mathbb R}_\zeta$.

Note that, in general, for some $x\in{\mathbb R}_\zeta$, $\zeta=\omega^{\omega^\mu}$, $0\leq\mu<\Omega$, $\omega^x\notin{\mathbb R}_\zeta$. Simple examples: when $\mu=0$ then only for one number $x=0\in{\mathbb R}_\omega$ one has $\omega^x\in{\mathbb R}_\omega$. For all other numbers it is wrong. For $\mu=1$, $x\in{\mathbb R}_{\omega^\omega}$ whose birthday is $0<\alpha<\omega$, we have $\omega^x\in{\mathbb R}_{\omega^{\omega}}$ but for $x=\omega^2$ (say) $\omega^{\omega^2}\notin{\mathbb R}_{\omega^{\omega}}$.

Nevertheless, if $\zeta$ is an $\varepsilon$-number, i.e., $\omega^\zeta=\zeta$, then for each $x\in{\mathbb R}_\zeta$ one has $\omega^x\in{\mathbb R}_\zeta$. Indeed, $\zeta=\omega^\zeta=\omega^{\omega^\zeta}$ and hence $\mu=\zeta$. Thus, the proof was above. In particular, each initial ordinal number $\omega_\alpha$, $0<\alpha<\Omega$, is an $\varepsilon$-number and thus for every $x\in{\mathbb R}_{\omega_\alpha}$ one has $\omega^x\in{\mathbb R}_{\omega_\alpha}$. $\Box$

{\bf Lemma 6}. {\it Let $g_m(x)=b_0x^m+b_1x^{m-1}+...+b_{m-1}x$, $0<m<\omega$, be an arbitrary polynomial with coefficients $b_0,b_1,...,b_{m-1}$ in ${\mathbb R}_\zeta$, $b_0\not=0$ and $\bar x$ be any fixed number in {\bf No}.  If $|\bar x|<\frac{1}{\alpha}$ for all $0<\alpha<\zeta$, then  there is an inequality
\begin{equation}
\label{f1409}
|g(\bar x)|<\frac{1}{\alpha}
\end{equation}
 for all $0<\alpha<\zeta$}.

{\bf Proof}. 
It is clear that $|b_{n}\bar x^n|<\frac{1}{\alpha}$ for all $0<\alpha<\zeta$, $0\leq n<m$. Otherwise if $|b_{n}\bar x^n|\geq\frac{1}{\alpha_0}$
for some $0\leq n<m$ such that $b_n\not=0$, and $0<\alpha_0<\zeta$, then $|\bar x^n|\geq\frac{1}{\alpha_0|b_n|}$. Take any ordinal $\alpha_1$ such that $\frac{1}{\alpha^n_1}<\frac{1}{\alpha_0|b_n|}\leq|\bar x^n|$ and hence $\frac{1}{\alpha_1}<|\bar x|$, what is in contradiction with supposition of Lemma 3.

Consider now the following inequality
\begin{equation}
\label{f2409}
|b_0x^m+b_1x^{m-1}+...+b_{m-1}x|\leq |b_0x^m|+|b_1x^{m-1}|+...+|b_{m-1}x|\leq m\cdot\max\limits_{0\leq n<m}|b_nx^{m-n}|,
\end{equation}
which implies the following inequality
\begin{equation}
\label{f3409}
|b_0x^m+b_1x^{m-1}+...+b_{m-1}x|<\frac{1}{\alpha}
\end{equation}
for all $0<\alpha<\zeta$ because $m\cdot|b_n\bar x^n|=|m\cdot b_n\bar x^n|<\frac{1}{\alpha}$ for all $0<\alpha<\zeta$, where $|b_n\bar x^n|=\max\limits_{0\leq n<m}\{|b_0x^m|,|b_1x^{m-1}|,...,|b_{m-1}x|\}$.

 {\bf Theorem 8.} {\it Let ${\mathbb R}_\zeta$ be a $\zeta$-field with $\zeta=\omega^{\omega^\mu}$, $0\leq\mu<\Omega$. Then every number $\bar x\in{\bf No}\setminus{\mathbb R}_\zeta$ is transcendental over ${\mathbb R}_\zeta$. e. i., for every $P(x)=a_0x^n+a_1x^{n-1}+...+a_{n-1}x+a_n$ with $a_i\in{\mathbb R}_\zeta$, $0\leq i\leq n$, $0<n<\omega$, $a_0\not=0$, one has $P(\bar x)\not=0$}.

 {\bf Proof.} Since $\bar x\in{\bf No}\setminus {\mathbb R}_\zeta$ then, by Lemma 2, there are three possible situations.

$1)$ $\bar x\in\{y\in{\bf No}\,|\,y>x\, \forall x\in{\mathbb R}_\zeta\}$ or $\bar x\in\{y\in{\bf No}\,|\,y<x\, \forall x\in{\mathbb R}_\zeta\}$. Then $|\bar x|>\alpha$ for all $0<\alpha<\zeta$ and, by Lemma 5, $|P_n(\bar x)|>\alpha$ for all $0<\alpha<\zeta$. Thus $P_n(\bar x)\not=0$.

$2)$  $a_\alpha<\bar x<a$, for all $\alpha\in(0,\lambda)$ and all $a\in X^+$, where $(a_\alpha)_{0<\alpha<\lambda}$   is a strictly increasing $\lambda$-sequence in ${\mathbb R}_\zeta$ $(\lambda$ is a limit ordinal $\omega\leq\lambda\leq\zeta)$ and the set $X^+=\{a\in{\mathbb R}_\zeta\,|\, a>a_\alpha\, \forall\alpha\in(0,\lambda)\}$.

There can be two differend cases: 

$a)$ $(x_\alpha)_{0<\alpha<\lambda}$ is $\zeta$-fundamental and 

$b)$ $(x_\alpha)_{0<\alpha<\lambda}$ is not $\zeta$-fundamental.

In the first case there are the following inequalities: $a_{\lambda'}<\bar x<\lim\limits_{0<\lambda'<\lambda}a_{\lambda'}$ for all $0<\lambda'<\lambda$. Denote $\lim\limits_{0<\lambda'<\lambda} a_{\lambda'}=\hat a$. Notice that $\hat a\in{\mathbb R}_\zeta$ and $\hat a-\bar x<\frac{1}{\alpha}$ for all $0<\alpha<\zeta$. Then $\bar x=\hat a-\varepsilon$ and $0<\varepsilon<\frac{1}{\alpha}$ for all $0<\alpha<\lambda$. We shall prove that $P_n(\bar x)\not=0$. Indeed,  $P_n(\bar x)=\bar x^n+a_1\bar x^{n-1}+...+a_{n-1}\bar x+a_n=(\hat x-\varepsilon)^n+a_1(\hat x-\varepsilon)^{n-1}+...+a_{n-1}(\hat x-\varepsilon)+a_n=\hat x^n+C^1_n\hat x^{n-1}(-\varepsilon)+...+C^{n-1}_n\hat x(-\varepsilon)^{n-1}+(-\varepsilon)^n+a_1\hat x^{n-1}+a_1C^1_{n-1}\hat x^{n-2}(-\varepsilon)+...+a_1C^{n-2}_{n-1}\hat x(-\varepsilon)^{n-2}+a_1(-\varepsilon)^{n-1}+...+a_{n-1}\hat x+a_{n-1}(-\varepsilon)=A_0\varepsilon^n+A_1\varepsilon^{n-1}+..+A_{n-1}\varepsilon+\hat x^n+a_n$, where $A_i\in{\mathbb R}_\zeta$, $0\leq i\leq n$. Then to prove $P(\bar x)\not=0$ is the same as to prove that $\varepsilon^n+\frac{A_1}{A_0}\varepsilon^{n-1}+..+\frac{A_{n-1}}{A_0}\varepsilon+\frac{\hat x^n+a_n}{A_0}\not=0$. Note that $A_0=\pm 1$ and depends on $n$, i.e., if $n=2k$, then $A_0=1$ and if $n=2k+1$, then $A_0=-1$.

On one hand,  $|\varepsilon^n+\frac{A_1}{A_0}\varepsilon^{n-1}+..+\frac{A_{n-1}}{A_0}\varepsilon|\leq|\varepsilon^n|+|\frac{A_1}{A_0}\varepsilon^{n-1}|+..+|\frac{A_{n-1}}{A_0}\varepsilon|\leq n\cdot\max_{0<i<n}|\frac{A_i}{A_0}\varepsilon^{n-i}|<\frac{1}{\alpha}$ for all $0<\alpha<\zeta$.

On the other hand, since $\varepsilon\notin{\mathbb R}_\zeta$ then $\varepsilon\not=-\frac{A_1}{A_0}$ and thus $\varepsilon+\frac{A_1}{A_0}\not=0$. Moreover, as above $0<|\varepsilon^2+\frac{A_1}{A_0}\varepsilon|<\frac{1}{\alpha}$   for all $0<\alpha<\zeta$ and analogously $\varepsilon^2+\frac{A_1}{A_0}\varepsilon+\frac{A_2}{A_0}\not=0$ and $0<|\varepsilon^3+\frac{A_1}{A_0}\varepsilon^2+\frac{A_2}{A_0}\varepsilon|<\frac{1}{\alpha}$ for all $0<\alpha<\zeta$ and so on up to $0<|\varepsilon^n+\frac{A_1}{A_0}\varepsilon^{n-1}+\frac{A_2}{A_0}+...+\frac{A_{n-1}}{A_0}\varepsilon+\frac{A_n}{A_0}|<\frac{1}{\alpha}$ for all $0<\alpha<\zeta$. Thus, $\varepsilon^n+\frac{A_1}{A_0}\varepsilon^{n-1}+\frac{A_2}{A_0}+...+\frac{A_{n-1}}{A_0}\varepsilon+\frac{A_n}{A_0}\not=0$, where $A_n=\frac{\hat x^n+a_n}{A_0}$.

 A case when $a<\bar x<a_{\lambda'}$, for all $\lambda'\in(0,\lambda)$ and all $a\in X^-$, where $(a_\alpha)_{0<\lambda'<\lambda}$   is a strictly decreasing $\lambda$-sequence in ${\mathbb R}_\zeta$ $(\lambda$ is a limit ordinal $\omega\leq\lambda\leq\zeta)$ and the set $X^-=\{a\in{\mathbb R}_\zeta\,|\, a<a_\alpha\, \forall\alpha\in(0,\lambda)\}$, is absolutely analogous.
 
 Consider now a case when a strictly increasing $\lambda$-sequence $(x_{\lambda'})_{0<\lambda'<\lambda}$ in ${\mathbb R}_\zeta$ $(\lambda$ is a limit ordinal $\omega\leq\lambda\leq\zeta)$ is not $\zeta$-fundamental.

 Let now $\bar x$ be a number such that $x_{\lambda'}<\bar x<a$ for all $0<\lambda'<\lambda$ and all $a\in X^+$. We shall prove that $P_n(\bar x)\not=0$. Without loss of generality we can suppose that $P_n(x)=x^n+a_1x^{n-1}+...+a_{n-1}x+a_n$ and prove that $P_n(\bar x^n)=\bar x^n+a_1\bar x^{n-1}+...+a_{n-1}\bar x+a_n\not=0$.

 Consider now the smaller index $0\leq\beta_0<\gamma$ of a term $\omega^{y_{\beta_0}}\cdot r_{\beta_0}$ in the normal form $\sum\limits_{\beta<\gamma}\omega^{y_\beta}\cdot r_\beta$ of $\bar x$ such that $\omega^{y_{\beta_0}}\notin{\mathbb R}_\zeta$, i. e., $\bar x=\sum\limits_{0\leq\beta<\beta_0}\omega^{y_\beta}\cdot r_\beta+\sum\limits_{\beta_0\leq\beta<\gamma}\omega^{y_\beta}\cdot r_\beta$ and hence $\bar x=A+\omega^{y_{\beta_0}}\cdot r_{\beta_0}+\bar x_1$, where there are the following inequalities $-\bar x+A<\bar x_1<\bar x-A$.
 
 Clearly, such item $\omega^{y_{\beta_0}}\cdot r_{\beta_0}$ exists because otherwise $\bar x$ should be an element of ${\mathbb R}_\zeta$.
 
 The idea to prove it  is the following. We also consider the first terms $\omega^{z^{k}_0}\cdot q^{k}_0$, $1\leq k\leq n$, in the normal forms of all items $a_1x^{n-1}$,...,$a_{n-1}x$, $a_n$ of $P_n(\bar x)$ and show that its sum cannot delete $\omega^{y_{\beta_0}}\cdot r_{\beta_0}$.
 
 Indeed, first of all notice that for all $\omega^{y_{\beta_0}}$, $\omega^{2\cdot y_{\beta_0}}$, ..., $\omega^{n\cdot y_{\beta_0}}$ there are the following inequalities: $x_{\lambda'}-A<\omega^{k\cdot y_{\beta_0}}<a-A$ for all $0<\lambda'<\lambda$ and all $a\in X^+$, because $y_{\beta_0}\not=-\frac{b}{k}$ for all $1<k<\omega$ and all $b\in{\mathbb R}_\zeta$ such that $b<\omega^\mu$. Otherwise, $\omega^{k\cdot y_{\beta_0}}\in{\mathbb R}_\zeta$.
 
 Now $P_n(\bar x)=P_n(\sum\limits_{0\leq\beta<\beta_0}\omega^{y_\beta}\cdot r_\beta+\sum\limits_{\beta_0\leq\beta<\gamma}\omega^{y_\beta}\cdot r_\beta)=P_n(B+\omega^{k\cdot y_{\beta_0}}+\bar x_1)=C+\omega^{n\cdot y_{\beta_0}}+b_1\omega^{(n-1)\cdot y_{\beta_0}}+...+b_{n-1}\omega^{ y_{\beta_0}}+b_n$+ smaller, where $B=\sum\limits_{0\leq\beta<\beta_0}\omega^{y_\beta}\cdot r_\beta$, $C=B^n+a_1B^{n-1}+...+a_{n-1}B+a_n\in{\mathbb R}_\zeta$ and $b_1,...,b_n$ are adduced coefficients which are evidently in ${\mathbb R}_\zeta$.
 
 Since all $\omega^{z^{k}_0}\cdot\omega^{m\cdot y_{\beta_0}}=\omega^{z^{k}_0+m\cdot y_{\beta_0}}\not=\omega^{t\cdot y_{\beta_0}}$ for all $0<k,m,t<\omega$ we conclude that $P_n(\bar x)\not=0$. $\Box$

{\bf Lemma 7}. Let $x$ be a number $x\in{\bf No}$ whose birthday is $\gamma$. Then the birthday of number $\omega^x$ is $\omega^\gamma$. 

{\bf Proof}. We shall prove first that Lemma 4 is true when $x$ is an ordinal number. Clear that if $x=0$, then $\omega^0=1$ and thus it is the simplest case of Lemma 4 because the birthday of $0$ is $0$ and the birthday of $1$ is $1$. 

Since the birthday of every ordinal $\alpha$ is $\alpha$ because $\alpha=\{\mbox{ordinals}\,\,\beta<\alpha\,|\,\}$ and, by definition of power of ordinal numbers $\gamma^\alpha$ ($\gamma^0=1$; $\gamma^{\alpha+1}=\gamma^\alpha\cdot\gamma$; $\gamma^\alpha=\lim\limits_{0\leq\alpha'<\alpha}\gamma^{\alpha'}$, where $\alpha$ is a limit ordinal and $\lim\limits_{0\leq\alpha'<\alpha}\gamma^{\alpha'}$ is the smallest ordinal which is greater than $\gamma^{\alpha'}$, for all $0\leq\alpha'<\alpha$), $\omega^\alpha$ is also an ordinal number and thus its birthday is $\omega^\alpha$. Notice also that the set $\{\omega^{\alpha'}\}_{0\leq\alpha'<\alpha}$ is co-final in the set $\{\beta\}_{0\leq\beta<\omega^\alpha}$. Indeed, since $\omega^\alpha$ is the smallest ordinal which is greater than $\omega^{\alpha'}$, for all $0\leq\alpha'<\alpha$), then for any $\beta<\omega^\alpha$ there is $\alpha'<\alpha$ such that $\beta<\omega^{\alpha'}$.

We repeat it for the definition of power operation $\omega^x$, $x\in{\bf No}$ and $x\geq 0$, given by Conway. Indeed, $\omega^\alpha=\{0,r\omega^\beta\,|\,\}$ where $r$ denotes a variable ranging over all positive real numbers. In particular, $\omega^\alpha=\{0,n\omega^\beta\,|\,\}$ where $n$ denotes a variable ranging over all natural numbers. We have to show that the set $\{n\omega^\beta\}_{0\leq n<\omega}$ is co-final in the set $\{\qquad\mbox{ordinals}\,\,\,\beta'<\omega^\alpha\}$ and thus $\omega^\alpha$ is an ordinal whose birthday is $\omega^\alpha$.

If $\alpha=\alpha'+1$ for some ordinal $\alpha_0\leq\alpha'<\alpha$, then the birthday of $\omega^\alpha$ is $\omega^\alpha$, because $\{\omega^{\alpha'}, 2\omega^{\alpha'},...,n\omega^{\alpha'},...\,|\,\}=\omega^\alpha$.

If $\alpha$ is a limit ordinal, then the birthday of $\omega^\alpha$ is also $\omega^\alpha$, because, by definition of power $\omega^\alpha$ of ordinal numbers, $\omega^\alpha=\lim\limits_{0\leq\alpha'<\alpha}\omega^{\alpha'}$, where the latter is the smallest ordinal number which is greater rthan $\omega^{\alpha'}$, $0\leq\alpha'<\alpha$. But all ordinals $\omega^{\alpha'}$, $0\leq\alpha'<\alpha$ are in the Left option of $\omega^\alpha=\{0,r\omega^{\alpha'}\,|\,\}$ and co-final in the set $\{r\omega^{\alpha'}\}_{0\leq\alpha'<\alpha}$. Thus, the birthday of $\omega^\alpha$ is exactly $\omega^\alpha$.

Let $\{x^L\,|\,x^R\}$ be a birthday form of $x$, i.e., all $x^L,x^R\in O_\gamma$ and $(\{x^L\},\{x^R\})$ is a Dedekind section of $O_\gamma$. It exists because the birthday of $x$ is a day $\gamma$. By formula $(\ref{f314})$, $\omega^\gamma=\{0,r\omega^{x^L}\,|\,r\omega^{x^R}\}$ where $r$ denotes a variable ranging over all positive real numbers.

It is clear that there exists an ordinal number $\alpha_0$ such that for every ordinal $\alpha_0\leq\alpha<\gamma$ we have $\alpha\in O_\gamma$ and $\alpha$ is in a Right option of $x$. Hence for every positive real number $r$ the number $r\omega^\alpha$ is in a Right option of $\omega^x$. In particular, $n\omega^\alpha$ is in a Right option of $\omega^x$ for every positive natural number $n$. Thus, if $\gamma=\alpha'+1$ for some ordinal $\alpha_0\leq\alpha'<\gamma$, then the birthday of $\omega^x$ is $\omega^\gamma$, because $\{\omega^{\alpha'}, 2\omega^{\alpha'},...,n\omega^{\alpha'},...\,|\,\}=\omega^\gamma$.

If $\gamma$ is a limit ordinal, then the birthday of $\omega^x$ is also $\omega^\gamma$, because for all ordinals $\omega^\alpha$, $\alpha_0\leq\alpha<\gamma$, we have $\{\omega^{\alpha_0}, \omega^{\alpha_0+1},...,\omega^{\alpha},...\,|\,\}=\omega^\gamma$, i.e., the birthday of each such number $\omega^\alpha$ from Right option of $\omega^x$ is less than $\gamma$ as well as for each ordinal $\gamma'<\gamma$ there is an ordinal number $\alpha$ such that $\gamma'<\omega^\alpha<\gamma$ and thus the birthday of $\omega^x$ is exactly $\gamma$.

{\bf Proposition 5.} For each number $x\in{\mathbb R}_\zeta$ whose birthday is $0\leq\alpha<\omega^\mu$, the number $\omega^x$ is also in ${\mathbb R}_\zeta$. Moreover, if $\zeta$ is an $\varepsilon$-number, i.e., $\omega^\varepsilon=\varepsilon$, then for each number $x\in{\mathbb R}_\zeta$ whose birthday is less than $\zeta$ the power $\omega^x$ is always in ${\mathbb R}_\zeta$.

The proof is evident.

{\bf Lemma 8}. Let $x$ and $y$ be any positive numbers in ${\bf No}$ such that $x\not=y$. Then $\omega^x$ is not commensurate with $\omega^y$.

{\bf Proof}. Suppose the opposite, i.e., there exists a natural number $0<n<\omega$ such that $\omega^x<n\cdot\omega^y$ and $\omega^y<n\cdot\omega^x$. The latter inequalities imply the following once $\omega^{x-y}<n$ and $\omega^{y-x}<n$.

Since for $x\not=y$ there are only two possibilities: either $x<y$ or $y<x$ what is the same that either $y-x>0$ or $x-y>0$. Thus one of the inequalities $\omega^{x-y}<n$ and $\omega^{y-x}<n$ is in contradiction with Lemma 3. $\Box$

{\bf Corollary 3.} If $\omega^x\not=\omega^y$, then $\omega^x$ is not commensurate with $\omega^y$.

{\bf Proof}.  $\omega^x\not=\omega^y$ implies  $x\not=y$.

{\bf Lemma 9}. ${\mathbb R}_\zeta=S^{-1}\bar{O}_\zeta$, where $\bar O_\zeta$ is the closure of $O_\zeta$ in ${\mathbb R}_\zeta$.

{\bf Proof}. By Definition, ${\mathbb R}_\zeta$ is a completion $\bar P_\zeta$ of the localization $P_\zeta=S^{-1}O_\zeta$ of $O_\zeta$ at zero, where $S=O_\zeta\setminus\{0\}$. 

Consider $S'=\bar O_\zeta\setminus\{0\}\subseteq{\mathbb R}_\zeta$. We shall show that $S'^{-1}\bar O_\zeta\subseteq{\mathbb R}_\zeta$ and ${\mathbb R}^+_\zeta\subseteq S'^{-1}\bar O^+_\zeta$. 

If $x=0$, then it is a member of both sets.

Let $x\not=0$ be a member of $S'^{-1}\bar O_\zeta$. 

If $x\in O_\zeta$, then $x\in P_\zeta$ and hence $x\in{\mathbb R}_\zeta$. 

If $x\notin O_\zeta$, then there is a $\zeta$-sequence $(x_\alpha)_{0<\alpha<\zeta}$ which converges to $x$ and $x_\alpha\in O_\zeta$ for all $0<\alpha<\zeta$. Thus $x_\alpha\in P_\alpha$ for all $0<\alpha<\zeta$ and hence $x\in{\mathbb R}_\zeta$.

Let now $x$ be a member of ${\mathbb R}_\zeta$. 

If $x\in O_\zeta$, then $x\in\bar O_\zeta$ and hence $x\in S'^{-1}\bar O_\zeta$. If $x\in P_\zeta$ and $x\notin O_\zeta$, then $\frac{1}{x}\in O_\zeta$ and hence $\frac{1}{x}\in \bar O_\zeta$ as well as $x\in \bar O_\zeta$. Hence $x\in S'^{-1}\bar O_\zeta$. 

If $x\notin P_\zeta$, then then there is a $\zeta$-sequence $(x_\alpha)_{0<\alpha<\zeta}$ which converges to $x$ and $x_\alpha\in P_\zeta$ for all $0<\alpha<\zeta$. 

Further, if there is a $\zeta$-subsequence of $(x_\alpha)_{0<\alpha<\zeta}$, say for short it itself, such that $x_\alpha\in O_\zeta$,  then $x\in S'^{-1}\bar O_\zeta$.

If there is a $\zeta$-subsequence of $(x_\alpha)_{0<\alpha<\zeta}$, say for short it itself, such that $x_\alpha\notin O_\zeta$, then $\frac{1}{x}\in  O_\zeta$ for all $0<\alpha<\zeta$. Moreover $(\frac{1}{x_\alpha})_{0<\alpha<\zeta}$ is $\zeta$-fundamental sequence and hence, by Proposition 3 and Theorem 3, $x'=\lim\limits_{0<\alpha<\zeta}\frac{1}{x_\alpha}\in\bar O_\zeta\subseteq{\mathbb R}_\zeta$. One can see that $x'=\frac{1}{x}.$ Thus $x\in{\mathbb R}_\zeta$. $\Box$

{\bf Corollary 4.} The birthday of each $x\in{\mathbb R}_\zeta$ is less or equal to $\zeta$

{\bf Proof}. It is enough to consider $x\in{\mathbb R}_\zeta$ such that $x\notin\bar O_\zeta$. Thus, its birthday can be only great or equal to $\zeta$. Suppose that it is $\zeta+1$. Then $x=x'\pm\frac{1}{\zeta}$, if $x'\in{\mathbb R}_\zeta$ and the birthday of $x'$ is $\zeta$. If the birthday of $x'$ is less than $\zeta$, then $x=x'+\frac{2}{\zeta}$, or $x=x'-\frac{1}{2\zeta}$, or $x=x'-\frac{2}{\zeta}$, or $x=x'+\frac{1}{2\zeta}$, respectively, when $x'=\pm\frac{1}{\zeta}$. In all these cases $x\notin{\mathbb R}_\zeta$. Otherwise, $\pm\frac{1}{\zeta}$, $\pm\frac{2}{\zeta}$, $\pm\frac{1}{2\zeta}$ would be members of ${\mathbb R}_\zeta$ what is wrong. If $x=x'\pm\frac{1}{\zeta}$ and the birthday of $x'$ is $\zeta$ and $x'\notin{\mathbb R}_\zeta$. Then evidently $x\notin{\mathbb R}_\zeta$. The same argument for supposition that the birthday of $x$ is $\zeta+2$ and so on.

{\bf Corollary 5.} $\bar O_\zeta={\mathbb R}_\zeta$.

{\bf Proof}. Indeed,
by Corollary 4, each $x\in{\mathbb R}_\zeta$ whose birthday is $\zeta$ has a birthday form, i.e., $x=\{x^L\,|\,x^R\}$ and the birthdays of all $x^L$, $x^R$ are less than $\zeta$. Thus, all $x^L$, $x^R$ are members of $O_\zeta$. $\Box$

{\bf Lemma 10}. For every number $x\in{\mathbb R}_\zeta$, $\zeta=\omega^{\omega^\mu}$, $0\leq\mu<\Omega$ whose  birthday is $\alpha=\omega^\mu$, the following relation occurs: $\omega^x\notin{\mathbb R}_\zeta$.

{\bf Proof}. Suppose the opposite, i.e., $\omega^x\in{\mathbb R}_\zeta$. Denote this number by $y$
and consider the form of it as $\{0,y^L\,|\,y^R\}$, where $\{y^L\}$ and $\{y^R\}$ are sets of all positive numbers such that $y^L\in O_\zeta$ and $y^R\in O_\zeta$ with the following properties $y^L<y$ and $y^R>y$.

Denote by $\omega^{z^L}$ and $\omega^{z^R}$ the numbers which commensurate with $y^L$ and $y^R$, respectively.

Since the birthdays of $y^L$ and $y^R$ are less than $\zeta$ then, by Lemma 4, $\omega^{z^L}$ and $\omega^{z^R}$ are members of ${\mathbb R}_\zeta$ and their birthdays are also less than $\zeta$ because the zero-terms $\omega^{z^L_0}r_0$ and $\omega^{z^R_0}s_0$ in normal forms of numbers $y^L$ and $y^R$, respectively, are in $M_{\gamma_L}$ and $M_{\gamma_R}$, respectively, where $\gamma_L$ and $\gamma_R$ are the birthdays of $y^L$ and $y^R$, respectively, and thus less than $\zeta$. Consequently, every $\omega^{z^L}$ and $\omega^{z^R}$ are numbers of ${\mathbb R}_\zeta$. Hence, by Corollary 3 to Lemma 5, no $\omega^{z^L}$ and $\omega^{z^R}$ commensurate with $y=\omega^x$.

By the proof of Theorem 3, $y=\{0,y^L,r\cdot\omega^{z^L}\,|\,y^R,r\cdot\omega^{z^R}\}=\omega^z$, where $z=\{0,z^L\,|\,z^R\}$ and $r$ denotes a variable ranging over all positive reals. 

Since ${\mathbb R}_\zeta$ is $\zeta$-complete we can  choose an increasing $\zeta$-sequence $(y^L_\alpha)_{0<\alpha<\zeta}$ in the set $\{y^L\}$ such that $\lim\limits_{0<\alpha<\zeta}y^L_\alpha=y$.

Thus $\lim\limits_{0<\alpha<\zeta}y^L_\alpha=y=\omega^z$ then for each positive number $\varepsilon\in{\mathbb R}_\zeta$ there is an ordinal number $\alpha_0<\zeta$ such that $|y_\alpha-\omega^z|<\varepsilon$ for all $\alpha_0\leq\alpha<\zeta$. 

For each $y_\alpha$ there is a natural number $n_\alpha$ such that $y_\alpha<n_\alpha\cdot\omega^{z^L_\alpha}$ because $y_\alpha$ is commensurate with $\omega^{z_\alpha^L}$. Then $|y_\alpha-\omega^z|<\varepsilon$ implies $|n_\alpha\cdot\omega^{z_\alpha^L}-\omega^z|<\varepsilon$ for all $\alpha_0\leq\alpha<\zeta$. The latter inequality implies $|n_\alpha\cdot\omega^{z^L_\alpha-z}-1|<\frac{\varepsilon}{\omega^z}$ and thus there is a stronger inequality $|\omega^{z^L_\alpha-z}-1|<\frac{\varepsilon}{\omega^z}$. But by Lemma 3, $\omega^{z^L_\alpha}<\frac{1}{n}$ for all natural $0<n<\omega$ because $z^L_\alpha-z<0$. Thus $|\omega^{z^L_\alpha-z}-1|>|\frac{1}{n}-1|$ for all naturals $0<n<\omega$. In particular, for $n=2$ we obtain $|\omega^{z^L_\alpha-z}-1|>\frac{1}{2}$. Then for $\varepsilon=\frac{\omega^y}{2}$ the above inequality $|\omega^{z^L_\alpha-z}-1|<\frac{\varepsilon}{\omega^z}$ is in contradiction with $|\omega^{z^L_\alpha-z}-1|>\frac{1}{2}$. What proves that $x$ is commensurate with some of $\omega^{z^L}$ or $\omega^{z^R}$ and also $y\not=\omega^z$ for any $z$. 

{\bf Lemma 11}. Let $\zeta$ be the birthday of number $x\in{\mathbb R}_\zeta$. Then $x$  is reducible in the sense of Conway, i.e., $x=\sum\limits_{\beta<\zeta}\omega^{y_\beta}r_\beta$ and all $\beta$-terms $\omega^{y_\beta}\cdot r_\beta$ of $x$ are simpler numbers (whose birthdays are less than $\zeta$).

{\bf Proof}. Let $x$ be a positive number of $x\in{\mathbb R}_\zeta$ whose birthday is equal to $\zeta$. Since ${\mathbb R}_\zeta$ is $\zeta$-complete we can consider $x=\{L\,|\,R\}$, where $L=\{x^L\in{O_\zeta}: 0\leq x^L<x\}$ and $R=\{x^R\in{O_\zeta}: x< x^R\}$, i.e., $x^L,x^R$ are nonnegative \grqq dyadic\grqq\, numbers whose birthdays are less that $\zeta$. Moreover, we can choose an increasing $\zeta$-sequence $(x_\alpha)^L_{0<\alpha<\zeta}$ in the set $\{x^L\}$ such that $\lim\limits_{0<\alpha<\zeta}x^L_\alpha=x$.

 By Theorem 3, each $x^L$ is commesurate with some $\omega^{y^L}$ (say) and each $x^R$ is commensurate with some $\omega^{y^R}$. Moreover, the birthdays of $\omega^{x^L}$ and $\omega^{x^R}$ are less than $\zeta$ because the zero-terms $\omega^{y_0}r_0$ and $\omega^{z_0}s_0$ in normal forms of numbers $x^L$ and $x^R$, respectively, are in $M_{\gamma_L}$ and $M_{\gamma_R}$, respectively, where $\gamma_L$ and $\gamma_R$ are the birthdays of $x^L$ and $x^R$, respectively, and thus less than $\zeta$. Consequently, every $\omega^{y^L}$ and $\omega^{y^R}$ are numbers of ${\mathbb R}_\zeta$. 

We shall prove that $x$ must be commensurate with some $\omega^{y^L}$ or $\omega^{y^R}$. Otherwise, by the proof of Theorem 3, $x=\{0,x^L,r\cdot\omega^{y^L}\,|\,x^R,r\cdot\omega^{y^R}\}=\omega^y$, where $y=\{0,y^L\,|\,y^R\}$ and $r$ denotes a variable ranging over all positive reals. 

Since $\lim\limits_{0<\alpha<}x^L_\alpha=x=\omega^y$ then for each positive number $\varepsilon\in{\mathbb R}_\zeta$ there is an ordinal number $\alpha_0<\zeta$ such that $|x_\alpha-\omega^y|<\varepsilon$ for all $\alpha_0\leq\alpha<\zeta$. 

For each $x_\alpha$ there is a natural number $n_\alpha$ such that $x_\alpha<n_\alpha\cdot\omega^{y^L_\alpha}$ because $x_\alpha$ is commensurate with $\omega^{y_\alpha^L}$. Then $|x_\alpha-\omega^y|<\varepsilon$ implies $|n_\alpha\cdot\omega^{y_\alpha^L}-\omega^y|<\varepsilon$ for all $\alpha_0\leq\alpha<\zeta$. The latter inequality implies $|n_\alpha\cdot\omega^{y^L_\alpha-y}-1|<\frac{\varepsilon}{\omega^y}$ and thus there is a stronger inequality $|\omega^{y^L_\alpha-y}-1|<\frac{\varepsilon}{\omega^y}$. But by Lemma 3, $\omega^{y^L_\alpha}<\frac{1}{n}$ for all natural $0<n<\omega$ because $y^L_\alpha-y<0$. Thus $|\omega^{y^L_\alpha-y}-1|>|\frac{1}{n}-1|$ for all naturals $0<n<\omega$. In particular, for $n=2$ we obtain $|\omega^{y^L_\alpha-y}-1|>\frac{1}{2}$. Then for $\varepsilon=\frac{\omega^y}{2}$ the above inequality $|\omega^{y^L_\alpha-y}-1|<\frac{\varepsilon}{\omega^y}$ is in contradiction with $|\omega^{y^L_\alpha-y}-1|>\frac{1}{2}$. What proves that $x$ is commensurate with some of $\omega^{y^L}$ or $\omega^{y^R}$ and also $x\not=\omega^y$ for any $y$.

{\bf Proof of Theorem 7}.

When $\mu=0$  it is a classical result that every positive real number $x\in{\mathbb R}$ has an $n$th root. We recall the proof to show what new problems arise in the cases when $\mu>0$.

Let $\mu=0$. Then $P_\omega={\mathbb Q}$ and ${\mathbb R}_\omega={\mathbb R}.$ This is a well-known result of Classical Analysis. Consider the section $A|B$ of rationals ${\mathbb Q}$, where $A$ contains all negative numbers, zero number and those positive numbers $a$ such that $a^n<x$, and $B$ is the set of all rational numbers $b$ such that $a^n\geq x$. If there is a rational number $b$ such that $b^n=x$ we put $y=b$ and Theorem 7 is proved.

 If  not, then $B=\{b\,\,|\,\,b^n>x\}$. Clearly, that $A$ and $B$ are not empty because for each natural number $l$ such that $\frac{1}{l}<x<l$ one has $\frac{1}{l^n}<x<l^n$, thus, $\frac{1}{l}\in A$ and $l\in B$. It is also clear that $A\cap B=\emptyset$, $A\cup B={\mathbb Q}$ and, for each $a\in A$, $b\in B$, one has $a<b$ because otherwise, $b\leq a$ implies $b^n\leq a^n<x$ what is in contradityion with $b^n>x$. So the section $A|B$ defines a real number $y$. One can show that $y^n=x$, i.e., $y=\sqrt[n]{x}$.

Indeed, since for every $0<a<y<b$ one has $a^n<y^n<b^n$ and, for every natural number $l>0$, there are $a\in A$, $b\in B$ such $b-a<\frac{1}{l}$ (choose $y-\frac{1}{2l}<a<y$ and $y<b<y+\frac{1}{2l}$ what is possible because ${\mathbb Q}$ is dense in ${\mathbb R}$) we conclude that $b^n-a^n=(b-a)(b^{n-1}+b^{n-2}a+..+a^{n-1})<\frac{1}{l}\cdot n\cdot b^{n-1}_0$, where $b_0$ is a fixed number such that $b<b_0$, when $b-a<\frac{1}{l}$, and hence $A'|B'$ is a section of ${\mathbb Q}$, where $B'=\{b^n\,\,|\,\, b\in B\}$ and $A'={\mathbb Q}\setminus B'$ which defines number $x$. Thus, $y^n=x$. 

Suppose now that $\zeta=\omega^{\omega^\mu}$, $0<\mu<\Omega$ and consider the similar section $A|B$ of $P_\zeta$, then it is not evident that for every ordinal number $\alpha<\zeta$ there are $a\in A$, $b\in B$ such $b-a<\frac{1}{\alpha}$ in spite of the fact that they exist because a priori there are many other sections of $P_\zeta$ which do not satisfy this property. Indeed, consider the set $B$ of all numbers $b\in P_\zeta$ such that $b>\frac{1}{l}$, for each natural number $l>0$, and $A=P_\zeta\setminus B$. Clearly, $A|B$ is a section of $P_\zeta$ but for the Conway's number $c=\{A|B\}$ the number $c+\frac{1}{\omega}$ satisfied that following inequalities: $a<c<c+\frac{1}{\omega}<b$, for every $a\in A$, $b\in B$, and hence $b-a>c+\frac{1}{\omega}-c=\frac{1}{\omega}$, although $\omega<\zeta$. 
 
 Thus, we have to show that it is not the case. There are several ways to do it. We use Conway's  Theorem 21 in $\cite{l3}$, p. 33, that each number $y\in{\bf No}$ defines a unique expression
\begin{equation}
\label{f0201}
y=\sum\limits_{0\leq\beta<\alpha}\omega^{y_\beta}r_\beta,
\end{equation}
in which $\alpha$ denotes some ordinal, the numbers $r_\beta$ $(0\leq\beta<\alpha)$ are non-zero reals, and the numbers $y_\beta$ form a descending sequence of numbers. Moreover, normal forms for distinct $y$ are distinct, and every form satisfying these conditios occurs.

There are no any hint in  Conway's proof that $y=\sqrt[n]{x}$  in $(\ref{f0201})$ is a number in ${\mathbb R}_\zeta$ when $x\in{\mathbb R}_\zeta$. We shall show that $y=\sqrt[n]{x}\in{\mathbb R}_\zeta$ indeed.

Prove by contradiction. Suppose that it is not so and $y=\sqrt[n]{x}\notin{\mathbb R}_\zeta$. 

Notice first that $\omega^{y_0}r_0<\alpha_0$, for some ordinal $0<\alpha_0<\zeta$. Otherwise, if $\omega^{y_0}r_0>\alpha$ for all $0<\alpha<\zeta$, and hence $y>\alpha$ for all $0<\alpha<\zeta$, because $y_1=\sum\limits_{0<\beta<\alpha}\omega^{y_\beta}r_\beta$ is small compared to $\omega^{y_0}r_0$. Thus, by Lemma 2, $y^n>\alpha$ for all $0<\alpha<\zeta$ what is in contradiction with $y^n=x\in{\mathbb R}_\zeta$.
 
 Moreover, $\omega^{y_0}r_0\in{\mathbb R}_\zeta$. Otherwise, $\omega^{y_0}r_0<\frac{1}{\alpha}$ for all $0<\alpha<\zeta$ and hence $y<\frac{1}{\alpha}$ for all $0<\alpha<\zeta$, because $y_1=\sum\limits_{0<\beta<\alpha}\omega^{y_\beta}r_\beta$ is small compared to $\omega^{y_0}r_0$. Thus, by Lemma 3, $y^n=x<\frac{1}{\alpha}$ for all $0<\alpha<\zeta$ what is in contradiction with $y^n=x\in{\mathbb R}_\zeta$.
 
Then there is a $\beta$-term in $(\ref{f0201})$ such that $\omega^{y_{\beta}}r_{\beta}\notin{\mathbb R}_\zeta$, otherwise, it would be contrary to our assumption. Moreover, there is the smallest ordinal number $\beta_0$ with such property because the set of all ordinals is well ordered and hence all $\beta$ with above property has the smallest one. Thus,
\begin{equation}
\label{33334}
y=\sqrt[n]{x}=\sum\limits_{0\leq\beta<\beta_0}\omega^{y_{\beta}}r_{\beta}+\sum\limits_{\beta_0\leq\beta<\alpha}\omega^{y_{\beta}}r_{\beta}=a+b, 
\end{equation}
 where $a\in{\mathbb R}_\zeta$, and $\sum\limits_{0\leq\beta<\beta_0}\omega^{y_{\beta}}r_{\beta}$ is its normal form, $b\notin{\mathbb R}_\zeta$, and $\sum\limits_{\beta_0\leq\beta<\alpha}\omega^{y_{\beta}}r_{\beta}$ is its normal form. Clearly, that $|b|$ is small compared to any positive number in ${\mathbb R}_\zeta$, i.e., for each number $c>0$ and $c\in{\mathbb R}_\zeta$ one has $c\cdot b\notin{\mathbb R}_\zeta$. Otherwise, it should be contrary to the minimality of $\beta_0$.
 
 Consider $x=(a+b)^n=a^n+C_n^1a^{n-1}b+C^2_na^{n-2}b^2+...+b^n$. Clearly, $x-a^n\in{\mathbb R}_\zeta$ and $C_n^1a^{n-1}b+...+b^n$ is small compared to any number in ${\mathbb R}_\zeta$ and cannot be equal to $x-a^n\in{\mathbb R}_\zeta$. The only way to avoid contradiction is to conclude that $b=0$ and $a=y$. Thus $y\in{\mathbb R}\zeta$ is an $n$th root of $x$. $\Box$

{\bf Corollary 6} {\it Let $x$ be a positive number in $\mathbb R_\zeta$, $\zeta=\omega^{\omega^\mu}$, $0<\mu<\Omega$, $A=\{a\in P_\zeta\,|\, a\leq 0\,\vee\,a^n<x\}$, $B=\{b\in P_\zeta\,|\, b> 0\,\&\,b^n>x\}$ and $A'=\{a'\in P_\zeta\,|\, a'\leq 0\,\vee\,s'=a^n,\,\,a>0\,\&\,a\in A\}$, $B'=\{b'\in P_\zeta\,|\, b'=b^n,\,\,b> 0\,\&\,b\in B\}$. Then  $\{A\,|\,B\}$ defined $y=\sqrt[n]{x}$ and $\{A'\,|\,B'\}$ defined $x$ such that $y^n=x$. }

{\bf Proof}. By Theorem 7, $y=\sqrt[n]{x}\in{\mathbb R}_\zeta$. Consider the following $\zeta$-sequences  $(y-\frac{1}{2\alpha})_{0<\alpha<\zeta}$ and $(y+\frac{1}{2\alpha})_{0<\alpha<\zeta}$. Since $P_\zeta$ is dense in ${\mathbb R}_\zeta$, for each $0<\alpha<\zeta$, we can find $a_\alpha\in A$ and $b_\alpha\in B$ such that $y-\frac{1}{2\alpha}<a_\alpha<b+\frac{1}{2\alpha}$. Moreover, in addition we can choose $a_\alpha$ and $b_\alpha$ such that $a_\alpha<a_{\alpha+1}$ and $b_\alpha>b_{\alpha+1}$. Then $[a_1,b_1]\supset[a_1,b_1]\supset...\supset[a_\alpha,b_\alpha]\supset[a_{\alpha+1},b_{\alpha+1}]\supset...$ is an transfinite $\zeta$-sequence of embedded intervals $[a_\alpha,b_\alpha]$ in ${\mathbb R}_\zeta$ such that  $\lim\limits_{0<\alpha<\zeta}(b_\alpha-a_\alpha)=0$, because $b_\alpha-a_\alpha<y+\frac{1}{2\alpha}- y+\frac{1}{2\alpha}=\frac{1}{\alpha}$. Then there is a unique number $c\in {\mathbb R}_\zeta$ which belongs to each interval $[a_\alpha,b_\alpha]$, $0<\alpha<\zeta$, and thus $c=y$.

Consider now $\zeta$-sequences $(x-\frac{1}{2\alpha})_{0<\alpha<\zeta}$ and $(x+\frac{1}{2\alpha})_{0<\alpha<\zeta}$. Since $P_\zeta$ is dense in ${\mathbb R}_\zeta$, for each $0<\alpha<\zeta$, we can find $a'_\alpha\in A'$ and $b'_\alpha\in B'$ such that $x-\frac{1}{2\alpha}<a'_\alpha<b'+\frac{1}{2\alpha}$. Indeed, by Theorem 7, there are numbers $v_\alpha,w_\alpha\in {\mathbb R}_\zeta$ such that $v^n=x-\frac{1}{2\alpha}$ and $w^n=x+\frac{1}{2\alpha}$. By density of $P_\zeta$ in ${\mathbb R}_\zeta$, we can find numbers $a_\alpha,b_\alpha\in P_\zeta$ such thar $v_\alpha<a_\alpha<y<b_\alpha<w_\alpha$. Moreover, in addition we can choose $a_\alpha$ and $b_\alpha$ such that $a_\alpha<a_{\alpha+1}$ and $b_\alpha>b_{\alpha+1}$. Then we put $a'_\alpha=a^n_\alpha$ and $b'_\alpha=b^n_\alpha$ and see that $[a'_1,b'_1]\supset[a'_2,b'_2]\supset...\supset[a'_\alpha,b'_\alpha]\supset[a'_{\alpha+1},b'_{\alpha+1}]\supset...$ is an transfinite $\zeta$-sequence of embedded intervals $[a'_\alpha,b'_\alpha]$ in ${\mathbb R}_\zeta$ such that  $\lim\limits_{0<\alpha<\zeta}(b'_\alpha-a'_\alpha)=0$, because $b'_\alpha-a'_\alpha<x+\frac{1}{2\alpha}- x+\frac{1}{2\alpha}=\frac{1}{\alpha}$. Then, by Lemma 5, there is a unique number $c'\in {\mathbb R}_\zeta$ which belongs to each interval $[a'_\alpha,b'_\alpha]$, $0<\alpha<\zeta$,  and thus $c'=x$. $\Box$

\begin{center}
{\bf 9. A root in ${\mathbb R}_\zeta$ of odd-degree polynomial with coefficients in ${\mathbb R}_\zeta$}
\end{center}

{\bf Theorem 9}. {\it Every odd-degree polynomial with coefficients in ${\mathbb R}_\zeta$ has a root in ${\mathbb R}_\zeta$.}
 
{\bf Proof}. By Theorem 5, proved by Conway in $\cite{l3}$, p. 40-41, every odd-degree polynomial with coefficients in ${\bf No}$ has a root in ${\bf No}$. We have to show that by supposition of Theorem 9, i.e., when an odd-degree polynomial, say without loss of generality $P_n(x)=a_0x^n+a_1x^{n-1}+a_{2}x^{n-2}+...+a_{n-1}x+a_n$, $n=2k+1$, $0<k<\omega$, with coefficients $a_0,a_1,a_{2},...,a_{n-1},a_n\in{\mathbb R}_\zeta$, $a_0\not=0$, has a root, say $\bar x$, in ${\mathbb R}_\zeta$. For now we know only, by Theorem 4, that $\bar x\in{\bf No}$.

We omit trivial case $k=0$ when $f_1(x)=a_0x+a_1$ because it has an evident root $\bar x=-\frac{a_1}{a_0}\in{\mathbb R}_\zeta$. Thus, we consider only cases $0<k<\omega$.

First of all, clearly $a_n\not=0$, and we shall show that there are no inequalities $|\bar x|>\alpha$, for all $0<\alpha<\zeta$, and $|\bar x|<\frac{1}{\alpha}$, for all $0<\alpha<\zeta$, and hence there are numbers $a,b\in{\mathbb R}_\zeta$ such that $a\leq\bar x\leq b$.

Indeed, consider the {\it first case} when $|\bar x|>\alpha$, for all $0<\alpha<\zeta$. Then 
\begin{equation}
\label{f0309}
|\bar x^n+a_1\bar x^{n-1}+a_2\bar x^{n-2}+...+a_{n-1}\bar x|>|-a_n|.
\end{equation}
Contradiction, because $\bar x^n+a_1\bar x^{n-1}+a_2\bar x^{n-2}+...+a_{n-1}\bar x+a_n=0$ implies $|\bar x^n+a_1\bar x^{n-1}+a_2\bar x^{n-2}+...+a_{n-1}\bar x|=|-a_n|$.

Formula $(\ref{f0309})$ follows from  Lemmas 3 and 4. Indeed, in the first case when $|\bar x|>\alpha$ for all $0<\alpha<\zeta$ Lemma 2 implies formula $(\ref{f0309})$ because $|\bar x^n+a_1\bar x^{n-1}+a_2\bar x^{n-2}+...+a_{n-1}\bar x|=|\bar x|\cdot|\bar x^{n-1}+a_1\bar x^{n-2}+a_2\bar x^{n-3}+...+a_{n-1}|>|\bar x\cdot 1|>|-a_n|.$

Consider now the {\it second case} when 
$|\bar x|<\frac{1}{\alpha}$ for all $0<\alpha<\zeta$.

Then we shall prove that 
\begin{equation}
\label{f1309}
|\bar x^n+a_1\bar x^{n-1}+a_2\bar x^{n-2}+...+a_{n-1}\bar x|<|-a_n|.
\end{equation}
Contradiction because as above it should be an equality.

We continue now the proving of Theorem 9. We have already proved that there are two numbres $a,b\in{\mathbb R}_\zeta$ such that $a\leq\bar x\leq b$. Really, by above provings, $|\bar x|\leq\alpha_0$ or $|\bar x|\geq\frac{1}{\alpha_0}$ for some ordinal number $\alpha_0\in{\mathbb R}_\zeta$. In the first case $a=-\alpha_0$ and $b=\alpha_0$ and in the second case $a=\frac{1}{\alpha_0}$ and $b=\alpha_1$ for some ordinal $\alpha_1>\bar x$, if of course $\bar x\geq\frac{1}{\alpha_0}$; and $a=-\alpha_1$ for some ordinal $\alpha_1>|\bar x|$ and $b=-\frac{1}{\alpha_0}$ if of course $\bar x\leq-\frac{1}{\alpha_0}$.

Let $P_n(x)=a_0x^n+a_1x^{n-1}+...+a_{n-1}x+a_n$, $a_0\not=0$, be a polynomial function such that $n=2k+1$, $0<k<\omega$ and $a_n\not=0$ because if $a_n=0$, then  the polynomial $P_n(x)$ has an evident root $\bar x=0\in{\mathbb R}_\zeta$.

 By the linear transformation $\varphi:[0,1]\rightarrow [a,b]$, given by formula $\varphi(x)=(b-a)x+a$ we can consider another polynomial function $Q_n(x)=P_n(\varphi(x))=b_0x^n+b_1x^{n-1}+...+b_{n-1}x+b_n$, $b_0=a_0\cdot(b-a)\not=0$, and $b_n=a_0a^n+a_1a^{n-1}+...+a_{n-1}a+a_n\not=0$ because $a$ is not a root of $P_n(x)$, otherwise. Theorem 9  had be proved.
 
 Since $\varphi:[0,1]\rightarrow [a,b]$ can be also considered in ${\bf No}$ is a bijection, because the inverse map $\varphi^{-1}:[a,b]\rightarrow [0,1]$ is given by $\varphi^{-1}(y)=\frac{1}{b-a}y+\frac{a}{a-b}$ which preserves the oders on $[0,1]$ and $[a,b]$, i.e., if $x_1<x_2$, then evidently $\varphi(x_1)=x_1(b-a)+a<x_2(b-a)+a=\varphi(x_2)$ as well as $y_1<y_2$ implies $\varphi^{-1}(y_1)=\frac{1}{b-a}y_1+\frac{a}{a-b}<\frac{1}{b-a}y_2+\frac{a}{a-b}=\varphi^{-1}(y_2)$. Consequently, $a\leq\bar x\leq b$ implies $0\leq \hat x\leq 1$, where $\hat x=\varphi^{-1}(\bar x)$.
 
 One can easily check that $\hat x$ is a root of $g_n(x)$ because $Q_n(\hat x)=P_n(\varphi(\hat x))=P_n(\bar x)=0$ and if we prove that $\hat x\in{\mathbb R}_\zeta$, then evidently $\bar x=\varphi(\hat x)\in{\mathbb R}_\zeta$ and Theorem 9 will be proved.
 
We shall show that $\hat x=c+\varepsilon$, where $c\in[0,1]\subset{\mathbb R}_\zeta$ and $\varepsilon\notin[0,1]$ and smaller than each number in $(0,1]$.

Indeed, consider a normal form of $\hat x$:
\begin{equation}
\label{f13333}
\hat x=\sum\limits_{0\leq\beta<\alpha}\omega^{y_{\beta}}r_{\beta}. 
\end{equation}

 Suppose that $\hat x\notin[0,1]\subset{\mathbb R}_\zeta$. Then there is a $\beta$-term in $(\ref{f13333})$ such that $\omega^{y_{\beta}}r_{\beta}\notin[0,1]\subset{\mathbb R}_\zeta$, otherwise, it would be contrary to our assumption. Moreover, there is the smallest ordinal number $\beta_0$ with such property because the set of all ordinals is well ordered and hence all $\beta$ with above property has the smallest one. Thus,
\begin{equation}
\label{13334}
\hat x=\sum\limits_{0\leq\beta<\beta_0}\omega^{y_{\beta}}r_{\beta}+\sum\limits_{\beta_0\leq\beta<\alpha}\omega^{y_{\beta}}r_{\beta}=c+\varepsilon, 
\end{equation}
 where $c\in[0,1]\subset{\mathbb R}_\zeta$, and $\sum\limits_{0\leq\beta<\beta_0}\omega^{y_{\beta}}r_{\beta}$ is its normal form, $\varepsilon\notin[0,1]\subset{\mathbb R}_\zeta$, and $\sum\limits_{\beta_0\leq\beta<\alpha}\omega^{y_{\beta}}r_{\beta}$ is its normal form. 
 Clearly, that $|\varepsilon|$ is small compared to any number in $(0,]\subset{\mathbb R}_\zeta$, i.e., $|\varepsilon|<\frac{1}{\alpha}$ for all $0<\alpha<\zeta$.

 We can suppose that $0<\hat x<b$, otherwise, $0$ or $1$ will be a root of $g_n(x)$ and $a$ or $b$
will be a root of $P_n(x)$ and also Theorem 7 will be proved.

By Lemma 2, $\hat x=c+\varepsilon$, where $c\in[0,1]\subset{\mathbb R}_\zeta$ and $|\varepsilon|<\frac{1}{\alpha}$ for all $0<\alpha<\zeta$, moreover, $c\not=0$, because for $c=0$, by Lemma 2, $|b_0\varepsilon^n+b_1\varepsilon^{n-1}+...+b_{n-1}\varepsilon|<|-b_n|$.

By the same Lemma, $|b_0\varepsilon^n+b_0\varepsilon^{n-1}c+...+nb_0\varepsilon c^n+b_1\varepsilon^{n-1}+b_1(n-1)\varepsilon^{n-2}c+...+b_1(n-1)\varepsilon^{n-1})+...+b_{n-1}\varepsilon|<|-b_0c^n-b_1c^{n-1}-..-b_{n-1}c-b_n|$.

Thus, $\hat x=c\in{\mathbb R}_\zeta$ and hence $\bar x=\varphi(\hat x)=\varphi(c)=\bar x\in{\mathbb R}_\zeta$ and Theorem 7 has been finally proved. 

Consider now the {\it second case} when 
$|\bar x|<\frac{1}{\alpha}$ for all $0<\alpha<\zeta$.

Then we shall prove that 
\begin{equation}
\label{f1309}
|\bar x^n+a_1\bar x^{n-1}+a_2\bar x^{n-2}+...+a_{n-1}\bar x|<|-a_n|.
\end{equation}
Contradiction because as above it should be an equality.
Lemma 3 implies formula $(\ref{f1309})$. Indeed, if contrary
\begin{equation}
\label{f1319}
|\bar x^n+a_1\bar x^{n-1}+a_2\bar x^{n-2}+...+a_{n-1}\bar x|\geq|-a_n|,
\end{equation}
then take any ordinal number $\alpha_0$ such that $\frac{1}{\alpha_0}<|-a_n|$ and obtain an inequality $|\bar x^n+a_1\bar x^{n-1}+a_2\bar x^{n-2}+...+a_{n-1}\bar x|>\frac{1}{\alpha_0}$ what is in contradiction with Lemma 3. $\Box$

\begin{center}
{\bf 9. Definition of exponential and logarithmic functions}
\end{center}

Conway did not know of any power operation $x^y$ defined for all numbers $y$ and for all positive $x$ but he defined the exponential function $exp\,\, x$ provided $-n<x<n$ for some integer $n$ and also he defined a natural function $\omega^x$ which play a vital role in his theory of numbers. Moreover, he expressed his doubt about the possibility of giving a natural definition of the function $x^y$ for infinite $y$ and said that \grqq Nor does there seem to be any particular point of making these definition\grqq (see $\cite{l3}$, p. 43). In this paragraph we dispel Conway's doubts.

First of all, here are two definitions  and some remarks we need for.

{\bf Definition 15}. A Conway's number $x\in{\bf No}$ is called {\it finite} if for some natural number $n$ there are inequalities $-n<x<n$. Otherwise, $x\in{\bf No}$ is called {\it infinite}. The set of all infinite numbers in ${\mathbb R}_\zeta\subset{\bf No}$ we denote by ${\mathbb R}^{inf}_\zeta$ and the set of all finite numbers in ${\mathbb R}_\zeta\subset{\bf No}$ we denote by ${\mathbb R}^{f}_\zeta$. In our notation ${\mathbb R}^{inf}_\zeta=\{x\in{\bf No}\,|\,x<-\infty_\omega \vee x>+\infty_\omega\}$ and ${\mathbb R}^{f}_\zeta=\{x\in{\bf No}\,|\,-\infty_\omega<x<+\infty_\omega\}$.

{\bf Definition 16}.  By an {\it infinite initial integer}    we understand a number $x\in{\bf No}$ whose normal form is the following: $x=\sum\limits_{0\leq\alpha'<\alpha}\omega^{y_{\alpha'}}r_{\alpha'}$, where $y_{\alpha'}>0$, for all $0\leq\alpha'<\alpha$. We also call  null $0$ the only initial finite integer number.

 One can see that each number $x\in{\bf No}$ is a sum of an infinite initial integer $x'$ or the finite initial integer $x'=0$ and a finite number $x''$. Indeed, let $\sum\limits_{0\leq\beta<\alpha}\omega^{y_\beta}r_{\beta}$ be a normal form of $x$. Then we put $x'=\sum\limits_{0\leq\beta<\alpha'}\omega^{y_\beta}r_{\beta}$, where $y_\beta>0$ for all ${0\leq\beta<\alpha'}$ in the normal form of $x$ and $\alpha'$ is the greatest ordinal with such property; if $\alpha'=0$, then we put $x'=0$. Notice that,  by Conway's definition, the formal sum $\sum\limits_{0\leq\beta<\alpha'}\omega^{y_\beta}r_{\beta}$   is the {\it simplest} number (i.e., the first-born number) whose $\beta$-term is $\omega^{y_\beta}r_\beta$ for all $0\leq\beta<\alpha'$. We also put  $x''=\sum\limits_{\alpha'\leq\beta<\alpha}\omega^{y_\beta}r_{\beta}$ in the normal form of $x$ and it is the  simplest number  whose $\beta$-term is $\omega^{y_\beta}r_\beta$ for all $\alpha'\leq\beta<\alpha$ as well. 
 
 Then it is obvious that $x=x'+x''$ and we call $x'$ an infinite part of $x$ and $x''$ a finite part of $x$.  Clearly, finite numbers $x\in{\mathbb R}_\zeta$, $\zeta=\omega^{\omega^\mu}$, $0\leq\mu<\Omega$, have no infinite parts and infinite inegers have no finite parts. Moreover, all these partial sums as well as $x$ itself must belong to the set $M_\gamma$, where $\gamma$ is the birthday of $x$, and hence if $x\in{\mathbb R}_\zeta$, $\zeta=\omega^{\omega^\mu}$, $\mu>0$, then its infinite part $x'$ and finite part $x-x'$ are also in ${\mathbb R}_\zeta$. 
 
Denote by ${\bf Z'}$ the set of all initial intergers of ${\bf No}$, i.e., ${\bf No}'=\{x'\,|\,x\in{\bf No}\,\&\, x'=0\}$ and by ${\bf No}''=\{x''\,|\, x\in{\bf No}\}$. Note that ${\bf No}''={\bf No}^{f}$ but ${\bf No}''\not={\bf No}^{\inf}$ at all. It is clear that ${\bf Z}'$ and ${\bf No}''$ are additive subgroups of ${\bf No}$.

{\bf Definition 17.} We say that numbers $x$ and $y$ are congruent if $x-y\in {\bf No}''$. 

Plainly this is an equivalent relation whose equivalent Classes are {\it convex}, i.e., if $x<z<y$ and $x$ and $y$ are congruent, then $z$ is compatible with both. Indeed, this relation is, evidently, reflexive since 
$$x-x=0\in{\bf No}'';$$
it is also symmetrical because  
$$x-y\in{\bf No}''\Rightarrow y-x\in{\bf No}'';$$
it is also transitive because 
$$(x-y\in {\bf No}'')\,\&\,(y-z\in {\bf No}'')\Rightarrow (x-y+y-z\in {\bf No}'')\Rightarrow (x-z\in {\bf No}'').$$

It is really convex because if $x<z<y$ and
$x-y\in {\bf No}''$, then $x-z+z-y\in {\bf No}''$ and hence $x-z\in {\bf No}''$ and $z-y\in {\bf No}''$. If $x,z,y\in{\bf No}''$ it is clear because ${\bf No}''$ is a subgroup of ${\bf No}$. If $x,z,y\notin{\bf No}''$, then $z-x,y-z\in{\bf No}''$ because $z-x<y-x$ and $y-z<y-z$ since $x<z<y$.

Now we give a definition of  exponential functions and logarithmic functions with domains $X={\bf No}$ and $X\subset {\bf No}^+$, where ${\bf No}^+=\{x\in{\bf No}\,|\,x>0\}$, respectively,  and hence a definition of  expotential functions and logarithmic functions with domains $X_\zeta=(-\infty_\zeta,+\infty_\zeta)$ and $X_\zeta\subset (\frac{1}{+\infty_\zeta},+\infty_\zeta)$ of ${\mathbb R}_\zeta$, respectively,  for all $\zeta=\omega^{\omega^\mu}$, where $0<\mu<\Omega$. It can be compared

First, as mentioned above, Conway defined the expotential function  $y=e^x$ and logarithmic function $y=\ln x$ with domains $X_\omega=(-\infty_\omega,+\infty_\omega)$ and $X_\omega=(\frac{1}{+\infty_\omega},+\infty_\omega)$, respectively, by power-series which are convergent in these intervals (see $\cite{l3}$, p. 43) using his result that power-series with real coefficients is always absolutely convergent for all infinitesimal valuses of a variable (see Theorem 23 in $\cite{l3}$, p. 40). Scheme of this construction is the following. At the beginning $y=e^x$ and $\ln(1+x)$ are defined on the interval $(-\frac{1}{+\infty_\omega},+\frac{1}{+\infty_\omega})$, or for infinitesimals $\delta$, i.e., $-\frac{1}{n}<\delta<+\frac{1}{n}$ for all natural numbers $n>0$, by the following formulas:

\begin{equation}
\label{f5151}
e^x=1+x+\frac{x^2}{2!}+...+\frac{x^n}{n!}+...
\end{equation}
and
\begin{equation}
\label{f5152}
\ln (1+x)=x-\frac{x^2}{2}+...+(-1)^{n-1}\frac{x^n}{n}+...,
\end{equation}
respectively, which uniquely define Conway's numbers because of the fact that power-series with real coefficients is always absolutely convergent for all infinitesimal valuses of a variable.

Indeed, for each  infinitesimal $\delta$ formulars $\ref{f5151}$ and $\ref{f5152}$ define uniquely numbers $e^\delta-1$ and $\ln(1+\delta)$, respectively. Really, if $\delta>0$, then
\begin{equation}
\label{f5203}
e^\delta-1=\delta+\frac{\delta^2}{2!}+...+\frac{\delta^n}{n!}+...=\{\delta,\delta+\frac{\delta^2}{2!},...,\delta+\frac{\delta^2}{2!}+...+\frac{\delta^2}{2!},...\,|\,1,\frac{1}{2},...,\frac{1}{n},...\}.
\end{equation}
If $\delta<0$, then
\begin{equation}
\label{f5203}
e^\delta-1=\delta+\frac{\delta^2}{2!}+...+\frac{\delta^n}{n!}+...=\{-1,-\frac{1}{2},...,-\frac{1}{n},...|\,\delta,\delta+\frac{\delta^2}{2!},...,\delta+\frac{\delta^2}{2!}+...+\frac{\delta^2}{2!},...\}.
\end{equation}

Similarly, if $\delta>0$, then
\begin{equation}
\label{f4203}
\ln(1+\delta)=\delta-\frac{\delta^2}{2}+...+(-1)^{n-1}\frac{\delta^n}{n}+...=\{\delta,\delta-\frac{\delta^2}{2},...,\delta-\frac{\delta^2}{2}+...+(-1)^{n-1}\frac{\delta^n}{n},...\,|\,1,\frac{1}{2},...,\frac{1}{n},...\}.
\end{equation}
If $\delta<0$, then
\begin{equation}
\label{f4203}
\ln(1+\delta)=\delta-\frac{\delta^2}{2}+...+(-1)^{n-1}\frac{\delta^n}{n}+...=\{-1,-\frac{1}{2},...,-\frac{1}{n},...\,|\,\delta,\delta-\frac{\delta^2}{2},...,\delta-\frac{\delta^2}{2}+...+(-1)^{n-1}\frac{\delta^n}{n},...\}..
\end{equation}

For arbitrary $x\in{\bf No}$ such that $-\infty_\omega<x<+\infty_\omega$ and $\frac{1}{+\infty_\omega}<1+x<+\infty_\omega$, respectively, there are real numbers $\bar x$ such that $x=\bar x+\delta$, where $\delta$ is an infinitesimal $\delta$, and we put $e^x=e^{\bar x+\delta}=e^{\bar x}\cdot e^\delta$ and $\ln(1+x)=\ln(1+\bar x+\delta)=\ln((1+\bar x)(1+\frac{\delta}{1+\bar x}))=\ln(1+\bar x)+\ln(1+\frac{\delta}{1+\bar x})$, respectively. So there is a bijection between $(-\infty_\omega,+\infty_\omega)$ and $(\frac{1}{+\infty_\omega},+\infty_\omega)$, given by $y=e^x$ or by $\ln y$.

Actually there is much more broad extension of functions $y=a^x$ and $y=\log_ax$  for all numbers $x\in{\bf No}$ and for all $x\in X\subset{\bf No}^+$, respectively,  where $\frac{1}{+\infty_\omega}<a<+\infty_\omega$ such that $a\not= 1$, in particular, for all $x\in{\mathbb R}_\zeta$ and for all $x\in X_\zeta\subset{\mathbb R}^+_\zeta$, where $a\in (\frac{1}{+\infty_\omega},+\infty_\omega)\subset{\mathbb R}_\zeta$ such that $a\not= 1$.

{\bf Theorem 10.} {\it Let $a\in{\bf No}$ be a fixed  number such that $\frac{1}{+\infty_\omega}<a<+\infty_\omega$ and $a\not=1$. Then there are functions $y=a^x$  and $y=\log_ax$, defined for each $x\in{\bf No}$ and for each $x\in X\subset {\bf No}^+$, where $X$ is some definite subset  of ${\bf No}^+$, repectively, which are inverse with each other and when $a\in{\mathbb R}$, then the restrictions of these functions on ${\mathbb R}$ and ${\mathbb R}^+$ are the real-valued expotential function and logarithmic function, respectively.  }

{\bf Proof.}  As mentioned above, Conway defined the functions $y=e^x$ and $y=\ln x$ on $(-\infty_\omega,+\infty_\omega)\subset{\bf No}$ and $(\frac{1}{+\infty_\omega},+\infty_\omega)\subset{\bf No}^+$, respectively, which are extensions of real-valued functions $y=e^x$ and $y=\ln x$,  when $x\in{\mathbb R}$ and $x\in{\mathbb R}^+$, respectively.

 A light extension is for a fixed number $a$ such that $\frac{1}{+\infty_\omega}<a<+\infty_\omega$ and $a\not= 1$, when we put $y=a^x=(e^{\ln a})^x=e^{\ln a\cdot x}$ and $y=\log_ax=\frac{\ln x}{\ln a}$, which are well-defined above for all $x\in(-\infty_\omega,+\infty_\omega)\subset{\bf No}$ and all $x\in(\frac{1}{+\infty_\omega},+\infty_\omega)\subset{\bf No}^+$, respectively, and are inverse with each other; moreover, when 
 $a\in{\mathbb R}$, then the restrictions of these functions on ${\mathbb R}$ and ${\mathbb R}^+$ are the real-valued expotential function and logarithmic function, respectively.

We have to extend them to the fucntions $y=a^x$ and $y=\log_ax$ for all numbers $x\in{\bf No}$ and $x\in X\subset {\bf No}^+$, respectively, where $X=X^-\cup X^0\cup X^+$, $X^0=(\frac{1}{+\infty_\omega},+\infty_\omega)$ and $X^-\cup X^+$ will be defined further.

Suppose at first that $a>1$. At the beginning we define $y=a^x$ for every number $x\in{\bf No}$ such that $x>+\infty_\omega$. For $x=\omega$ we put $a^\omega\stackrel{def}{=}\omega$. For all other $x>+\infty_\omega$ we consider $x=x'+x''$, where $x'$ and $x''$  are infinite and finite part of $x$ and put, by definition, the value of the function $a^x$, given by the following formulas:
\begin{equation}
\label{f131}
a^x=a^{x'+x''}\stackrel{def}{=}a^{x'}\cdot a^{x''}\stackrel{def}{=}\omega^{\frac{x'}{\omega}}\cdot a^{x''},
\end{equation}
 where the latter expression $\omega^{\frac{x'}{\omega}}$ is uniquely defined by formula (\ref{f314}).

Clearly, if $x_1,x_2\in{\bf No}$ such that $x_1,x_2>+\infty_\omega$, then 
\begin{equation}
\label{f132}
\begin{array}{l}
a^{x_1+x_2}=a^{x'_1+x''_1+x'_2+x''_2}=a^{(x'_1+x'_2)+(x''_1+x''_2)}=\omega^{\frac{x'_1+x'_2}{\omega}}\cdot a^{x''_1+x''_2}=\\
=\omega^{\frac{x'_1}{\omega}+\frac{x'_2}{\omega}}\cdot a^{x''_1+x''_2}=\omega^{\frac{x'_1}{\omega}}\cdot\omega^{\frac{x''_2}{\omega}}\cdot a^{x''_1}\cdot a^{x''_2}=a^{x'_1}\cdot a^{x'_2}\cdot a^{x''_1}\cdot a^{x''_2}=\\
=a^{x'_1}\cdot a^{x''_1}\cdot a^{x'_2}\cdot a^{x''_2}=a^{x'_1+x''_1}\cdot a^{x'_2+x''_2}=a^{x_1}\cdot a^{x_2},
\end{array}
\end{equation}
 by a correct formula $\omega^{\frac{x'_1}{\omega}+\frac{x'_2}{\omega}}=\omega^{\frac{x'_1}{\omega}}\cdot\omega^{\frac{x'_2}{\omega}}$ (see Theorem 20 in $\cite{l3}$, Theorem 21, p. 33.)

If $x\in{\bf No}$ such that $x<-\infty_\omega$, then we put $a^x=\frac{1}{a^{-x}}$. It is clear that for all $x_1,x_2<-\infty_\omega$ we obtain $a^{x_1+x_2}= \frac{1}{a^{-x_1-x_2}}=\frac{1}{a^{-x_1}\cdot a^{-x_2}}=\frac{1}{a^{-x_1}}\cdot\frac{1}{a^{-x_2}}=a^{x_1}\cdot a^{x_2}$.

At last, if $a<1$, then for each $x\in{\bf No}$ we put $a^x=\frac{1}{(\frac{1}{a})^x}$. Indeed, $a^0=\frac{1}{\frac{1}{a}^0}=1$, $a^{-1}=\frac{1}{(\frac{1}{a})^{-1}}=\frac{1}{a}$ and for every $x_1,x_2\in{\bf No}$ such that $x_1,x_2>+\infty_\omega$ we obtain

\begin{equation}
\label{f232}
\begin{array}{l}
a^{x_1+x_2}=\frac{1}{(\frac{1}{a})^{x_1+x_2}}=\frac{1}{(\frac{1}{a})^{x_1}\cdot(\frac{1}{a})^{x_2}}=a^{x_1}\cdot a^{x_2}.
\end{array}
\end{equation}

Notice also that for any two numbers $x_1,x_2\in{\bf No}$ such that $a^{x_1}=a^{x_2}$, then $\frac{a^{x_1}}{a^{x_2}}=1$ and also $a^{x_1-x_2}=1$ thus $x_1-x_2=0$ and $x_1=x_2$ whenever $a>1$ or $a<1$. Thus the exponent function $y=a^x$ is injective.

{\bf Remark 2}. In the definition of $a^{x'}=\omega^{\frac{x'}{\omega}}$ we did not use here usual formula in real numbers $\mathbb R$ as ${(a^b)}^c=a^{b\cdot c}$, which does not work for ininfinite number, i.e., $a^{x'}=a^{\omega\cdot\frac{x'}{\omega}}=(a^{\omega})^{{\frac{x'}{\omega}}}=\omega^{\frac{x'}{\omega}}$, because, e.g., $a^{2\cdot\omega}\not=(a^2)^\omega$, otherwise, for $x=2\cdot\omega=x'$ we should obtain $a^{2\cdot\omega}=(a^2)^\omega=\omega\not=\omega^\frac{2\cdot\omega}{\omega}=\omega^2$.

Let $a\in{\bf No}^+$ such that $a\not=1$. Now we are going to define $y=\log_a x$ for some (not all) numbers $x>+\infty_\omega$. 

{\bf Definition 18.} By $X^+\subset{\bf No}^+$ we denote the set of all $x>+\infty_\omega$ such that in a canonical form $\sum\limits_{0\leq\beta<\alpha}\omega^{y_\beta}\cdot r_\beta$ of $x$ the number $y_0=\sum\limits_{0\leq\gamma<\alpha_0}\omega^{z_\gamma}\cdot s_\gamma$ has the following property: $z_\gamma>-1$,  for all $0\leq\gamma<\alpha_0$, where $\sum\limits_{0\leq\gamma<\alpha_0}\omega^{z_\gamma}\cdot s_\gamma$ is a canonical form of $y_0$. By $X^-\subset{\bf No}^+$ we denote the set of all $x\in{\bf No}^+$ such that $\frac{1}{x}\in X^+$. And at last, by $X^0$ we denote the set $(\frac{1}{+\infty_\omega},+\infty_\omega)$.

{\bf Lemma 12.} $x\in X^+$ if and only if there is a number $y\in{\bf No}$ such that $x=a^y$. Also $x\in X^-$ if and only if there is a number $y\in{\bf No}$ such that $x=\frac{1}{a^y}$.

{\bf Proof.}
If $x\in X^+$, then we consider a canonical form $x=\sum\limits_{0\leq\beta<\alpha}\omega^{y_\beta}\cdot r_\beta$, where  numbers $y_\beta$ form a descending sequence, i.e., $y_\beta>y_{\beta'}$ for all $0\leq\beta<\beta'<\alpha$, and $r_\beta$ are real numbers. Since $x$ is a positive number hence $r_0>0$  and thus we can write
\begin{equation}
\label{f4748}
x=\omega^{y_0}\cdot r_0\cdot(1+\delta)
\end{equation}
where 
\begin{equation}
\label{f4749}
\delta=\sum\limits_{0<\beta<\alpha}\omega^{y_\beta-y_0}\cdot \frac{r_\beta}{r_0}
\end{equation}
and $\delta$ is an infinitesimal number.

 Now we can define $y$ as a value of the logarithmic function in a point $x$ by the following formulas:
\begin{equation}
\label{f133}
y=\log_ax\stackrel{def}{=}\omega\cdot{y_0}+\log_ar_0+\log_a(1+\delta),
\end{equation}
 which is well-defined because $\log_ar_0$ real-valued logarithm and $\log_a(1+\delta)$ is defined above because $1+\delta\in X^0$. It is really an exponent whose power $a^y=x$. Indeed, 
\begin{equation}
\label{f4750}
a^y=a^{\omega\cdot{y_0}+\log_ar_0+\log_a(1+\delta)}=a^{\omega\cdot y_0}\cdot a^{\log_ar_0}\cdot a^{\log_a(1+\delta)}=
\omega^{y_0}\cdot r_0\cdot(1+\delta)=x,
\end{equation}
because $a^{\omega\cdot y_0}=\omega^{y_0}$ since $\omega\cdot y_0=\sum\limits_{0\leq\gamma<\alpha_0}\omega^{z_\gamma+1}\cdot s_\gamma$, by our definition of $X^+$, and $z_\gamma+1>0$ for all $0\leq\gamma<\alpha_0$, and thus the formula $a^{\omega\cdot y_0}=\omega^{y_0}$ is correct.

Consider now an arbitrary number $y\in{\bf No}$ such that $y>+\infty_\omega$. We have to prove that $a^y=x\in X^+$. We know that $y=y'+y''$, where $y'$ is an infinite part of $y$ and $y''$ is a finite part of $y$. Then $a^y=a^{y'+y''}=\omega^{\frac{y'}{\omega}}\cdot a^{y''}=x\in X^+$. Indeed, consider a canonical form $\sum\limits_{0\leq\beta<\alpha}\omega^{y_\beta}\cdot r_\beta$ of $x$. It is clear that $y_0=\frac{y'}{\omega}$ because the canonical form of $\omega^{\frac{y'}{\omega}}$ is $\omega^{\frac{y'}{\omega}}$ itself and the canonical form of $\frac{1}{+\infty}<a^{y''}<+\infty_\omega$ is the following: $\omega^0\cdot r_0+\sum\limits_{0<\beta<\alpha}\omega^{z_\beta}\cdot s_\beta$. Thus the caninical form of $x$ is the following: $\omega^{\frac{y'}{\omega}}\cdot(\omega^0\cdot r_0+\sum\limits_{0<\beta<\alpha}\omega^{z_\beta}\cdot s_\gamma)=\omega^{\frac{y'}{\omega}}\cdot r_0+\sum\limits_{0<\gamma<\alpha}\omega^{z_\gamma+\frac{y'}{\omega}}\cdot s_\gamma$. And since $x=\sum\limits_{0\leq\beta<\alpha}\omega^{y_\beta}\cdot r_\beta=\omega^{\frac{y'}{\omega}}\cdot r_0+\sum\limits_{0<\gamma<\alpha}\omega^{z_\gamma+\frac{y'}{\omega}}\cdot s_\gamma$ we obtain the following equality: $y_0=\frac{y'}{\omega}$ and thus the canonical form of $\sum\limits_{0\leq\gamma<\alpha_1}\omega^{z_\gamma}\cdot s_\gamma$ of $y_0$  has the following property: $z_\gamma>-1$ for all $0\leq\gamma<\alpha_1$. Indeed, let $\sum\limits_{0\leq\gamma<\alpha_1}\omega^{t_\gamma}\cdot s_\gamma$ be a canonical form of $y'$. Since $y'$ is an ininite  initial integer then, by Definition 18, $t_\gamma>0$ for all $0\leq\gamma<\alpha_1$ and thus $z_\gamma>-1$ for all $0\leq\gamma<\alpha_1$ because $z_\gamma=t_\gamma-1$ for all $0\leq\gamma<\alpha_1$. 

If $x\in X^-$, then $\frac{1}{x}\in X^+$ and, by the first part of proving, it is iff there is a number $y$ such that $\frac{1}{x}=a^y$ and hence $x=\frac{1}{a^y}$.
$\Box$

And it is the end of proving Lemma 12.

Now we continue the proof of Theorem 8. For each $x\in X^-$ we put $y=\log_a x=-\log_a\frac{1}{x}$. Really it is needed logarithm because $a^{-\log_a\frac{1}{x}}=\frac{1}{a^{\log_a\frac{1}{x}}}=\frac{1}{\frac{1}{x}}=x$.

Thus, we have defined a logarithmic function $y=\log_ax$ for all $x\in X=X^-\cup X^0\cup X^+$ by three formulas on each component of $X$.

At last, if $a<1$, then we define $y=\log_ax\stackrel{def}{=}-\log_{\frac{1}{a}}x$ for all $x\in X=X^-\cup X^0\cup X^+$.

Notice that for every $x_1,x_2\in X^-\cup X^0\cup X^+$ there is a natural formula
\begin{equation}
\label{f10}
\log_a(x_1\cdot x_2)=\log_ax_1+\log_ax_2.
\end{equation}
Indeed,  consider canonical forms of numbers  $x_1$ and $x_2$, respectively, in informed kind $(\ref{f4748})$  $x_1=\omega^{y^{(1)}_0}\cdot r^{(1)}_0\cdot(1+\delta^{(1)})$ and $x_2=\omega^{y^{(2)}_0}\cdot r^{(2)}_0\cdot(1+\delta^{(2)})$, respectively. Then by formula $(\ref{f133})$, we obtain the following
formulas:
 \begin{equation}
 \label{f321}
 \begin{array}{l}
 \log_ax_1=\omega\cdot{y^{(1)}_0}+\log_ar^{(1)}_0+\log_a(1+\delta^{(1)}_0)
\end{array}
\end{equation} 
 and
 \begin{equation}
 \label{f421}
 \begin{array}{l}
\log_ax_2=\omega\cdot{y^{(2)}_0}+\log_ar^{(2)}_0+\log_a(1+\delta^{(2)}),
\end{array}
\end{equation} respectively.
 
 Consider now 
 \begin{equation}
 \label{f121}
 \begin{array}{l}
 \log_ax_1+\log_ax_2=\\\omega\cdot{y^{(1)}_0}+\log_ar^{(1)}_0+\log_a(1+\delta^{(1)})_+ \omega\cdot{y^{(2)}_0}+\log_ar^{(2)}_0+\log_a(1+\delta^{(2)})=\\\omega\cdot({y^{(1)}_0}+{y^{(2)}_0})+\log_ar^{(1)}_0+\log_ar^{(2)}_0+\log_a(1+\delta^{(1)})+\log_a(1+\delta^{(2)})=\\
 \omega\cdot({y^{(1)}_0}+{y^{(2)}_0})+\log_a(r^{(1)}_0\cdot r^{(2)}_0\cdot(1+\delta^{(1)})\cdot(1+\delta^{(2)}))
 \end{array}
 \end{equation}
One can see that
\begin{equation}
 \label{f221}
 \begin{array}{l}
a^{\log_ax_1+\log_ax_2}=a^{\omega\cdot({y^{(1)}_0}+{y^{(2)}_0})+\log_a(r^{(1)}_0\cdot r^{(2)}_0\cdot(1+\delta^{(1)})\cdot(1+\delta^{(2)}))}=\\
\omega^{y^{(1)}_0}\cdot\omega^{y^{(2)}_0}\cdot r^{(1)}_0\cdot r^{(2)}_0\cdot(1+\delta^{(1)})\cdot(1+\delta^{(2)})=\\
\omega^{y^{(1)}_0}\cdot r^{(1)}_0\cdot(1+\delta^{(1)})\cdot\omega^{y^{(2)}_0}\cdot r^{(2)}_0\cdot(1+\delta^{(2)})=\\
x_1\cdot x_2=a^{\log_a(x_1\cdot x_2)}.
\end{array}
 \end{equation}
Thus the following formula
\begin{equation}
\label{f333}
\log_ax_1+\log_ax_2=\log_ax_1\cdot x_2.
\end{equation}
$\Box$

{\bf Corollary 7.} {\it For each fixed number $a\in{\bf No}^+$ such that  $a\in X^{-1}\cup X^+$  there is an exponential function $y=a^x$ whose domain is ${\bf No}$  and whose rang of values is $X^-\cup X^+$ there is  a logarithmic function $y=\log_ax$}.

{\bf Proof}. If $a\in X^{-1}\cup X^+$, then, by Theorem 20 in $\cite{l3}$, there is a well-defined value $\log_ba$ for any number $b\in(\frac{1}{+\infty_\omega},+\infty_\omega)$ such that $b\not=1$. Then for each $x\in {\bf No}$ we put $y=a^x\stackrel{def}{=}b^{x\cdot\log_ba}$. This number $b^{x\cdot\log_ba}$ uniquely defined by Theorem 20, and the definition of $a^x=b^{x\cdot\log_ba}$ is correct. It is clear that such definition of $a^x$ does not depend of a choice of $b\in (\frac{1}{+\infty_\omega},+\infty_\omega)$, $b\not=1$. Indeed, for any $c\in (\frac{1}{+\infty_\omega},+\infty_\omega)$, $c\not= 1$ we obtain $a^x=c^{x\cdot \log_ca}=c^{\log_cb\cdot x\cdot\frac{\log_ca}{\log_cb}}=(c^{\log_cb})^{x\cdot\frac{\log_ca}{\log_cb}}=b^{x\cdot\log_ba}$.
 Moreover, for any $x_1, x_2\in {\bf No}$ we have $a^{x_1+x_2}=b^{(x_1+x_2)\cdot\log_ba}=b^{x_1\cdot\log_ba}\cdot b^{x_2\cdot\log_ba}=a^{x_1}\cdot a^{x_2}$.

Then for each  $x\in  X^+$ we  put $\log_ax\stackrel{def}{=}\frac{\log_bx}{\log_ba}$. This definition is also correct since numbers $\log_bx$ and $\log_ba$ are unuquely defined by Theorem 20, and $\log_ba\not=0$ because $a>+\infty_\omega$ or $0<a<\frac{1}{+\infty_\omega}$.

One can see, by formula $(\ref{f131})$, that $a^{\frac{\log_bx}{\log_ba}}=b^{\log_ba\cdot\frac{\log_bx}{\log_ba}}=b^{\log_bx}=x$. Moreover, if $x_1,x_2\in X^+$, then $\log_a(x_1\cdot x_2)=\frac{\log_b(x_1\cdot x_2)}{\log_ba}=\frac{\log_bx_1+\log_bx_2}{\log_ba}=\frac{\log_bx_1}{\log_ba}+
\frac{ x_2}{\log_ba}=\log_ax_1+\log_ax_2$. Clearly, a definition and properties of $\log_ax$ do not depend on a choice of $b$. Indeed, for any $c\in (\frac{1}{+\infty_\omega},+\infty_\omega)$, $c\not= 1$ we obtain $\log_ax=\frac{\log_cx}{\log_ca}=\frac{\log_cx}{\log_cb}\cdot\frac{\log_cb}{\log_ca}=\frac{\log_bx}{\log_ba}$.

If $x\in X^-$, then $\frac{1}{x}\in X^+$ and  we put $\log_ax=-\log_a\frac{1}{x}$. $\Box$

{\bf Remark 3}. Notice that for $a=\omega$ and $x\in {\bf No}$ such that $x\cdot \omega\in {\bf Z}'$ the value $\omega^x$ in Conway's definition coincides with ours here. Indeed,   $\omega^x=b^{x\cdot\log_b\omega}=b^{x\cdot\omega}$ and $b^{x\cdot\omega}\stackrel{def}{=}b^{\omega\cdot\frac{x\cdot\omega}{\omega}}=\omega^{ x}$.
$\Box$

{\bf Corollary 8.} {\it Let ${\mathbb R}_\zeta$, $\zeta=\omega^{\omega^\mu}$, $0<\mu<\Omega$, be a topological space with ${\mathbb R}_\zeta$-topology, i.e., a lineary-ordered topology, and  $a\in{\mathbb R}_\zeta^+\cap X$ be a fixed number such that  $a\not= 1$, Where $X=X^-\cup X^0\cup X^+$. Then there are continuous  functions  $ y=a^x$ whose domain is $\bar X_1\subseteq{\mathbb R}_\zeta$ and whose range of values is $\bar X_2= a^{(X_1)}$,  and $y=\log_ax$ whose domain is $\bar X_2$ and whose range of values is $\bar X_1$, which are inverse with each other, and when $a\in{\mathbb R}^+$, then the restrictions of these functions on ${\mathbb R}$ and ${\mathbb R}^+$ are the real-valued exponential function and logarithmic function, respectively. }

{\bf Proof.}  
Let $a\in{\mathbb R}_\zeta\cap X$ be a number such that $a\not=1$, $\zeta=\omega^{\omega^\mu}$, $0<\mu<\Omega$. By Theorem 20 and Corollary 7,
there are functions $y=a^x$  and $y=\log_ax$, defined for all $x\in{\bf No}$ and for all $x\in X=X^-\cup X^0\cup X^+$ (when $\frac{1}{+\infty_\omega}<a<+\infty_\omega$) and for all $x\in X^-\cup\{1\}\cup X^+$ (when $a>+\infty_\omega$ or $a<\frac{1}{+\infty_\omega}$),  respectively, which are inverse with each other and when $a\in{\mathbb R}$ and each $x\in{\mathbb R}$ and $x\in{\mathbb R}^+$. 

Consider now the restrictions of them on ${\mathbb R}_\zeta$ and ${\mathbb R}_\zeta^+\cap X$, when $a\in X^0$, and on ${\mathbb R}_\zeta\cap(X^-\cup\{1\}\cup X^+)$, when $a>+\infty_\omega$ or $a<\frac{1}{+\infty_\omega}$, respectively. We define $\bar X_1=\{x\in {\mathbb R}_\zeta\,|\, a^x\in {\mathbb R}_\zeta^+\cap X\}$ as a domain of the funcion $y=a^x$ in ${\mathbb R}_\zeta$ and its range of values is $\bar X_2=a^{(\bar X_1)}$. Then we define a function $y=\log_ax$ in ${\mathbb R}_\zeta$ as a restriction of $y=\log_ax$ on $\bar X_2$, where, evidently, $\bar X_2={\mathbb R}_\zeta^+\cap X$ is a domain of the function $y=\log_ax$ in ${\mathbb R}_\zeta$ and its range of values is $\bar X_1$. Clearly,
the restrictions of these functions on ${\mathbb R}$ and ${\mathbb R}^+$ are the real-valued expotential function and logarithmic function, respectively.

We have to prove now that $y=a^x$ and $y=\log_ax$ are $\zeta$-continuous functions.

Notice that functions $y=a^x$ and $y=\log_ax$ are monotone, more exactly, for $a>1$ their are monotonic incresing and for $a< 1$ their are monotonic decreasing. Indeed, if  $x_1,x_2\in \bar X_1$ such that $x_1<x_2$ and $a>1$, then $a^{x_2}-a^{x_1}=a^{x_2}\cdot(1-\frac{a^{x_1}}{a^{x_2}})=a^{x_2}\cdot(1-a^{x_1-x_2})>0$ and hence $a^{x_1}<a^{x_2}$. If
 $x_1,x_2\in \bar X_1$ such that $x_1<x_2$ and $a<1$, then $a^{x_2}-a^{x_1}=a^{x_2}\cdot(1-\frac{a^{x_1}}{a^{x_2}})=a^{x_2}\cdot(1-a^{x_1-x_2})>0$ and hence $a^{x_1}>a^{x_2}$. 
 
 Since $y=\log_ax$ is an inverse function to $y=a^x$, then it is also monotonic increasing for $a>1$ and monotonic decreasing for $a<1$.
 
Now if we suppose that $y=a^x$  is not continuous (say) in a point $x_0\in \bar X_1$, then there is a number $\varepsilon\in(0,1)\subset {\mathbb R}_\zeta$ such that for each $\frac{1}{\alpha}\in{\mathbb R}_\zeta$, $0<\alpha<\zeta$, there is a number $x_\alpha\in \bar X_1$ with $|x_\alpha-x_0|<\frac{1}{\alpha}$ and $|a^{x_\alpha}- a^{x_0}|>\varepsilon$. We shall show that there is a number $\bar y\in \bar X_2$ such that there is no $\bar x\in \bar X_1$ with $a^{\bar x}=\bar y$. 

Really,  if $x_\alpha>x_0$, for all $0<\alpha<\zeta$, then for all $x>x_0$ we obtain $|a^x-a^{x_0}|=a^x-a^{x_0}>\varepsilon$, by monotonic increasing function $y=a^x$, for $a>1$ and some $0<\alpha<\zeta$ such that $x_0<x_\alpha<x$, because $a^{x_\alpha}<a^x$; and $|a^x-a^{x_0}|=a^{x_0}-a^x>\varepsilon$, by monotonic decreasing function $y=a^x$, for $a<1$ and some $0<\alpha<\zeta$ such that $x_0<x_\alpha<x$, because $a^{x_\alpha}>a^x$. Thus, a desired number is  $\bar y=a^{x_0}+\frac{\varepsilon}{2}$ (a possible, in general, case for a discontinuous monotone function when it is continuous from the left and break from the right).

If $x_\alpha>x_0$, for all $0<\alpha<\zeta$, then for all $x<x_0$ we obtain $|a^x-a^{x_0}|=a^{x_0}-a^x>\varepsilon$, by monotonic increasing function $y=a^x$, for $a>1$ and some $0<\alpha<\zeta$ such that $x<x_\alpha<x_0$, because $a^xa^{x_\alpha}<a^{x_\alpha}$; and $|a^x-a^{x_0}|=a^x-a^{x_0}>\varepsilon$, by monotonic decreasing function $y=a^x$, for $a<1$ and some $0<\alpha<\zeta$ such that $x<x_\alpha<x_0<x$, because $a^{x_\alpha}<a^x$. Thus, a desired number is  $\bar y=a^{x_0}-\frac{\varepsilon}{2}$ (a possible, in general, case for a discontinuous monotone  function when it is continuous from the right and break from the left).

If $x_{\alpha'}>x_0$ and $x_{\alpha''}<x_0$ for all $\alpha',\alpha''\in(0,\zeta)$ such that $\{\alpha'\}$ and $\{\alpha''\}$ are co-final in $(0,\zeta)$, then  for all $x>x_0$ we obtain $|a^x-a^{x_0}|=a^x-a^{x_0}>\varepsilon$, by monotonic increasing function $y=a^x$, for $a>1$ and some $0<\alpha'<\zeta$ such that $x_0<x_{\alpha'}<x$, because $a^{x_{\alpha'}}<a^x$; and $|a^x-a^{x_0}|=a^{x_0}-a^x>\varepsilon$, by monotonic decreasing function $y=a^x$, for $a<1$ and some $0<\alpha''<\zeta$ such that $x_0<x_{\alpha''}<x$, because $a^{x_{\alpha''}}>a^x$. Thus, a desired number is  $\bar y=a^{x_o}\pm\frac{\varepsilon}{2}$, (a possible, in general, case for a discontinuous monotone function when it is break from the right and break from the left). 

In all cases the existence of $\bar y$ contradicts to the fact that $\bar y\in \bar X_2$ since $a^{x_0}\pm\varepsilon=a^{x_0}\cdot(1\pm\frac{\varepsilon}{a^{x_0}})=a^{x_0}\cdot a^{\log_a(1\pm\frac{\varepsilon}{a^{x_0}})}=a^{x_0+\log_a(1\pm\frac{\varepsilon}{a^{x_0}})}$ and, by Lemma 12, $\bar x=x_0+\log_a(1\pm\frac{\varepsilon}{a^{x_0}})\in \bar X_1$. Moreover, $\bar X_1$ and $\bar X_2$ are bijective and this bijection if given by $y=a^x$. Hence there is a number $\bar x$ such that $a^{\bar x}=\bar y$ . Contradiction.

Thus, $y=a^x$ is continuous and by strictly monotonity of it. 

Absolutely analogously, $y=\log_ax$ is  a strictly monotone continuous function in ${\mathbb R}_\zeta$, where $\frac{1}{+\infty_\zeta}<a<+\infty_\zeta$ and we omit this exposure.

\begin{center}
{\bf 10. Trigonometry in ${\bf No}$ and $\zeta$-trigonometry in ${\mathbb R}_\zeta$, $\zeta=\omega^{\omega^\mu}$, $0<\mu<\Omega$}
\end{center}

Conway defined trigonometry functions $y=\sin x$ and $\cos x$ for all numbers $x\in(-\infty_\omega,+\infty_\omega)$ by the following steps.

First, he used the following formulas for all infinitesimal $-\frac{1}{+\infty_\omega}\delta<\frac{1}{+\infty_\omega}$
\begin{equation}
\label{f6223}
\sin x=x-\frac{x^3}{3!}+...+(-1)^{n-1}\frac{x^{2n-1}}{(2n-1)!}+...
\end{equation}
and
\begin{equation}
\label{f6213}
\cos x=1-\frac{x^2}{2!}+...+(-1)^n\frac{x^{2n}}{2n!}+...
\end{equation}

Then for each  infinitesimal $\delta$ formulars $\ref{f6213}$ and $\ref{f6213}$ one can define uniquely numbers $\sin\delta$ and $\cos\delta)$, respectively. Really, if $\delta>0$, then
\begin{equation}
\label{f5203}
\begin{array}{l}
\sin\delta=\delta-\frac{\delta^3}{3!}+...+\frac{\delta^{2n-1}}{(2n-1)!}+...=\\=\{\delta,\delta-\frac{\delta^3}{3!},...,\delta-\frac{\delta^3}{3!}+...+(-1)^{n-1}\frac{\delta^{2n-1}}{(2n-1)!},...\,|\,1,\frac{1}{2},...,\frac{1}{n},...\}.
\end{array}
\end{equation}
If $\delta<0$, then
\begin{equation}
\label{f6203}
\begin{array}{l}
\sin\delta=\delta-\frac{\delta^3}{3!}+...+\frac{\delta^{2n-1}}{(2n-1)!}+...=\\=\{1,\frac{1}{2},...,\frac{1}{n},...\,|\,\delta,\delta-\frac{\delta^3}{3!},...,\delta-\frac{\delta^3}{3!}+...+(-1)^{n-1}\frac{\delta^{2n-1}}{(2n-1)!},...\}.
\end{array}
\end{equation}

Similarly, if $\delta$ is an arbitrary such that $\delta\not=0$, then
\begin{equation}
\label{f6204}
\begin{array}{l}
\cos\delta-1=-\frac{\delta^2}{2!}+\frac{\delta^4}{4!}-...+(-1)^n\frac{\delta^{2n}}{2n!}+...=\\=\{-1,-\frac{1}{2},...,-\frac{1}{n},...\,|\,-\frac{\delta^2}{2!},-\frac{\delta^2}{2!}+\frac{\delta^4}{4!},...,-\frac{\delta^2}{2!}+\frac{\delta^4}{4!}+...+(-1)^n\frac{\delta^{2n}}{2n!},...\}.
\end{array}
\end{equation}

A second step is the following. Since every number $x\in(-\infty_\omega,+\infty_\omega)$ has a form $x=x'+\delta$, where $x'\in{\mathbb R}$, then $\sin x\stackrel{def}{=}\sin (x'+\delta)=\sin x'\cdot \cos\delta+\cos x'\cdot\sin\delta$ and $\cos x\stackrel{def}{=}\cos (x'+\delta)=\cos x'\cdot \cos\delta+\sin x'\cdot\sin\delta$.

A third step is an extension of $\sin x$, $\cos x$  and $\tan x$ from$(-\infty_\omega,+\infty_\omega)$ to $(-\infty_\zeta,+\infty_\zeta)$. For each $x\in{\mathbb R}_\zeta$ there is a decoposition $x=x'+x''$, where $x'$ is an initial ineger and $x''$ is a finite number. Then we put $\sin x=\sin x''$, $\cos x=\cos x''$ and $\tan x=\tan x''$, respectively, if in latter case $\cos x''\not=0$.

We denine these trigonometry functions on ${\mathbb R}_\zeta$, $\zeta=\omega^{\omega^\mu}$, $0<\mu<\Omega$, as $\sin x=\{\sin x\}\cap{\mathbb R}_\zeta$, $\sin x=\{\sin x\}\cap{\mathbb R}_\zeta$ and $\sin x=\{\sin x\}\cap{\mathbb R}_\zeta$, where $\{\sin x\}$, $\{\cos x\}$ and $\{\tan x\}$ are already defined in ${\bf No}$. It's not so obvious that $\sin({\mathbb R}_\zeta)\subset{\mathbb R}_\zeta$, $\sin({\mathbb R}_\zeta)\subset{\mathbb R}_\zeta$ and $\sin({\mathbb R}_\zeta)\subset{\mathbb R}_\zeta$, respectively. But it is so if one defines them in ${\mathbb R}_\zeta$ geometrically as functions of  angles $-360^\circ\leq\theta\leq 360^\circ$ on the  $\zeta$-circle $S_\zeta(0,0)$ with radius $r>0$ and the center in the origin (0,0) as the set of all points $A(x,y)$ in ${\mathbb R}_\zeta\times{\mathbb R}_\zeta$ satisfying the following equation: $x^2+y^2-r^2=0$ and replace the angles with {\it radians}=$x$, $0\leq x<2\pi$ for $\sin x$ and $\cos x$.

 Consider now ${\mathbb C}_\zeta=\{z\,\,|\,\,z=a+b  i,\,a,b\in{\mathbb R}_\zeta,\,\, i^2=-1\}$ in Theorem 16 which is a closed field with the following operations: 
 \begin{equation}
 \label{8090}
z_1+z_2=(a_1+a_2)+(b_1+b_2) i, 
 \end{equation}
 \begin{equation}
 \label{8091}
 z_1\cdot z_2=(a_1a_2-b_1b_2)+(a_1b_2+a_2b_1) i, 
 \end{equation}
 \begin{equation}
 \label{8092}
 \frac{z_1}{z_2}=\frac{a_1+b_1 i}{a_2+b_2 i}=\frac{a_1a_2+b_1b_2}{a^2+b^2}+\frac{a_1b_2-a_2b_1}{a^2+b^2} i,
 \end{equation}
  where $z_2\not=0$.
 
 The trigonometry form of $z=a+b i$ is clear: $z=r(\cos \theta +\sin\theta i)$, where $r=|z|=\sqrt{a^2+b^2}>0$, $\theta\in[0,2\pi)$, $\cos\theta=\frac{a}{r}$ and $\sin\theta=\frac{b}{r}$. Then $z_1\cdot z_2=r_1r_2(\cos(\theta_1+\theta_2)+\sin(\theta_1+\theta_2) i)$.
 
 Let $S_\zeta$ be the unit circle, i.e. $S_\zeta\stackrel{def}{=}\{z|\,\,x\in{\mathbb C}_\zeta,\,\,|z|=1\}$. Then $S_\zeta$ is a multiplicative group $(S_\zeta,\cdot)$. 
 
 {\bf Proposition 6.} There is a continuous homomorphism $ex:{\mathbb R}_\zeta\rightarrow S_\zeta$ of the additive group $({\mathbb R}_\zeta,+)$ to the multiplicative group $(S_\zeta,\cdot)$, given by the following formula:
 \begin{equation}
 \label{f1507}
 ex(x)=\cos{2\pi x}+i\sin{2\pi x},\,\,\, x\in{\mathbb R}_\zeta.
 \end{equation}
 
 {\bf Proof.} Let $x,y$ be arbitrary numbers in ${\mathbb R}_\zeta$. By formula $(\ref{f1507})$ and formulas $\cos(x+y)$ and $\sin(x+y)$ , we obtain $ex(x+y)=\cos 2\pi(x+y)+i\sin2\pi(x+y)=(\cos 2\pi x\cos 2\pi y-\sin 2\pi x\sin 2\pi y)+i(\sin 2\pi x\cos 2\pi y+\cos 2\pi x\sin 2\pi y)=(\cos 2\pi x+i\sin 2\pi x)\cdot(\cos 2\pi y+i\sin 2\pi y)=ex(x)\cdot ex(y)$. Thus, $ex$ is a homomorphism.
 
 Since ${\mathbb R}_\zeta$ and $S_\zeta$ are groups it is enough to verify the continuous of $ex$ at the point $x=0$.  We consider the topology on $S_\zeta$ as the   topology induced by the topology of Cartesian product topology of ${\mathbb R}_\zeta$ on ${\mathbb C}_\zeta={\mathbb R}_\zeta\times{\mathbb R}_\zeta$ and the topology on ${\mathbb R}_\zeta$ is a topology defined by the linear oder on it. Let $U$ be any  neighborhood of the number $1\in S_\zeta$. Then there are neighborhoods, say $(a,b)$, $a<x<b$, of $1$ and $(c,d)$ of $0$ in ${\mathbb R}_\zeta$ such that $(a,b)\times(c,d)\cap S\subset U$. Let $A$ and $B$ are the intersection of $S_\zeta$ and the boundary $Fr([a,b]\times[c,d])$ of $[a,b]\times[c,d]$ in ${\mathbb C}_\zeta$, say $A$ with negative ordinate and $B$ with positive ordinate. And let $a'$ be the length of an angle $\angle AOC$ and $b'$ the length of an angle $\angle BOC$ in radians, where $C(1,0)$. Consider the following open neighborhood $(-a',b')$ of $0$ in ${\mathbb R}_\zeta$. Then for each $x\in(-a',b')$ we have $ex(x)\in U$. Then for any neighborhood $U$ of $1\in S_\zeta$ we see that for $V=(a,b)$ the inclusion $ex(V)\subset U$ fulfills. Thus, $\varphi$ is a continuous map at the point $0$.
 
 Another proof is to show that $ex:{\mathbb R}_\zeta\rightarrow S_\zeta$ is a local homeomorphism. Indeed, for each point $x\in{\mathbb R}_\zeta$ the restriction $ex|_{(x-\frac{1}{2},x+\frac{1}{2})}$ is a homeomorphism $(x-\frac{1}{2},x+\frac{1}{2})$ onto $S_\zeta\setminus\{ex(x)\}$ because $ex(x)=ex(y)$ iff $x-y-tr(x-y)$ is an integer in ${\mathbb Z}\subset{\mathbb R}_\zeta$, where $tr(x-y)$ is a transfinite part of $x-y$. But each local homeomorphism is a continuous mapping. $\Box$
 
 {\bf Theorem 11.} {\it Every polynomial $p(z)=z^n+Az^{n-1}+Bz^{n-2}+...+K$, where $A,B,...,K\in {\mathbb C}_\zeta$,  has a root in ${\mathbb C}_\zeta$}.
 
 {\bf Proof}. Consider the case when $K\not=0$, otherwise, the root $z_0=0$. For every $p$ which has no root on $S_\zeta\subset{\mathbb C}_\zeta$ we define a map $\hat p:S_\zeta\rightarrow {\mathbb U}_\zeta$, where ${\mathbb U}_\zeta\stackrel{def}{=}\{z\,|\, z\in{\mathbb C}_\zeta\,\,\&\,\,|z|=1\}$, by the following formula:
 \begin{equation}
 \label{f333444}
 \hat p(z)=\frac{p(z)}{|p(z)|},
 \end{equation}
 and prove Theorem 11 in two steps:
 
 $1).$ If $p$ has no root $z$ with $|z|\leq 1$ then deg$(\hat p)=0$.
 
 $2).$ If $p$ has no root $z$ with $|z|\geq 1$ then deg$(\hat p)=n$.
 
 For case $1).$ we consider the deformation $\hat p_t:{\mathbb U}_\zeta\rightarrow {\mathbb U}_\zeta$ given by
\begin{equation}
 \label{f333445}
 \hat p_t(z)=\frac{p(tz)}{|p(tz)|}.
 \end{equation}
 Clearly $\hat p_1=\hat p$ and $\hat p_0=$constant which is equal $1$ or $-1$, hence deg$(\hat p)=0$. 
 
 For case $2).$ we consider the deformation $\hat p_t:{\mathbb U}_\zeta\rightarrow {\mathbb U}_\zeta$ given by
\begin{equation}
 \label{f333446}
 \hat p_t(z)=\frac{q(z,t)}{|q(z,t))|},
 \end{equation}
where
\begin{equation}
 \label{f333447}
 q(z,t)=t^np(\frac{z}{t})=z^n+t(Az^{n-1}+tBz^{n-2}+...+t^{n-1}K).
 \end{equation}
 The right side of $(\ref{f333447})$ shows that $q(z,t)$ is continuous (even when $t=0$).
 Cleraly $\hat p_1=\hat p$ and $\hat p_0=z^n$, hence, deg$(\hat p)=$deg$(\hat p_0)=n$, what is in contradiction with case $1).$ where deg$(\hat p)=0$. This contradiction can be removed only when there is $z_0\in{\mathbb C}_\zeta$
such that $p(z_0)=0$, i.e., $p(z)$ has a root,
 and the above deformation in case $1).$ is impossible.

By deg$(\hat p)$, i.e., degree of a map, we understand the following. Consider a continuous map $s:[0,1]\rightarrow{\mathbb U}_\zeta$ such that given by formula $s(t)=\cos(2\pi t)+i\sin(2\pi t)$. Then for composition $\hat p\circ s$ there exists a continuous mapping $\tau:[0,1]\rightarrow{\mathbb R}''_\zeta$ such that $ex\circ \tau=\hat p\circ s$. Indeed, $\tau$ is given by the following formula:
$\tau(0)=exp(x_0),\,\,\,\tau(t)=\frac{1}{2\pi i}\ln(\hat p( s(t))),\,\,t\in[0,1],$
where $x_0$ is any element in ${\mathbb R}''_\zeta$ such that $ex(x_0)=\hat p( s(0))$.
Now we define a degree deg$(\hat p)$ as an integer $\tau(1)-\tau(0)$. This deg$(\hat p)$ is homotopy invariant. We omit details.

\begin{center}
{\bf 11. Appendix and last remarks: $\neg\exists X(\forall x) (x\in X)$ vs $\exists X\neg\{X\}$}
 \end{center}
 
Finally we consider a linear ordering topology on ${\bf No}$ and elucidate what transfinite $\Omega$-fundamental sequences of numbers in ${\bf No}$ can be called {\it convergent} and what transfinite $\Omega$-fundamental  sequences of numbers in ${\bf No}$ can be called {\it non-convergent} in this ordering topology on ${\bf No}$.

{\bf Definition 19}. A mapping $x:(0,\lambda)\rightarrow {\bf No}$ is called a {\it transfinite sequence} of a type $\lambda$ of Conway's numbers in ${\bf No}$, or $\lambda$-{\it transfinite sequence} in ${\bf No}$, or shortly a $\lambda$-{\it sequence} in ${\bf No}$, when  $\lambda$ is a {\it limit ordinal number} such that $\omega\leq\lambda\leq\Omega$. We denote it as above by $(x_\alpha)_{0<\alpha<\lambda}$.

{\bf Definition 20}.
We say that $\Omega$-sequence $(x_\alpha)_{0<\alpha<\Omega}$ in  ${\bf No}$, {\it converges} to $a\in {\bf No}$, and we write $\lim\limits_{0<\alpha<\Omega}x_\alpha=a$,  if for each positive  number $\varepsilon\in {\bf No}$ there is an ordinal number $\alpha_0\in(0,\Omega)$ such that $|x_\alpha-a|<\varepsilon$ for all $\alpha_0\leq\alpha<\Omega$.  
In this case we also say that $\Omega$-sequence $(x_\alpha)_{0<\alpha<\Omega}$ is {\it convergent} in  ${\bf No}$. If $\lim\limits_{0<\alpha<\Omega}x_\alpha=0$, then $(x_\alpha)_{0<\alpha<\Omega}$ is called an $\Omega$-{\it infinitesimal infinite sequence} in ${\bf No}$, or shortly  an $\Omega$-{\it infinitesimal} in ${\bf No}$.

{\bf Definition 21.}
 An $\Omega$-sequence $(x_\alpha)_{0<\alpha<\Omega}$ in ${\bf No}$,  is called  {\it fundamental}  or a {\it Cauchy} $\Omega$-sequence,  if for each positive  number $\varepsilon\in {\bf No}$  there is an ordinal number $\alpha_0$ such that $|x_\alpha-x_{\alpha'}|<\varepsilon$, for all $\alpha_0\leq\alpha\leq\alpha'<\Omega$. 

{\bf Definition 22.}
Two $\Omega$-fundamental sequences $(x_\alpha)_{0<\alpha<\Omega}$ and $(y_\alpha)_{0<\alpha<\Omega}$ in ${\bf No}$ are $\Omega$-{\it equivalent}, denoted by as  $(x_\alpha)_{0<\alpha<\Omega}\sim(y_{\alpha})_{0<\alpha<\Omega}$,  if for each  positive number $\varepsilon\in {\bf No}$  there are ordinal numbers $\alpha_0$ and $\alpha'_0$ such that $|x_\alpha-y_{\alpha'}|<\varepsilon$, for all $\alpha_0\leq\alpha<\Omega$ and all $\alpha'_0\leq\alpha'<\Omega$. 

One can see that if $\lim\limits_{0<\alpha<\Omega}x_\alpha=a$ and $\lim\limits_{0<\alpha<\Omega}y_\alpha=b$, then $\lim\limits_{0<\alpha<\Omega}(x_\alpha+y_\alpha)=a+b$, $\lim\limits_{0<\alpha<\Omega}x_\alpha\cdot y_\alpha=a\cdot b$ and $\lim\limits_{0<\alpha<\Omega}\frac{x_\alpha}{y_\alpha}=\frac{a}{b}$ (latter when $b\not= 0)$.

Recall here the definition of a Dedekind section $(L,R)$, induced by an $\Omega$-fundamental sequence $(x_\alpha)_{0<\alpha<\Omega}$. We put $L$ as  the subset of ${\bf No}$ of all $l\in {\bf No}$ such that there exists an $\alpha_0$ and inequalities $ l<x_\alpha$ for all $\alpha_0\leq\alpha<\Omega$ and put $R= {\bf No}\setminus L$.

By the above experience, we can denote ${\bf No}$ as $O_\Omega=R_\Omega={\mathbb P}_\Omega$ and consider only $\Omega$-sequences $(x_\alpha)_{0\leq\alpha<\Omega}$ in ${\mathbb P}_\Omega$ and try to conceive what ${\tilde {\mathbb P}}_\Omega$ could be.

If $\Omega$-fundamental sequences $(x_\alpha)_{0<\alpha<\Omega}$ in ${\mathbb P}_\Omega$ is non-convergent and $(L,R)$ is a Dedekind section in ${\mathbb P}_\Omega$, induced by $(x_\alpha)_{0<\alpha<\Omega}$, then we cannot consider a mathematical object like ${\mathbb R}_\Omega$ as a completion of ${\mathbb P}_\Omega$ in the given $\Omega$-topology one of its elements is $(L,R)$, because Dedekind sections in ${\mathbb P}_\Omega$, induced by all $\Omega$-fundamental sequences $(x_\alpha)_{0<\alpha<\Omega}$ as well as $(L,R)$ itself, are illegal in a von Neumann-Bernays-G\"{o}del-type set theory $NBG$ and thus $(L,R)\in{\mathbb R}_\Omega$ is not well-formed formula in $NBG$
since these sections are proper classes and thus they  are not elements of any sets or classes. In particular, there is no object $\{(L,R)\}$ in $NBG$ or what is the same $\exists X\neg\{X\}$. On the other hand,  $\exists X(\forall x)({\bf M}(x))(x\in X)$ is a true formular.

But, nevertheless, there are operations in $NBG$ like $X\cup Y$, $X\cap Y$, $X\setminus Y$ etc., which are legal in $NBG$ because they are defined only by elements of $X$ and $Y$. Thus we can define a sum,  product and  quotient of two $\Omega$-fundamental sequences $(x_\alpha)_{0<\alpha<\Omega}$ and $(y_\alpha)_{0<\alpha<\Omega}$ in ${\mathbb P}_\Omega$ as $(x_\alpha+y_\alpha)_{0<\alpha<\Omega}$, $(x_\alpha\cdot y_\alpha)_{0<\alpha<\Omega}$ and $(\frac {x_\alpha}{y_\alpha})_{0<\alpha<\Omega}$ (latter when $(y_\alpha)_{0<\alpha<\Omega}$ is not a $\Omega$-infinite small sequence), respectively, which are also $\Omega$-fundamental. Moreover, these proper classes (i.e., $\Omega$-sequences themselves and induced by them Dedekind sections in ${\mathbb P}_\Omega$, or even lower classes $L$ as well as  upper classes $R$ of sections $(L,R)$) are gaps of the first kind in the sense of Conway $\cite{l3}$, p. 37-38 and thus we have defined the operations $+$, $\cdot$ and $/$ on these gaps, which do not depend on a choice of equivalent classes of $\Omega$-equivalent $\Omega$-sequences, and thus these gaps can be understood somehow or other as new  numbers in spite of the fact that they are inhibited in $NBG$ and they together with ${\mathbb P}_\Omega$, whose elements (numbers) can be identified with Dedekind sections in ${\mathbb P}_\Omega$ with extreme numbers, form a linear ordered Field. One can see that if a Dedekind section $(L,R)$, induced by $\Omega$-fundamental $\Omega$-sequence $(x_\alpha)_{0<\alpha<\Omega}$, has extremal elements, i.e., the greatest number in $L$ or the smallest number in $R$, then these Dedekind sections define Conway numbers which are these extremal elements and operations on them (sums,  products and  quotients) of $\Omega$-fundamental $\Omega$-sequences coincide with corresponding operations of Conway numbers. Thus such sections are forms of all Conway numbers and  sections which, induced by an $\Omega$-fundamental $\Omega$-sequence $(x_\alpha)_{0<\alpha<\Omega}$, have no extremal elements are forms of new numbers which are out of ${\bf No}$. Note that the collection of all these gaps (numbers) is even improper class (which together with ${\mathbb P}_\Omega$ should be ${\mathbb R}_\Omega$) and thus beyond the border of the University of all sets in $NBG$. Notice also that the gaps of the second kind are also numbers but we cannot define operations by means of elements or subclasses of ${\mathbb P}_\Omega$. But it is possible by the following considerations.

One can consider more wide formal set theories which allows proper classes be elements of super classes (see \cite{l133}, p. 119-153), then ${\mathbb R}_\Omega$ is a legal object in such theories and it is a topological field (in ${\Omega}$-topology) and we can extend it to ${\mathbb R}_\zeta$, for $\zeta=\omega^{\omega^{\Omega+1}}$, which contains all gaps of ${\bf No}$ including the gaps of the second kind in the sense of Conway. Thus, considering the above cases, the formula $\neg\exists X(\forall x)(\in X)$ is true, in general.

If we recall Cantor's definition of the first derivative  set $P'$ of a fixed P as the set  of all its limiting points, then further derivatives of this set may in some cases either end in an empty set, or become a perfect set whose derivative set coincides with it. In the first case, the sequence of derived sets determines the ordinal number $\alpha$ when the $(\alpha+1)$-derivative $P^{\alpha+1}$ is equal to the empty set.

In our case ${\bf No}={\mathbb P}_\Omega$ with $\Omega$-topology we can construct an example of a proper class $P=\bigcup\limits_{0\leq\alpha<\Omega}P_\alpha\subset[0,1]\subset{\bf No}$ such that $P_\alpha^{(\alpha)}={\frac{1}{\alpha+1}}$, $0<\alpha<\Omega$,  and $P^{(0)}=P$ and $P^{(\Omega)}=\{0\}$ with the following order $P^{(\alpha)}\supset P^{(\alpha+1)}\Leftrightarrow \alpha<\alpha+1$, $0\leq\alpha\leq\Omega$.
Thus proper classes $P^{(\alpha)}\subset{\bf No}$, $0\leq\alpha<\Omega$, and the set $\{0\}$ determine all ordinal numbers including $\Omega$.

Here is a construction of this example $P$.

Put $Q_0=\{0\}$; put $P_0=Q_0$; clearly $P_0'=\{0\}'=\emptyset$ thus $P_0$ presents the ordinal $0$.

Put $Q_1=(Q_0+\{0,1,\frac{1}{2},\frac{1}{3},...,\frac{1}{\alpha},...\}_{0<\alpha<\Omega})\cap[0,1]$; put $P_1=Q_1\cap[\frac{1}{2},1]$; clearly $(P_1')'=\{\frac{1}{2}\}'=\emptyset$ thus $P_1$ presents the orninal $1$.

Put $Q_2=(Q_1+\{0,1,\frac{1}{2},\frac{1}{3},...,\frac{1}{\alpha},...\}_{0<\alpha<\Omega})\cap[0,\frac{1}{2}]$; put $P_2=Q_2\cap[\frac{1}{3},\frac{1}{2}]$; clearly $(P_2'')'=\{\frac{1}{3}\}'=\emptyset$ thus $P_2$ presents the ordinal $2$.

Put $Q_3=(Q_2+\{0,1,\frac{1}{2},\frac{1}{3},...,\frac{1}{\alpha},...\}_{0<\alpha<\Omega})\cap[0,\frac{1}{3}]$; put $P_3=Q_3\cap[\frac{1}{4},\frac{1}{3}]$; clearly $(P_3''')'=\{\frac{1}{4}\}'=\emptyset$ thus $P_3$ presents the ordinal $3$.

.............................

Put $Q_n=(Q_{n-1}+\{0,1,\frac{1}{2},\frac{1}{3},...,\frac{1}{\alpha},...\}_{0<\alpha<\Omega})\cap[0,\frac{1}n]$; put $P_n=Q_n\cap[\frac{1}{n+1},\frac{1}{n}]$; clearly $(P_n^{(n)})'=\{\frac{1}{n}\}'=\emptyset$ thus $P_n$ presents the orninal $n$.

.............................

Put $Q_\omega=(\bigcup\limits_{0\leq n<\omega}Q_n)\cap\bigcap\limits_{0\leq n<\omega}[0,\frac{1}{n}]$; put $P_\omega=Q_\omega\cap\bigcap\limits_{0\leq n<\omega}[\frac{1}{\omega},\frac{1}{n}]$; clearly $(P_\omega^{\omega})'=\{\frac{1}{\omega}\}'=\emptyset$ thus $P_\omega$ presents the orninal $\omega$.

Put $Q_{\omega+1}=Q_\omega+\{0,1,\frac{1}{2},\frac{1}{3},...,\frac{1}{\alpha},...\}_{0<\alpha<\Omega})\cap[0,\frac{1}{\omega}]$; put $P_{\omega+1}=Q_{\omega}\cap[\frac{1}{\omega+1},\frac{1}{\omega}]$; clearly $(P_{\omega+1}^{(\omega+1)})'=\{\frac{1}{\omega+1}\}'=\emptyset$ thus $P_{\omega+1}$ presents the ordinal $\omega+1$.

.............................

Put $Q_\beta=(\bigcup\limits_{0\leq \alpha<\beta}Q_\alpha)\cap\bigcap\limits_{0< \alpha<\beta}[0,\frac{1}{\alpha}]$; put $P_\beta=Q_\beta\cap\bigcap\limits_{0<\alpha<\beta}[\frac{1}{\beta},\frac{1}{\alpha}]$ (Here $\beta$ is a limit ordinal.); clearly $(P_\beta^{\beta})'=\{\frac{1}{\beta}\}'=\emptyset$ thus $P_\beta$ presents the ordinal $\beta$.

.............................

Put $Q_{\gamma}=Q_{\gamma-1}+\{0,1,\frac{1}{2},\frac{1}{3},...,\frac{1}{\alpha},...\}_{0<\alpha<\Omega})\cap[0,\frac{1}{\gamma}]$. (Here $\gamma$ is not a limit ordinal.); put $P_\gamma=Q_{\gamma}\cap[\frac{1}{\gamma+1},\frac{1}{\gamma}]$; clearly $(P_{\gamma1}^{(\gamma1)})'=\{\frac{1}{\\frac{1}{\gamma}}\}'=\emptyset$ thus $P_{\gamma}$ presents the ordinal $\gamma$.

.............................

$P=(\bigcup\limits_{0\leq \alpha<\Omega}P_\alpha)$; clearly $(P^{(\Omega)})'=\{0\}'=\emptyset$ thus $P$ presents the ordinal $\Omega$.

Notice, first, that by  a {\it sum} of of two classes $S$ and $T$ we understand $S+T\stackrel{def}{=}\{s+t\,\,|\,\,s\in S\,\&\,t\in T\}$ and, second, $U=P+Q_1$ has the following property: $U^{(\Omega)}=Q_1$ and thus $({U^{(\Omega)}})''=(P_1)''=\{0\}'=\emptyset$. Consequently the class $U$ and all its limit points determine an ordinal number $\Omega+1$ and such similar process has no boundary.

We are going to explain this state of affairs.

In spite of the fact that the formal logical system $NBG$ is convenient in many cases (e.g., Conway theory of numbers and games, as well as  Category Theory and many other mathematical theories and constructions) by the finiteness of its axioms, by existence of universal objects (proper classes) identified with \grqq properties\grqq\, or \grqq singular proposition functions\grqq, etc., there is a lack in it which, by thought of P. J. Cohen, the \grqq theory $NBG$ is a less intuitive system than $ZF$\grqq (\cite{l222}, Cpt. 3, $\S 5$). Indeed, such formal systems are very formal and they were inputed to avoid so-called set-theoretic paradoxes of Russell's type. Often these \grqq tricks\grqq\, (proofs by contradiction via paradoxes) were out of contents and concepts of mathematical notions and definitions like \grqq reflexive sets\grqq, \grqq universal objects\grqq, \grqq premises of false propositions\grqq, \grqq undetectable  objects\grqq, etc.

The only correct and meaningful conclusion from Russell's paradox is the following: \grqq Any family of non-reflexive elements is always non-reflexive and is not contained in the original family\grqq. As a consequence of it is the following: \grqq There is no universal families like all sets, all ordinals, all groups, all rings, all fields, etc.\grqq. 

An attempt to restrict some families (such as proper classes) to be members of other families is formal and contradicts with the first statement. For example, the first axiom of the existence of classes in the NBG theory is the following: 
\begin{equation}
\label{f0770}
\exists X\forall u\forall v(<u,v>\in X)\Leftrightarrow u\in v,
\end{equation}
where $<u,v>=\{\{u\},\{u,v\}\}$ is a Kuratowski ordered pair, exactly repeats Russell's paradox because $<X,\{X\}>=\{\{X\},\{X,\{X\}\}\}$ is a pair $<u,v>$, where  $u=X$, $v=\{X\}$ such that $<u,v>\notin X$. But Axiom $(\ref{f0770})$ requires a relation $<u,v>\in X$ what is evidently false when $NBG$ has a Regular axiom.

The authors of the monograph \cite{l98}, p. 59, cite the following formal theorem of first-order logic: 
\begin{equation}
\label{f100}
\neg\exists x\forall y[yEx\Leftrightarrow\neg(yEy)] 
\end{equation}
without any assumptions about the notion of set, other than that the displayed predicate determines a unique set, independent of what $E$ happens to mean. The authors of the mention monograph mean of course the substitution rule: if under the assumption that such a set $x$ exists, i.e. the formula $(\exists x\forall y[yEx\Leftrightarrow\neg(yEy)])$  is true, where $y$ is a variable and $x$ is a constant of above propositional form $(yEx\Leftrightarrow\neg(yEy))$, then we substitute $x$ instead of $y$ into it   and obtain the paradoxical formula \begin{equation}
\label{f101}
xEx\Leftrightarrow\neg(xEx) 
\end{equation}
which implies, by formal methods, that there is no such $x$ or in $NBG$ theory a premise that $x$ is a set, i.e., ${\bf M}(x)$, by the logical rule $B\Rightarrow A\& \neg A\Rightarrow \neg B$, obtain that $\neg {\bf M}(x)$, i.e., $x$ is a proper class. Due to a contradiction it is easily to deduce (see this conclusion in detail in \cite{l144}, Chapter III, \S 1, p. 173) that there is no such set, i.e.,  $x$ does not exist. In particular, if the Russell class $R\stackrel{def}{=}\{z\,|\,z\notin z\}$ is understood as $x$, and the membership relation $\in$ is understood as $E$, then the assumption $(\exists R\forall y[yER\Leftrightarrow\neg(yEy)])$ implies $(R\in R\Leftrightarrow\neg(R\in R))$, and thus the Russell's paradox is obtained. 

(Note that if the same Russell class $R$ is understood as $x$, and the non-affiliation relation $\notin$ is understood as $E$, then the assumption  $(\exists R\forall y[y\notin R\Leftrightarrow(y\in y)])$ implies formally the same Russell's paradox  $R\notin R\Leftrightarrow\neg(R\notin R))=R\in R$ is obtained, although, as we will see below, there is a qualitative difference between them.)

A formal axiomatic method is very important in  modern mathematics especially in its foundation.
\grqq It is also important that a proof satisfying our strictness condition should remain valid for any
interpretation of the logistics system, so that we ultimately get savings due to the fact that various things are proved
by a single reasoning. The amount of savings is determined by the fact that there is no need to repeat an indefinite number of times proofs
that coincide in form but differ in content, so how can they be held all at once once and for all\grqq\, (see \cite{l23}).  But often these advantages can also become a disadvantage due to its {\it formality}.
We can illustrate this with the following specific examples.

Let $E$ in (\ref{f100}) be a predicate $<_{Card}$, i.e., an inequality of cardinality (power) of sets: $y<_{Card} x\stackrel{def}{\Leftrightarrow}Card(y)<Card(x)$, where $Card (x)$ is a cardinal number corresponding to the  cardinality (power) of a set $x$.

Then (\ref{f100}) turns to the following formula
\begin{equation}
\label{f102}
\neg\exists x\forall y[y<_{Card}x\Leftrightarrow\neg(y<_{Card}y)] 
\end{equation}
what is the same as
\begin{equation}
\label{f103}
\neg\exists x\forall y[y<_{Card}x\Leftrightarrow(Card(y)\geq Card(y))].
\end{equation}
Suppose now inverse that  the formula $\exists x\forall y[y<_{Card}x\Leftrightarrow(Card(y)\geq Card(y))]$ is true. Then the propositional form $[y<_{Card}x\Leftrightarrow(Card(y)\geq Card(y))]$ has a region the class of all sets since the relation $(Card(y)\geq Card(y))$ fulfils for every set $y$, in particular, for a supposed set $x$, and thus, by the substitution rule, we can substitute $x$ instead of $y$ and obtain a paradoxical formula $x<_{Card}x\Leftrightarrow x\geq_{Card} x$ what of $NBG$ as above proves that $x$ is not a set but a proper class. Moreover, in $NBG$ it is a true conclusion since for any set $y$ and the universal proper class $x$ of all sets, which exists by axioms in $NBG$, we obtain a true formula $Card(y)<Card(x)$.

But everything is not so simple in the case when we consider a predicate $>_{Card}$, where $y>_{Card} x\stackrel{def}{\Leftrightarrow}Card(y)>Card(x)$, and suppose that such $x$ exists. Then the propositional form $[y>_{Card}x\Leftrightarrow(Card(y)\leq Card(y))]$ has a region the class of all sets since the relation $(Card(y)\leq Card(y))$ fulfills for every set $y$. After this, the similar paradoxical formula $x>_{Card}x\Leftrightarrow x\leq_{Card} x$ implies  that in $NBG$ $x$ is not only a set but it is not also a proper class because for supposed $x$ the true formula is the following: $\forall y (Card(y)>Card(x))$. Thus, such $x$ does not exist in $NBG$ even for the theory with individuals for which such notion exists, i.e., an individual ${\bf Ur}(x)=\neg{\bf Cls}(x)$ is not a class  at all.

If we consider $E$ as $\leq_{Card}$ or $\geq_{Card}$, then the corresponding propositional forms $[y\leq_{Card}x\Leftrightarrow(Card(y)>Card(y))]$ and $[y\geq_{Card}x\Leftrightarrow(Card(y)<Card(y))]$ with a variable $y$ and a constant $x$ have  empty regions. Thus, in these cases the substitution rule does not work at all.

All these arguments can be repeated in the cases when $E=(\sim)$ or $E=(\neg\sim)$, in particular $E=(=)$ or $E=(\not=)$, as well as $E=(\subset)$ or $E=(\neg\subset)$, $E=(\supset)$ or $E=(\neg\supset)$ .

Now we will consider another aspect of the unsatisfactoriness of the formal method.

Let us consider so-called a relative Russell's paradox for any fixed set or class $X$. Denote $R_X\stackrel{def}{=}\{x\in X\,|\, x\notin x\}$ and the following formal theorem of first-order logic: $R_x\notin X$, which can be proved by contradiction with a help of Russell's paradox. Indeed, suppose that $R_X\in X$. Since $x\in R_X\Leftrightarrow x\in X\,\&\,x\notin x$ we substitute $R_X$ instead of $x$ in above formula and obtain $R_X\in R_X\Leftrightarrow R_X\notin R_X$. Contradiction. Thus our assumption is false.

It is easy to see the redundancy of the  above  formal proof that $R_X\notin X$ if we consider a trivial case when $X=\emptyset$. In this case the above theorem is the following: $R_\emptyset\notin\emptyset$. Indeed, truth of this theorem is evident by the definition of the empty set. Nevertheless, following the ideology of formal systems and formal methods to give one proof for all possible cases in one step, we repeat it with a help of Russell's paradox. Suppose the opposite $R_\emptyset\in\emptyset$. Then, by definitions $R_\emptyset\in\emptyset\Leftrightarrow R_\emptyset\notin\emptyset$, or more simple, if we suppose that $\emptyset\in\emptyset$, then, by Russell's paradox, obtain $\emptyset\in\emptyset\Leftrightarrow \emptyset\notin\emptyset$. Contradiction. Thus $\emptyset\notin\emptyset$. Note here that {\it all proofs with a help of Russell's paradox have the same intension}, e.g., when $NBG$ or $ZF$ include Regular axiom and $X$ is always founded and thus there is no formulas $X\in X$ in $NBG$ or $ZF$ and thus $R_X\in R_X\Rightarrow R_X\notin R_X$ is a true implication; on the other hand, $R_X\notin R_X\Rightarrow R_X\in R_X$ is always false because of Regular axiom as well as the nature of reflective sets and at last in the initial Russell's paradox because of false premise of an existence of universal set (class or family).

As we noted above if the  Russell class $R$ is understood as $x$, and the non-affiliation relation $\notin$ is understood as $E$, then the assumption  $(\exists R\forall y[y\notin R\Leftrightarrow(y\in y)])$ implies formally the same Russell's paradox  $R\notin R\Leftrightarrow\neg(R\notin R))=R\in R$ is obtained, although, if we consider $NBG$ or $ZF$ formal theories with the Regular Axiom, then the propositional form $[y\notin R\Leftrightarrow(y\in y)]$ with a variable $y$ has an empty region and thus the formula $y\in y$ has no denotate and there is no possibility to substitute $x$ instead of $y$ in $y\in y$.

However, one can easily pass over Russell's paradox in the following way (see also \cite{l1333} and \cite{l1393}). $R_X\notin R_X$ because all elements of $R_X$ are non-reflexive sets. Thus if $R_X\in X$, then $R_X\cup\{R_X\}\subseteq X$ and $R_X\subset R_X\cup\{R_X\}$ what contradicts with maximality of $R_X$. Thus $R_X\notin X$. Notice that $\{R_X\}$ always exists either it is a set or a proper class. Only fictitious formal Russell's paradox deducts that proper class is not an element. But it is very artificial to identify \grqq property\grqq\, with a \grqq proper class\grqq. The model of $NBG$ as the existence of an inaccessible cardinal number contradicts with it since there is a one-element set with this inaccessible cardinal number.

Moreover, this circumvention of the Russell's paradox makes it possible to prove propositions, which are impossible by  formal methods. Indeed, denote by $S=\{z\,|\, z\in z\}$ the collection of all reflexive sets in $NBG$. Then in $NBG$ the following Proposition is true: $(\neg\exists x\forall y[y\in x\Leftrightarrow(y\in y)])$ and hence $x$ is not a set but a proper class. The proof by the above formal method is useless because, by opposite and subsistence rule, we obtain a tautology $x\in x\Leftrightarrow x\in x$.

But by an assumption that such set $S$ exists we consider a reflexive set $S'$ such that $S'\setminus\{S'\}=S$, i.e., $S'=S\cup\{ S'\}$ (see, e.g., in \cite{l1333}), and prove that $S'\notin S$. Really, for each element $z\in S$ the set $z\setminus\{z\}$ does not have $z$ as its element. But $z$ is an element of $S'$. Consequently, $S\subset S\cup\{S'\}$, i.e., $S$ is a proper subset of $S\cup\{S'\}$ what contradicts with maximality of $S$. Thus $S$ is a proper class and there is no a singleton $\{S\}$. Note of course that we mentioned here a theory of reflexive sets (see, e.g., \cite{l1333}) for which an equality $z=S'$ of two reflexive sets $z$ and $S'$ implies the identity of $z\setminus\{z\}=S'\setminus\{S'\}$. If it is not so (see, e.g., \cite{l93}), then one must give another proof that $S$ is a proper class.

So, choosing the axiom of the existence of a singleton set $\{X\}$ and getting the result of the non-existence of a universal class ${\bf V}$ for sets is preferable from the point of view of the state of affairs or subsistence than postulating the axiom $\forall X {\bf P}r(X)\neg\{X\}$ of denying the existence of a singleton set $\{X\}$ for one's own classes and getting a universal class ${\bf V}$ with questionable \grqq geography\grqq. Conway expressed his attitude to this very problem by the following deep intuition:

\grqq This appendix is in fact a cry for a Mathematicians' Liberation Movement!

$(i)$ Objects may be created from earlier objects in any reasonably constructive fashion.

$(ii)$ Equality among the created objects can be desired equivalence relation\grqq. (See $\cite{l3}$, p. 66.)


\begin{thebibliography}{13}
  \bibitem{l93} P. Aczel, {\it Non-Well-Founded Sets}. CSLI Lecture Notes Number 14. Stanford: CSLI Publication, 1988. 

\bibitem{l98} J. Barwise, L. Moss, {\it Vicious Circles. On the Mathematics of Non-Wellfounded Phenomena}. CSLI Lecture Notes Number 60, Stanford, California, 1996.

\bibitem{l222} P. J. Cohen, {\it Set Theory and the Continuum Hypothesis}. W. A. Benjamin, INC. New York, Amsterdam, 1966.


\bibitem{l3}
John H. Conway,    
On Numbers and Games,
 London Mathematical Society Monographs, No. 6, Academic Press, London-New-San Francisco, 1976.
  \bibitem{l23}
  A. Church, Introduction to Mathematical Logic, Pronceton, New Jersey, 1956.
\bibitem{l4}
 R. Dedekind, 
Was sind und was wollen die Zahlen?,
(Appeared originally in 1888) 6th ad.: Braunschweig,   1930.
\bibitem{l144}
A. A. Fraenkel, Y. Bar-Hillel,
{\it A. A. Fraenkel, Y. Bar-Hillel}, North-Holland Publication Company, Amsterdam, 1956. 
\bibitem{l133}
A. A. Fraenkel, Y. Bar-Hillel, A. Levy,
 Foundation of set theory,
Elsevier,  Amsterdam  London New York Oxford Paris Shannon Tokyo, 1973.
\bibitem{l15}
E. Jacobstahl, \"{U}ber der Aufbau der infinite Arithmetic, Math. Ann., 66 (1909), 145-194.
 

\bibitem{l1333}
 Ju. T. Lisica {\it Skand theory and its applications. (A new look at non-well-founded sets)} https://archive.org/details/arxiv-1207.2985/page/n25/mode/2up
2012, pp. 1-61.
\bibitem{l1393}
 Ju. T. Lisica {\it On one of Specker's theorems}
Topology and its Applications
Volume 329, 15 April 2023, pp. 108383-108414. 

\bibitem{l12}
D. Mirimanoff,
  Les antinomies de Russell et de Burali-Forti et le probl\`{e}me fondumentale de la th\'{e}orie des ensembles,
 L'Enseignement Math\'{e}matique, 19,  (1917), pp. 37-52.
\bibitem{l090}
 J. von Neuman, 
  Zur Einf\"{u}hrung der infiniten Zahlen,
Acta Litt. Ac. Sci. Hung. Fran. Joseph., {\bf 1} (1923), pp. 199-208.



\end{thebibliography}
 \end{document}